%% file: Brinkman_paper.tex
\documentclass[review]{elsarticle}
\makeatletter
\def\ps@pprintTitle{%
 \let\@oddhead\@empty
 \let\@evenhead\@empty
 \def\@oddfoot{\centerline{\thepage}}%
 \let\@evenfoot\@oddfoot}
\makeatother

\pdfoutput=1
\usepackage{graphicx}
\usepackage{amsfonts, amsmath, amssymb}
\usepackage{booktabs,siunitx}
\usepackage[svgnames,table]{xcolor}
\usepackage[tableposition=above]{caption}
\usepackage{pifont}
\usepackage{bm}
\usepackage{cancel}
\usepackage{caption}
\usepackage{color}
\usepackage{enumerate}
\usepackage{float}
\usepackage{hyperref}
\usepackage[english]{babel}
\usepackage[margin=1in]{geometry}
\usepackage{mathtools}
\usepackage{multirow}
\usepackage{setspace}
\usepackage{subfigure}
\usepackage{lineno}

\renewcommand \d [2]{\frac{{\rm d} #1}{{\rm d} #2}}
\newcommand \D [2]{\frac{\partial #1}{\partial #2}}

\renewcommand{\vec}[1]{\bm{\mathrm{#1}}}
\newcommand{\V}[1]{\bm{\mathrm{#1}}}

\def \div{\nabla \cdot \mbox{}}
\def \grad{\nabla}

\def \x{\vec{x}}

\def \qdi{q^{\text{d}}_i}
\def \n{\vec{n}}
\def \r{\vec{r}}
\def \u{\vec{u}}

\def \e{\vec{e}}

\def \nprop{n_{\rm prop}}

\def \vbeta{\V{\beta}}

\def \fb{f_\text{b}}

\def \Omegas{\Omega_{\text{s}}}
\def \Omegaf{\Omega_{\text{f}}}

\def \X{\vec{X}}

\def \e{\vec{e}}

\def \f{\vec{f}}

\def \half{\frac{1}{2}}
\def \3half{\frac{3}{2}}
\def \5half{\frac{5}{2}}

\def \n{\vec{n}}

\def \ncells{n_{\text{smear}}}

\def \nsmear{n_{\text{smear}}}

\def \qexact{q_\text{exact}}
\def \qinexact{\widetilde{q}}

\def \u{\vec{u}}
\def \ub{\u_{\text{b}}}

\def \x{\vec{x}}

\def \div{\V{\nabla} \cdot \mbox{}}
\def \grad{\V{\nabla}}

\def \Dr{{\mathrm d}r}

\def \dS{\,\mathrm{dS}}
\def \dV{\,\mathrm{dV}}

\def \deltaf{\delta_{\partial \Omegaf}}

\newcommand{\upperRomannumeral}[1]{\uppercase\expandafter{\romannumeral#1}}

\newcommand{\REVIEW}[1]{#1}
\newcommand{\REVIEWRED}[1]{{}}

\begin{document}
\let\today\relax

\begin{frontmatter}
	
\title{Handling Neumann and Robin boundary conditions in a fictitious domain volume penalization framework}
\author[SDSU]{Ramakrishnan Thirumalaisamy}
\author[Northwestern]{Neelesh A. Patankar}
\author[SDSU]{Amneet Pal Singh Bhalla\corref{mycorrespondingauthor}}
\ead{asbhalla@sdsu.edu}

\address[SDSU]{Department of Mechanical Engineering, San Diego State University, San Diego, CA}
\address[Northwestern]{Department of Mechanical Engineering, Northwestern University, Evanston, IL}
\cortext[mycorrespondingauthor]{Corresponding author}

\begin{abstract}
\input{Abstract}
\end{abstract}


\end{frontmatter}

\section{Introduction} \label{sec_intro}
\input{Introduction}

\section{Mathematical formulation}\label{sec_math_eqs}
\input{Mathematical_formulation}

\section{Software} \label{sec_software}

The flux-based volume penalization algorithms  described here are implemented within the IBAMR library~\cite{IBAMR-web-page}, which is an open-source C++ simulation software focused on immersed boundary and volume penalization methods with adaptive mesh refinement. The code and test cases presented in Sec.~\ref{sec_results_and_discussion} are publicly available at \url{https://github.com/IBAMR/IBAMR}. IBAMR relies on SAMRAI~\cite{HornungKohn02, samrai-web-page} for Cartesian grid management and the AMR framework. Linear and nonlinear solver support in IBAMR is provided by the PETSc library~\cite{petsc-efficient, petsc-user-ref, petsc-web-page}. All of the example cases in the present work made use of distributed-memory parallelism using the Message Passing Interface (MPI) library.

\section{Results and discussion}
\label{sec_results_and_discussion}
\input{Results_and_discussion}

 \section{Conclusions}
 \label{sec_conclusions}
 \input{Conclusions}


\section*{Acknowledgements}
R.T and A.P.S.B~acknowledge support from NSF award OAC 1931368.  R.T acknowledges support from San Diego State University Graduate Fellowship award. N.A.P~acknowledges support from NSF award OAC 1931372. Computational resources provided by Fermi high performance computing cluster at San Diego State are also acknowledged.


\appendix
\renewcommand\thesection{\Alph{section}}
\input{Appendix}

\newpage
\section*{Bibliography}
\begin{flushleft}
 \bibliography{Brinkman_paper_biblography}
\end{flushleft}

\end{document}

%% file: Abstract.tex
Sakurai et al. (J Comput Phys, 2019) presented a flux-based volume penalization (VP) approach for imposing inhomogeneous Neumann boundary conditions on embedded interfaces. The flux-based VP method modifies the diffusion coefficient of the original elliptic (Poisson) equation and uses a flux-forcing function as a source term in the equation to impose the Neumann boundary conditions. As such, the flux-based VP method can be easily incorporated into existing fictitious domain codes. Sakurai et al. relied on an analytical construction of flux-forcing functions, which limits the practicality of the approach. Because of the analytical approach taken in the prior work, only (spatially) constant flux values on simple interfaces were considered. In this paper, we present a numerical technique for constructing flux-forcing functions for arbitrarily complex  boundaries. The imposed flux values are also allowed to vary spatially in our approach.  Furthermore, the flux-based VP method is extended to include (spatially varying) Robin boundary conditions, which makes the flux-based VP method even more general. The numerical construction of the flux-forcing functions relies only on a signed distance function that describes the distance of a grid point from the interface and can be constructed for any irregular boundary. We consider several two- and three-dimensional test examples to access the spatial accuracy of the numerical solutions.  The method is also used to simulate flux-driven thermal convection in a concentric annular domain. We formally derive the flux-based volume penalized Poisson equation satisfying Neumann/Robin boundary condition in strong form; such a derivation was not presented in Sakurai et al., where the equation first appeared for the Neumann problem. The derivation reveals that the flux-based VP approach relies on a surface delta function to impose inhomogeneous Neumann/Robin boundary conditions. However, explicit construction of the delta function is not necessary for the flux-based VP method, which makes it different from other diffuse domain equations presented in the literature.    

%% file: Introduction.tex
Partial differential equations (PDEs) in complex domains describe many natural and engineering processes. Examples include heat and mass transfer across melting/solidifying fronts, aquatic locomotion, cellular phenomena like cellular blebbing and cell crawling, flow in internal combustion engines or left ventricular assist devices, energy harvesting using wind turbines and wave energy converters, etc. In order to obtain meaningful solutions to PDEs, appropriate boundary conditions are required on the domain boundaries. Traditionally, body-fitted grid approaches, in which a complex domain is triangulated using sophisticated grid generation software, have been employed to solve PDEs numerically. Although body-fitted grid approaches allow imposing various types of boundary conditions accurately, they pose a serious challenge when the solution domain changes its topology over time. Issues like constant remeshing of the computational domain, the high aspect ratio of the elements, etc., limit the feasibility of body-fitted grid methods for modeling challenging moving domain problems.    

To overcome the limitations of the body-fitted grid methods, fictitious domain (FD) methods have been proposed. In fictitious domain methods, an irregular region of interest is embedded into a larger, simpler computational domain and the original PDE is reformulated on the entire domain. FD methods typically employ regular Cartesian grids to mesh the computational domain. This allows simpler discretization of PDEs and fast linear solvers to solve the discrete system of equations.  Since the regular grid no longer adheres to the irregular interface, incorporating original boundary conditions in the reformulated equation is not straightforward. Nevertheless, several techniques to incorporate Dirichlet boundary conditions have been proposed for various variants of the FD method. Dirichlet boundary conditions are particularly relevant for modeling fluid-structure interaction (FSI) problems, where velocity matching condition on the fluid-structure interface is required. Fictitious domain methods such as the immersed boundary (IB) method~\cite{Peskin02} and the volume penalization (VP) method~\cite{Angot1999} have been successfully used to model several FSI problems, including wave energy converters~\cite{Bergmann2015,Dafnakis2020,Khedkar2020}, water entry/exit problems~\cite{BhallaBP2019},  fish swimming~\cite{Bergmann2011,Bhalla13}, esophageal transport~\cite{Kou2015}, cardiovascular flows~\cite{Griffith2012}, etc.  The IB method was introduced by Peskin to model flow in a human heart~\cite{Peskin1972} and is a two-grid approach to FSI modeling: Lagrangian mesh for describing the moving structure and an Eulerian grid for describing the fluid flow. In contrast, the VP method introduced by Angot et al.~\cite{Angot1999} is a single grid approach in which all quantities related to fluid and structure are described on the Eulerian grid.  The moving structure in the VP method is typically tracked using an indicator function. Since all quantities are described on a single grid, parallel implementation of VP methods on distributed memory systems is relatively easier compared to the two-grid IB methods.  



The original VP method introduced by Angot et al.~\cite{Angot1999} considered only Dirichlet boundary conditions. Later the VP method was generalized to Neumann and Robin boundary conditions by Rami\`ere, Angot, and Belliard~\cite{Ramiere2007a}. The authors in~\cite{Ramiere2007a} implemented their VP technique within a finite element framework. In their formulation, inhomogeneous Neumann and Robin boundary conditions were incorporated by introducing a surface delta function in the reformulated equation; the singular delta function was regularized in the numerical implementation.  Recently, in Kadoch et al.~\cite{Kadoch2012}, a volume penalization method for imposing homogeneous Neumann boundary conditions was presented. The authors in~\cite{Kadoch2012} implemented their method within a pseudo-spectral code and used it to simulate moving domain problems involving chemical mixers. Since homogenous (Neumann) boundary conditions were considered in Kadoch et al.~\cite{Kadoch2012}, the need for a surface delta function kernel was circumvented. More recently, Sakurai and co-workers~\cite{Sakurai2019} introduced the so-called flux-based VP method, which  extends Kadoch et al.'s approach to imposing inhomogeneous Neumann boundary conditions. The flux-based VP approach uses a flux-forcing function to impose the inhomogeneous Neumann boundary conditions on the interface. Sakurai et al. used second-order central finite differences and interpolation to implement the flux-based VP method and solved several one- and two-dimensional Poisson problems to assess the spatial convergence rate of the numerical solutions. Sakurai et al. considered simple interfaces in two-spatial dimensions (circles and rectangles) in their study, which allowed them to construct flux-forcing functions analytically. Moreover, the imposed flux values were considered spatially constant on the interface. The analytical construction of flux-forcing functions limits the feasibility of the flux-based VP method for practical applications. \REVIEW{Recently, Thirumalaisamy et al.~\cite{Thirumalaisamy2021} critiqued Sakurai et al. for some inconsistencies in  their results and conclusions, following which the authors of~\cite{Sakurai2019} published a corrigendum~\cite{Sakurai2021} to their original work. Similar to Sakurai et al.,  Thirumalaisamy et al. also relied on the analytical construction of flux-forcing functions for the flux-based VP method.}

One of the objectives of this work is to generalize the flux-based VP method to handle arbitrarily complex interfaces in two and three spatial dimensions. This is achieved through numerical construction of flux-forcing functions, as described in Sec.~\ref{sec_beta} of this paper. Moreover, the imposed flux values are allowed to vary spatially on the interface. The proposed numerical approach for constructing flux-forcing functions requires only a signed distance function that describes the distance of a grid point from the interface. The signed distance function can be constructed analytically for simple geometries, or through computational geometry techniques for complex interfaces~\cite{Baerentzen2005}. Another objective of this work is to extend the flux-based VP method to include (spatially varying) Robin boundary conditions. This allows imposing both types of boundary conditions (Neumann and Robin) through similar (numerical) flux-forcing functions.

Similar to Sakurai et al., we also discretize the volume penalized equations using second-order finite differences. Using the method of manufactured solution, we assess the accuracy of the proposed approach by solving two- and three-dimensional Poisson problems with constant and spatially varying Neumann/Robin boundary conditions. \REVIEW{We compare the performance of our approach using both continuous and discontinuous indicator functions in the test problems considered in Sec.~\ref{sec_results_and_discussion}. It is observed that the continuous indicator function performs better (in terms of order of accuracy and uniformity of convergence rate) for imposing the spatially constant Neumann/Robin boundary condition, whereas the discontinuous one performs better for the spatially varying Neumann/Robin problem.}

We also provide a formal derivation of the flux-based VP Poisson equation, which was not provided in Sakurai et al.~\cite{Sakurai2019}, where the equation first appeared for the Neumann problem. The derivation reveals that the flux-based volume penalization method also uses a surface delta function to impose inhomogeneous Neumann/Robin boundary conditions. Interestingly, explicit  construction of the delta function is not required in the flux-based approach, which is in contrast to the volume penalization approach of Rami\`ere et al.~\cite{Ramiere2007a}. We remark that on a formulation level the volume penalization approaches of  Rami\`ere et al., Kadoch et al., and Sakurai et al. (and the present work) are equivalent; minor differences in these works arise from the definition of the surface delta function. This insight is gained from Li et al.~\cite{Li2009} who derived phase field-based diffuse domain equations satisfying Dirichlet, Neumann,  and Robin boundary conditions. Li et al. used the method of matched asymptotic expansions to provide different diffuse domain approximations for the Neumann/Robin problem~\footnote{Different diffuse domain approximations for the Dirichlet problem are also provided in Li et al.~\cite{Li2009}.}; these approximations differ in the way how surface delta function is defined.  


Characteristic-based approaches to impose Neumann and Robin boundary conditions for the volume penalized PDEs have also been proposed in the literature; see, for example, Brown-Dymkoski et al.~\cite{Brown2014characteristic} and Hardy et al.~\cite{Hardy2019penalization} who used characteristic-based VP approach to model the energy transport equation satisfying Neumann and Robin boundary conditions in the context of compressible flows and low Mach formulation of compressible flows, respectively. The main limitation of the characteristic-based VP method is that it relies on having a time-derivative term in the PDE and as such cannot be applied to steady-state (i.e., having no temporal derivative term) elliptic equations.  \REVIEW{In addition to the volume penalization methods~\cite{Ramiere2007a,Kadoch2012,Sakurai2019,Brown2014characteristic,Hardy2019penalization,bensiali2015penalization,Kolomenskiy2015,Schneider2015}}, other fictitious domain techniques have also been proposed to impose flux boundary conditions on embedded interfaces. Notable ones include the flux-correction technique (FCT) of Ren et al.~\cite{Ren2013}, Wang et al.~\cite{Wang2016}, and Guo et al.~\cite{Guo2019} and the direct forcing method of Lou et al.~\cite{Nils2020}. FCT is a predictor-corrector scheme and is implemented using the Lagrangian-Eulerian machinery of the IB method. In the prediction step of FCT, an intermediate scalar field is computed on the Eulerian grid, which in general does not satisfy the flux boundary condition on the interface defined by the Lagrangian markers. Next, in the correction step, a Lagrangian forcing term is computed either implicitly~\cite{Wang2016,Guo2019} or explicitly~\cite{Ren2013} that corrects the intermediate scalar field to satisfy the Neumann boundary condition. In an essence, FCT is a time-splitting approach (similar to the characteristic-based VP approach), which requires having a time-derivative term in the scalar transport equation. Therefore, unlike the flux-based VP method, FCT cannot be used for time-independent elliptic equations. In the direct forcing method, the scalar field near the interface is reconstructed locally using second- or third-degree polynomials in order to satisfy the flux boundary condition. This is achieved by identifying ``forcing" points on the fictitious (solid) side of the interface, on which the reconstructed scalar field value is directly imposed. Direct forcing methods are also typically implemented as a predictor-corrector scheme, which avoids modifying the system of linear equations.  

In the following sections, we first describe the continuous form of the volume penalized equations and thereafter describe the numerical construction of the flux-forcing functions. Finally, various test cases are considered in two- and three-spatial dimensions to access the accuracy of the numerical solutions.


%% file: Mathematical_formulation.tex
\subsection{The Neumann problem}

\begin{figure}[]
\centering
   \includegraphics[scale = 0.09]{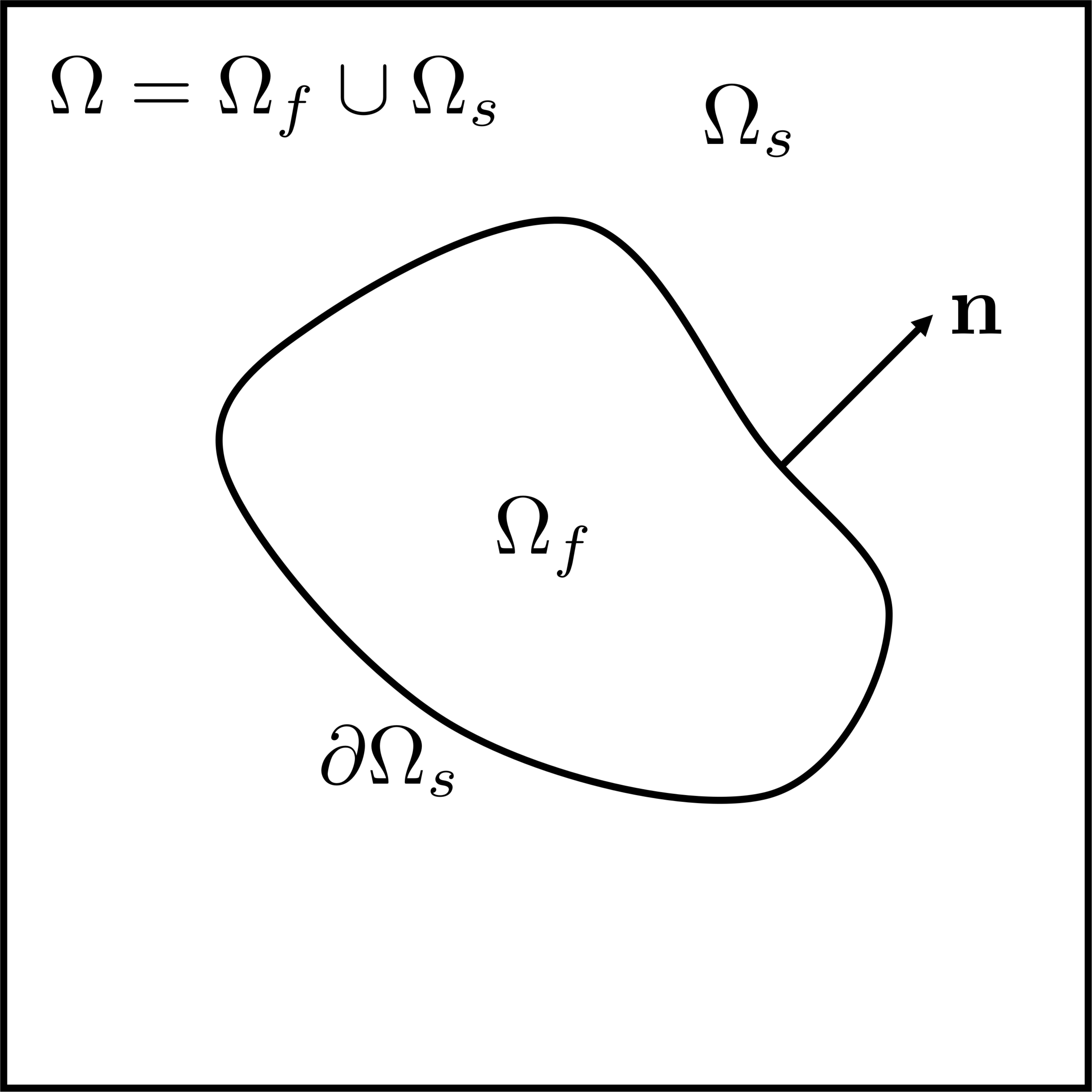}
   \caption{Schematic of a regular computational domain $\Omega$ with an embedded irregular fluid region $\Omegaf$. The solid domain is defined as $\Omegas = \Omega \setminus \Omegaf$. The fluid-solid interface $\partial \Omegas$ (or $\partial \Omegaf$) has the unit normal vector $\n$ pointing out from the fluid and into the solid. }
     \label{fig_domain_schematic}
\end{figure}

Consider an irregular fluid domain $\Omegaf$ embedded into a larger, regular computational domain $\Omega$, as shown in Fig.~\ref{fig_domain_schematic}. Define $\Omega \setminus \Omegaf = \Omegas$ as the fictitious solid domain and $\n$ as a unit outward normal of the fluid-solid interface $\partial \Omegas$ (or $\partial \Omegaf$).  With $q$ as the scalar quantity of interest, $\kappa$ as the diffusion coefficient, and $f$ as a source term, Sakurai et al.~\cite{Sakurai2019} extended the Poisson equation defined in the fluid region $\Omegaf$ 
\begin{equation}
-\nabla \cdot \kappa \; \grad q = f,
\label{eq_non_penalized_poisson}
\end{equation}
satisfying inhomogeneous Neumann/flux boundary conditions on $\partial \Omegas$
\begin{equation}
- \kappa  \; \n \cdot \grad q  = g,
\label{eq_neuman_bc}
\end{equation}
to the entire computational domain $\Omega$ using the flux-based VP  approach. The extended domain Poisson equation satisfying the inhomogeneous flux boundary conditions on the interface reads as
\begin{equation}
-\div \left[\left\{ \kappa \left(1 - \chi\right) + \eta \chi \right\} \grad q \right] = \left(1 - \chi\right) f + \div \left(\chi \vbeta\right) - \chi \div \vbeta. \label{eqn_vp_poisson}
\end{equation}
Here, $\eta$ is the penalization parameter, $\chi(\x)$ is an indicator function whose value is $1$ in the solid region and $0$ in the fluid region, and $\fb = \div \left(\chi \vbeta\right) - \chi \div \vbeta$ is an additional forcing term required to impose the flux boundary conditions on $\partial \Omegas$. The vector-valued flux-forcing function $\vbeta(\x)$ is selected such that $\vbeta \cdot \n = -g$ on the interface. In the limit of $\eta \rightarrow 0$, the solution to the volume penalized  (VP) Poisson equation converges to the solution of non-penalized Poisson equation (Eqs.~\eqref{eq_non_penalized_poisson} and~\eqref{eq_neuman_bc}). A formal derivation of Eq.~\eqref{eqn_vp_poisson} is provided in Appendix~\ref{sec_neumann_derivation}.  As noted in Thirumalaisamy et al.~\cite{Thirumalaisamy2021}, the flux-based VP approach allows $\kappa$ and $g$ to vary spatially as well.

\subsection{The Robin problem}

Next, we consider the inhomogeneous Robin boundary conditions of the type 
\begin{equation}
\zeta\; q + \kappa\; \n \cdot \grad q = -g
\label{eqn_robin}
\end{equation}
on the fluid-solid interface $\partial \Omegas$. Appendix~\ref{sec_robin_derivation} derives the flux-based VP Poisson equation for the Robin problem, which reads as
\begin{equation}
\zeta [ \div (\chi \n) - \chi \div \n ] q -\grad \cdot [\left\{ \kappa \left(1 - \chi\right) + \eta \chi \right\} \; \grad q] = (1- \chi) f + \div (\chi \vbeta) - \chi \div \vbeta.
\label{eqn_vp_robin}
\end{equation} 
In the equation above, the flux-forcing function satisfies the requirement of $\vbeta \cdot \n = -g$. The unit normal vector $\n$ appearing in the first term of  Eq.~\eqref{eqn_vp_robin} can be computed numerically using a signed distance function as explained later in Sec.~\ref{sec_level_set}. In our formulation, $\zeta$, $\kappa$, and $g$  are allowed to vary spatially. 

\subsection{Multiple interfaces and coupled volume penalized equations}

The VP Poisson equations (Eqs.~\eqref{eqn_vp_poisson} and~\eqref{eqn_vp_robin}) can also be generalized to handle multiple interfaces within the computational domain $\Omega$. For some of these interfaces, Dirichlet boundary conditions may also be prescribed.  Following Thirumalaisamy et al.~\cite{Thirumalaisamy2021}, the generalized form of the VP Poisson equation satisfying Neumann and Dirichlet boundary conditions reads as
\begin{equation}
-\div \left[\left\{ \kappa \left(1 - \sum_{j = 1}^{N} \chi^{\text{n}}_j\right) + \sum_{j = 1}^{N} \eta \chi^{\text{n}}_j \right\} \grad q \right] =
\left(1 - \sum_{j = 1}^{N}  \chi^{\text{n}}_j \right)f + \sum_{j = 1}^{N} \left\{\div \left(\chi^{\text{n}}_j \vbeta_j \right) - \chi^{\text{n}}_j \div \vbeta_j\right\} -
\sum_{i = 1}^{D} \frac{\chi^{\text{d}}_i \left(q - \qdi\right)}{\eta}.
 \label{eqn_general_vp_poisson}
\end{equation}
For the above equation to hold true, the computational domain $\Omega$ is assumed to consist of disjoint volumetric regions $\Omega_i^{\rm d}$ (for $i = 1, 2, \dots, D$) and $\Omega_j^{\rm n}$ (for $j = 1, 2, \dots, N$),  with imposed Dirichlet ($q = q_i^{\rm d}$) and Neumann ($-\kappa\; \grad q \cdot \n_j = g_j^{\rm n}$) boundary conditions, respectively. Furthermore, the union of $\Omega_i^{\rm d}$ and  $\Omega_j^{\rm n}$ regions defines the total solid domain, i.e., $\Omegas = \Omega_1^{\rm d} \cup \Omega_2^{\rm d} \cup \cdots \Omega_D^{\rm d} \cup \Omega_1^{\rm n} \cup \Omega_2^{\rm n} \cup \cdots \Omega_N^{\rm n}$. In Eq.~\eqref{eqn_general_vp_poisson} the indicator function $\chi^{\text{n}}(\x) (\text{respectively}, \chi^{\text{d}} (\x))$ is $1$ if
$\x \in \Omega^{\text{n}} (\text{respectively}, \Omega^{\text{d}})$ and $0$ if $\x \in \Omega \setminus \Omega^{\text{n}} (\text{respectively}, \Omega^{\text{d}})$. Note that Robin boundary conditions can be easily included in Eq.~\eqref{eqn_general_vp_poisson}, as their form is very similar to the Neumann problem. We omit Robin boundary conditions in the generalized equation written above for brevity.   

The volume penalization approach can also be extended to other governing equations that describe conservation of momentum, energy, species, etc. For example, the VP incompressible Navier-Stokes equations coupled to the flux-based VP advection-diffusion equation satisfying Neumann boundary condition reads as
\begin{align}
\D{\rho \u}{t} + \div \rho\u\u &= -\grad p + \div \left[\mu \left(\grad \u + \grad \u^T\right) \right] + \frac{\chi}{\eta}(\ub - \u) + \f (\x,q,t), \label{eqn_momentum} \\
\div \u &= 0, \label{eqn_continuity} \\
\D{q}{t} + \left(1 - \chi \right) \left(\u\cdot\grad q\right) &= \div \left[\left\{\kappa \left(1 - \chi\right) + \eta \chi \right\} \grad q \right]+\left(1 - \chi\right) f + \div \left(\chi \vbeta\right) - \chi \div \vbeta. \label{eqn_adv_diff}
\end{align} 
In the equations above, $\u(\x, t)$ is the fluid velocity, $\ub(\x,t)$ is the structure velocity, $p(\x, t)$ is the hydrodynamic pressure, $\f$ denotes the momentum body force, $\rho(\x)$ is the mass density, and $\mu(\x)$ is the dynamic viscosity. The equation set~\eqref{eqn_momentum}-\eqref{eqn_adv_diff} is written considering only a single interface in the domain; generalization to handle multiple interfaces is also possible following Eq.~\ref{eqn_general_vp_poisson}. We remark that in the context of fluid-structure interaction (FSI) problems, only velocity matching condition on the fluid-structure interface is required, i.e., $ \u = \ub$ on $\partial \Omegas$ is essential, whereas $ \u = \ub$ in $\Omegas$ is optional. In the volume penalization approach to FSI, both these conditions are imposed through the penalization term $\frac{\chi}{\eta}(\ub - \u)$. Therefore, in the VP momentum equation~\eqref{eqn_momentum}, only Dirichlet boundary conditions have been considered.


\subsection{Interface capturing} \label{sec_level_set}

We use a signed distance function $\phi(\x)$ to implicitly define the fluid-solid interface $\partial \Omegas$. The scalar field $\phi(\x)$ is defined to satisfy the following property: $\phi(\x) > 0$ if $\x\in\Omegaf$, $\phi(\x) < 0$ if $\x\in\Omegas$ and $\phi(\x) = 0$ is $\x \in \partial \Omegas$. Moreover, the negative gradient of the signed distance function $\phi(\x)$ gives the unit outward normal vector of the interface, i.e., $\n = -\nabla \phi$.  The signed distance function can also be used to define the indicator function $\chi(\x)$. In this work we use $\phi(\x)$ to define two types of indicator functions: one is smooth and continuous and written as
\begin{align}
\chi(\x)&=
\begin{cases}
       1,  & \phi(\x) < - \nsmear \; h,\\
       1 - \frac{1}{2}\left(1 + \frac{1}{\nsmear h} \phi(\x)+ \frac{1}{\pi} \sin\left(\frac{\pi}{ \nsmear h} \phi(\x)\right)\right) ,  & |\phi(\x)| \le \nsmear\; h,\\
        0,  & \textrm{otherwise},
\end{cases}       \label{eqn_continuous_chi}
\end{align}
and the other one is discontinuous, which reads as 
\begin{align}
\chi(\x)&=
\begin{cases}
       1,  & \phi(\x) < 0,\\
       \half,  & \phi(\x) = 0,\\
        0,  & \textrm{otherwise}.
\end{cases}       \label{eqn_discontinuous_chi}
\end{align}
In the Eq.~\eqref{eqn_continuous_chi} above, $\nsmear\in\mathbb{R}$ is the number of grid cells over which the indicator function is smoothed on either side of the interface and $h$ is the grid cell size.


\section{Discrete equations}

\begin{figure}[]
\centering
   \includegraphics[scale = 0.07]{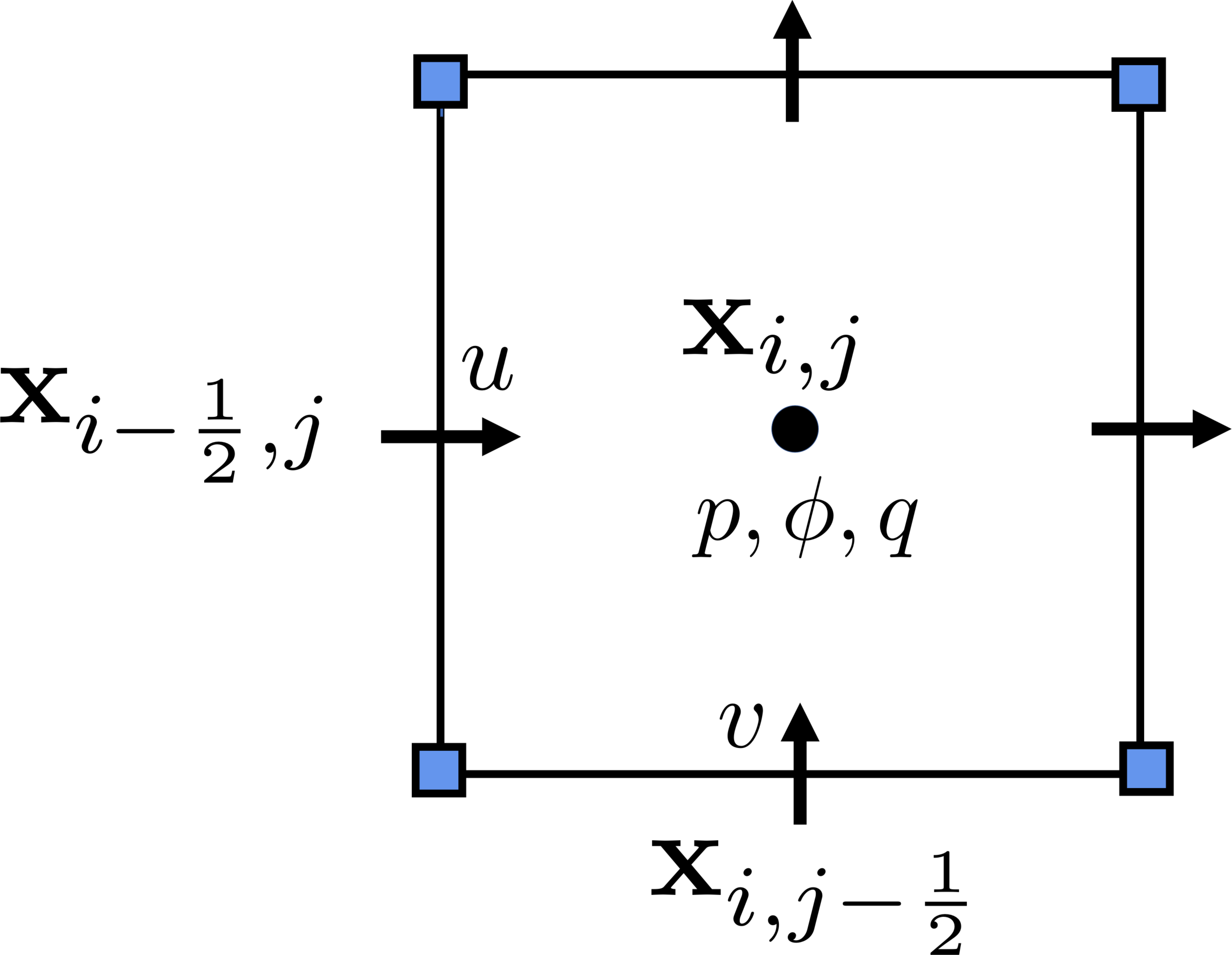}
   \caption{Schematic of a single Cartesian grid cell along with the placement of various variables: the velocity components are stored at the edge centers (black, $\rightarrow$); the fluid pressure $p$, the transported quantity $q$ and the signed distance function $\phi$ are stored at the cell centers (black, $\circ$).}
     \label{fig_cartesian_cell_schematic}
\end{figure}

We use second-order finite difference stencils to discretize the spatial derivative terms of cell-centered Poisson and face-centered momentum equations (Eqs.~\eqref{eqn_vp_poisson} and~\eqref{eqn_momentum}, respectively) on a Cartesian grid. Fig.~\ref{fig_cartesian_cell_schematic} shows a schematic representation of a two-dimensional Cartesian grid cell, in which the velocity components are stored on edge centers (face centers in three spatial dimensions), whereas the transported variable $q$, the fluid pressure $p$, and the signed distance function $\phi$ are stored at the cell center. The computational domain $\Omega$ is discretized into $N_x\times N_y$ Cartesian grid cells with  mesh spacing $\Delta x$ and  $\Delta y$ in  the $x$- and $y$-direction, respectively. In this work we use equal mesh spacing in the two directions, i.e., $\Delta x = \Delta y = h$. In what follows next, we primarily focus on the discretization of the VP Poisson Eq.~\eqref{eqn_vp_poisson} for the Neumann problem; details on the spatiotemporal discretization of the VP incompressible Navier-Stokes equations can be found in our prior works~\cite{Nangia2019MF,BhallaBP2019,Dafnakis2020}.

Referring to Fig.~\ref{fig_cartesian_cell_schematic}, let $(i,j)$ denote the cell index, $(i-\half,j)$ denote the lower $x$ edge index and $(i,j-\half)$ denote the lower $y$ edge index. Then the discretized form of the VP Poisson Eq.~\ref{eqn_vp_poisson} in two spatial dimensions reads as
\begin{align}
\label{eqn_discrete_poisson}
 \{ \psi_{i+\half,j} + \psi_{i-\half,j} +   \psi_{i,j+\half} + \psi_{i,j -\half} \} \; q_{i,j}  - \psi_{i+\half,j} \; q_{i+1,j} - \psi_{i-\half,j} \; q_{i-1,j} & - \psi_{i,j+\half}\; q_{i,j+1} \nonumber  \\ & - \psi_{i,j-\half}\; q_{i,j-1} = S_{i,j},
\end{align}
in which 
\begin{subequations}
\begin{align}
\psi_{i+\half,j} &= \frac{1}{\Delta x^2}\big\{  \kappa \;(1-\chi+\eta\chi) \big\}_{i+\half,j}   \\ 
\psi_{i-\half,j} &=  \frac{1}{\Delta x^2} \big\{ \kappa \; (1-\chi+\eta\chi) \big\}_{i-\half,j}   \\ 
\psi_{i,j+\half} &=  \frac{1}{\Delta y^2} \big\{  \kappa \;(1-\chi+\eta\chi) \big\}_{i, j+ \half}  \\ 
\psi_{i,j-\half} &= \frac{1}{\Delta y^2} \big\{  \kappa \;(1-\chi+\eta\chi) \big\}_{i, j- \half},
\end{align}
\end{subequations}
and the right hand side term $S_{i,j}$ is given by
\begin{align}
S_{i,j} = & \left(1-\chi_{i,j} \right )f_{i,j} + \frac{\left(\chi\beta \right )_{i+\half,j} - \left(\chi\beta \right )_{i-\half,j}}{\Delta x} +\frac{\left(\chi\beta \right )_{i,j+\half} - \left(\chi\beta \right )_{i,j-\half}}{\Delta y} \nonumber \\
 & \quad -\chi_{i,j}\left( \frac{\beta_{i+\half,j} - \beta_{i-\half,j}}{\Delta x} + \frac{\beta_{i,j+\half} - \beta_{i,j-\half}}{\Delta y}\right ).
\end{align}
Analogous discretization formulas can be written for the three-dimensional VP Poisson equation. In the discretized Eq.~\eqref{eqn_discrete_poisson} written above, the indicator function $\chi$ and the diffusion  coefficient $\kappa$ are required at the edge centers; these properties are first defined at the cell centers and then interpolated onto the edge centers using a second-order accurate linear interpolation scheme.  The flux-forcing function $\vbeta(\x)$ is also required at the edge centers; methods to construct $\vbeta$ are discussed next. 


\subsection{Construction of flux-forcing functions}  \label{sec_beta}

The vector-valued flux-forcing function $\vbeta(\x)$ plays a crucial role in imposing the desired inhomogeneous Neumann and Robin boundary conditions on the interface. In this section, we introduce three approaches to construct $\vbeta(\x)$, namely Approach A, B, and C. The three approaches are in increasing order of generality. While Approach A is specialized for the Neumann problem, Approaches B and C are equally applicable to the Robin problem. 

\subsubsection{Approach A: Analytical construction of spatially varying $g$}
 
 Consider for a moment that the solution to the non-penalized Poisson equation with inhomogeneous Neumann boundary conditions on $\partial \Omegas$ is known. Denote the exact solution by $\qexact$. If the flux-forcing function is taken to be of the form $\vbeta(\x) = \kappa\; \grad \qexact(\x)$, then it satisfies the requirement of $\vbeta \cdot \n = -g$ on $\partial \Omegas$. In practice the solution to the Poisson Eq.~\eqref{eq_non_penalized_poisson} with boundary condition~\eqref{eq_neuman_bc} is sought and not known \emph{a priori}. However, if an analytical approximation $\qinexact$ to the exact solution $\qexact$ exists, such that $  \n \cdot  \grad \qinexact  = \n \cdot \grad \qexact$ on $\partial \Omegas$, then  
 \begin{equation}
 \vbeta (\x) = \kappa\; \grad \qinexact (\x), 
 \label{eqn_approachA}
 \end{equation} 
 can be prescribed as a flux-forcing function. Away from the interface, the approximation $\qinexact$ can be close to or very different from $\qexact$, depending upon whether a continuous or a discontinuous indicator function $\chi$ is used. We denote the analytical construction of $\vbeta$ as Approach A. In component form,  Approach A is written as
\begin{subequations}
\begin{align}
 \beta_{i-\half,j} &= \left( \kappa \D{\qinexact} {x} \right)_{\x_{i-\half,j}},  \\
 \beta_{i,j-\half} &= \left( \kappa \D{\qinexact} {y} \right)_{\x_{i,j-\half}}.
\end{align}
\end{subequations}

Approach A was employed in Sakurai et al.~\cite{Sakurai2019} and Thirumalaisamy et al.~\cite{Thirumalaisamy2021} to demonstrate the feasibility of flux-based volume penalization method to solve PDEs with flux boundary conditions in complex domains, but is quite restrictive in practice as discussed next.  

\subsubsection{Approach B: Numerical construction of spatially constant $g$}

Although Approach A allows for imposing spatially varying $g$ values on the interface, approximating an analytical solution to the exact solution near the interface is a non-trivial task, especially if the interface is geometrically complex. However, if $g$ is spatially constant, then constructing $\vbeta$ is easy. This is achieved by taking $\vbeta(\x) = -g\; \n(\x)$, as it satisfies the requirement of $\vbeta  \cdot \n = -g$ on $\partial \Omegas$. Now, recalling from Sec.~\ref{sec_level_set} that the negative gradient of the signed distance function $\phi(\x)$ is the continuous normal vector field $\n(\x)$, the flux-forcing function can be constructed numerically for an irregular boundary as
\begin{equation}
\vbeta(\x) = -g\; \n(\x) = g\; \grad \phi(\x). 
\label{eqn_approachB}
\end{equation} 
We denote the numerical construction of spatially constant $g$ value on the interface as Approach B, which in component form is written as
\begin{subequations}
\begin{align}
 \beta_{i-\half,j} &= g \; \left(  \frac{ \phi_{i,j} - \phi_{i-1,j} } { \Delta x} \right), \\
 \beta_{i,j-\half} &= g \; \left(  \frac{\phi_{i,j} - \phi_{i, j-1}} { \Delta y} \right) .
\end{align}  
\end{subequations} 

\subsubsection{Approach C: Numerical construction of spatially varying $g$}

\begin{figure}[]
\centering
 \subfigure[Cartesian grid and an embedded interface]{
\includegraphics[scale = 0.09]{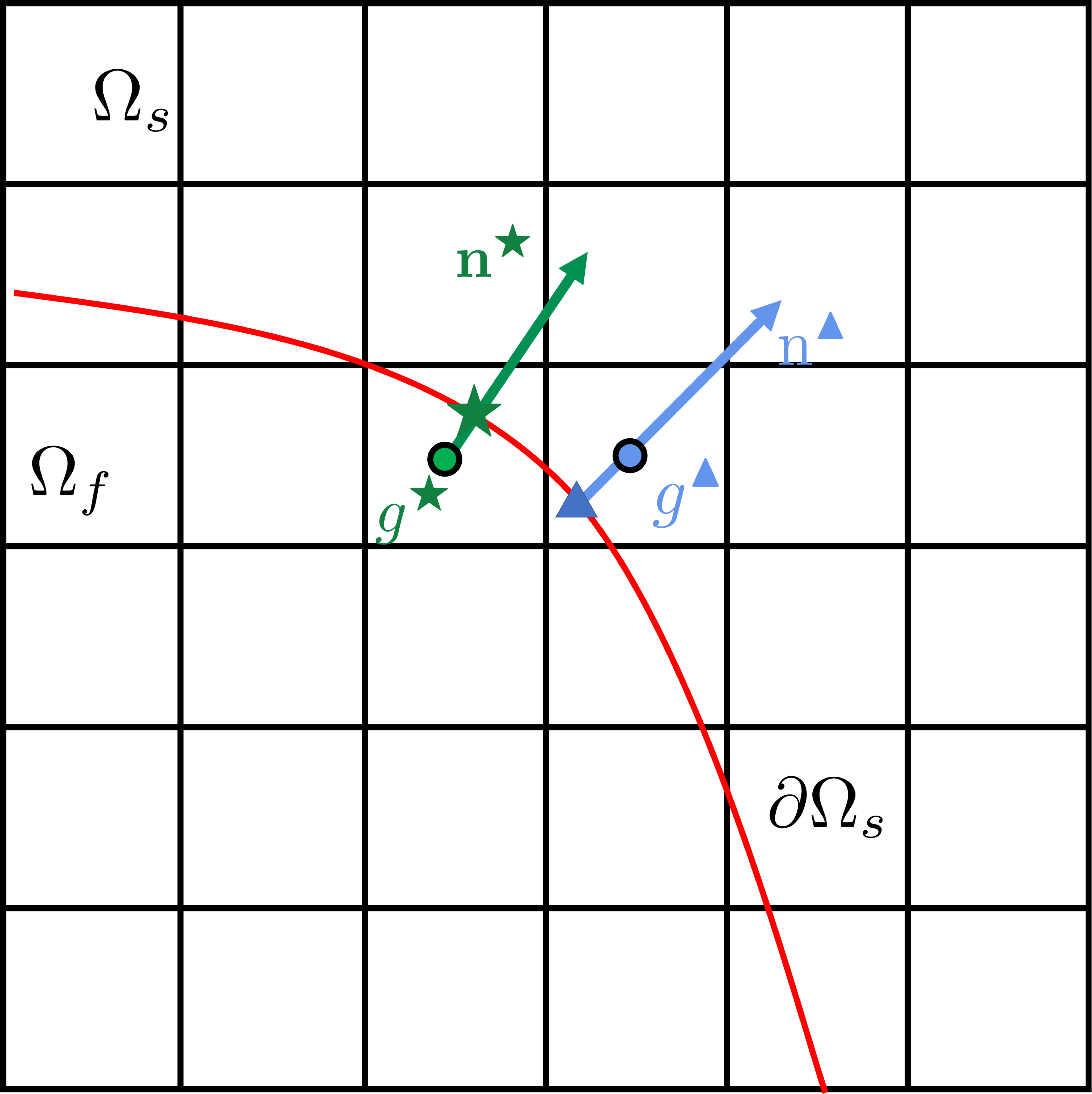}
\label{fig_g_c_schematic}
}
\subfigure[Propagation of interface cell $g$ values]{
\includegraphics[scale = 0.09]{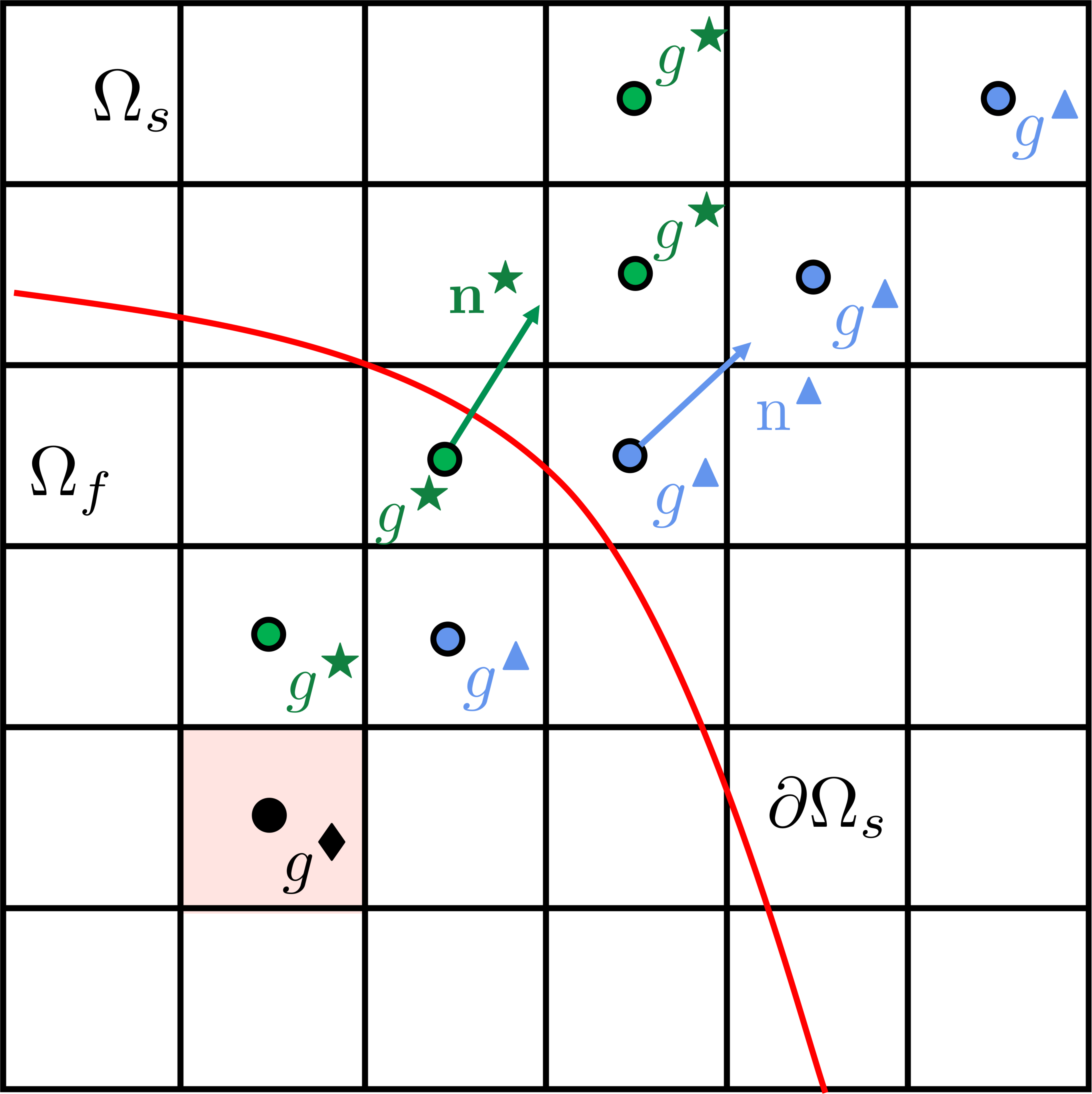}
\label{fig_g_prop_schematic}
}
\caption{Approach C for constructing the flux-forcing function $\vbeta$. \subref{fig_g_c_schematic} Schematic of a Cartesian grid (black lines) with an embedded fluid-solid interface $\partial \Omegas$ (red line). Two interface cells are highlighted in the figure: one whose cell center $\x^{\bigstar}$ lies in the fluid region  and the other whose cell center $\x^\blacktriangle$ lies in the the solid region. For these interface cells the corresponding function values $g^{\bigstar}$ and  $g^\blacktriangle$ and the outward unit normal vectors $\n^{\bigstar}$ and  $\n^\blacktriangle$ are shown. \subref{fig_g_prop_schematic} Propagation of the $g$ value into the domain following an interface cell normal and a grid cell where multiple normals intersect shown in light pink color. Out of the $g^{\bigstar}$ and  $g^\blacktriangle$ values arriving at the shaded cell, the one with the larger modulus is chosen. }
\label{fig_approach_c}
\end{figure}

As a generalization of Approach B, the flux-forcing function can be taken as 
\begin{equation}
\vbeta(\x) = -g(\x)\; \n(\x) = g(\x) \grad \phi(\x).
\label{eqn_approachC}
\end{equation}  
However, Eq.~\eqref{eqn_approachC} poses a challenge of extending the codimension-1 boundary condition function $g$ defined over the interface to a  codimension-0 function $g(\x)$ defined in the neighborhood of the interface. Although there are several ways  to achieve this function extension (in absence of a constraint), in this work we follow a simple strategy of propagating the interfacial $g$ values to the neighboring grid cells along the interface normal. More specifically, consider a fluid-solid interface $\partial \Omegas$ embedded into a Cartesian grid as shown in Fig.~\ref{fig_g_c_schematic}. The signed distance function $\phi(\x)$ can be used to identify the grid cells through which the interface passes. Denote these grid cells as \emph{interface cells}. Fig.~\ref{fig_g_c_schematic} highlights two such interface cells: one whose cell center $\x^\bigstar$ lies in the fluid region and the other whose cell center $x^\blacktriangle$ lies in the solid region. The normal vector of the interface cells is also known from the signed distance function: $\n^{\bigstar} = (-\grad \phi)^{\bigstar}$ and  $\n^{\blacktriangle} = (-\grad \phi)^{\blacktriangle}$.  Next, the $g$ value at the cell center of an interface cell is set equal to the closest interfacial $g$ value:
\begin{subequations}
\begin{align}
g^{\bigstar} &\leftarrow g(\x_{\partial \Omegas}^\bigstar), \\
g^{\blacktriangle} &\leftarrow g(\x_{\partial \Omegas}^\blacktriangle),
\end{align}
\end{subequations}
in which $\x_{\partial \Omegas}^\bigstar = \x^\bigstar + \phi^\bigstar \n^\bigstar$ and $\x_{\partial \Omegas}^\blacktriangle = \x^\blacktriangle + \phi^\blacktriangle \n^\blacktriangle$ are the \REVIEW{closest} points on the interface to the cell centers $\x^\bigstar$ and $\x^\blacktriangle$, respectively. Note that  the $g$ function on the interface is prescribed and therefore, $g(\x_{\partial \Omegas}^\bigstar)$ and  $g(\x_{\partial \Omegas}^\blacktriangle)$ are known \emph{a priori}. In the next part of the algorithm, $g^{\bigstar}$ and $g^{\blacktriangle}$ values are propagated to the grid cells that are within a distance of  $\nprop\;h$ to the interface cells along $\pm\; \n^{\bigstar}$ and $\pm\; \n^{\blacktriangle}$ directions, respectively. This procedure is pictorially described in Fig.~\ref{fig_g_prop_schematic}. The number of grid cells $\nprop$ to which $g$ values are propagated depends upon the choice of the indicator function $\chi$---we will explore the effect of $\nprop$ on the solution accuracy in Sec.~\ref{sec_results_and_discussion}. Note that propagating $g$ values along the normal directions may lead to a situation of conflict at a grid cell where two or more interface cell normals intersect. This situation is shown for the shaded cell in Fig.~\ref{fig_g_prop_schematic} where the two normals $\n^{\bigstar}$ and $\n^{\blacktriangle}$ intersect. For such cells, a $g$ value with the larger modulus is chosen: 
\begin{equation}
g^{\blacklozenge}= \mathrm{maxmodulus} (g^{\bigstar}, g^{\blacktriangle}).
\label{eq_g_prop_max}
\end{equation}
We also considered an average and the minimum modulus of $g$ at the conflicted cells; these choices however reduced the order of accuracy of the solution for the continuous/smoothed indicator function.  Note that for the discontinuous indicator function,  the $g$ value at the conflicted cells does not matter much for the solution accuracy. This is because such cells are generally located far away from the interface where the discontinuous indicator function is already zero. Nevertheless, we always make use of Eq.~\eqref{eq_g_prop_max} even for the discontinuous indicator function in this work. With $g$ values defined at the cell centers, the component form of $\vbeta$  reads as
\begin{subequations}
\begin{align}
 \beta_{i-\half,j} &= \left( \frac{g_{i-1,j} + g_{i,j}}{2} \right) \left(  \frac{ \phi_{i,j} - \phi_{i-1,j} } { \Delta x} \right), \\
 \beta_{i,j-\half} &= \left( \frac{g_{i,j-1} + g_{i,j}}{2} \right) \left(  \frac{\phi_{i,j} - \phi_{i, j-1}} { \Delta y} \right) .
\end{align}  
\end{subequations} 

\REVIEW{The propagation strategy of Approach C can also be implemented by solving a hyperbolic equation of the form
\begin{equation}
\D{g(\x)}{\tau} + \n(\x) \cdot \grad g(\x) = 0.
\label{eqn_hyperbolic_g}
 \end{equation}
The equation above can be integrated over a pseudo-time interval $\Delta \tau$ that is directly related to the propagation distance.  However, the test examples of Sec.~\ref{sec_results_and_discussion} show that the method of $g$ propagation described in Approach C is quite effective in imposing the spatially varying flux boundary conditions. Moreover, it does not require solving any additional partial differential equation.   
 
Note that there can be other ways of extending the flux-forcing function in the vicinity of the interface as discussed at the beginning of this section.  One straightforward approach is to extend the $\V{\beta}$ function defined over the interface $\Omegas$ to a flux-forcing function valid near the interface using the top hat  or a Gaussian bell-like function, which we refer to as Approach D. However, as demonstrated in Appendix~\ref{sec_other_ext_results}, this particular  function continuation approach does not produce satisfactory results; the numerical and actual solutions differ significantly and the numerical scheme does not converge under grid refinement in any norm. In contrast, Approach C produces the correct solution and a convergent numerical scheme. This also highlights the non-triviality in allowing spatially varying Neumann/Robin boundary conditions in the flux-based VP method.  
 }


%% file: Results_and_discussion.tex

In this section we discretely solve the volume penalized Poisson Eqs.~\ref{eqn_vp_poisson} and~\ref{eqn_vp_robin} satisfying inhomogeneous Neumann and Robin boundary conditions, respectively, to assess the accuracy of the numerical solutions. \REVIEW{We use the flexible GMRES (FGMRES) iterative solver  with a tight relative residual tolerance of $10^{-12}$ to solve the system of linear equations.} The order of accuracy results presented here are computed only in the fluid domain and are determined based on the $L^1$ and $L^\infty$ norm of the error (denoted \REVIEW{$\mathcal{E}^1$} and \REVIEW{$\mathcal{E}^\infty$}, respectively) between the numerical and analytical~\footnote{Analytical solution of the non-penalized equation is used for computing the error.} solutions. \REVIEW{Since the VP method is expected to produce a non-uniform convergence rate under grid refinement because of the delta function formulation (see Appendices~\ref{sec_neumann_derivation} and~\ref{sec_robin_derivation} for derivation), we curve-fit the error data and report the slope/convergence rate, denoted $m$ and the coefficient of determination, denoted $R^2$, in each case.  Appendix~\ref{sec_convergence_data} tabulates the error data.} The spatial convergence rate of the error is shown for both continuous \REVIEW{(denoted  $\mathcal{E}_\text{c}^{\infty}$ and  $\mathcal{E}_\text{c}^1$)} and discontinuous \REVIEW{(denoted  $\mathcal{E}_\text{d}^{\infty}$ and  $\mathcal{E}_\text{d}^1$)} indicator functions. We consider two- and three-dimensional examples involving constant and spatially varying flux boundary conditions on $\partial\Omegas$. In the test examples, the fluid region $\Omegaf$ is embedded into a larger computational domain $\Omega$ with Dirichlet boundary conditions imposed on the external boundary $\partial\Omega$ of the domain. The computational domain is discretized into $N\times N$ and $N\times N \times N$ grid cells for the two- and three-dimensional examples, respectively.  The penalization parameter $\eta$ is taken to be $10^{-8}$ \REVIEW{(Appendix~\ref{sec_eta_effect} considers the effect of $\eta$ on the convergence rate)} and the diffusion coefficient $\kappa$ is taken to be 1 for all of the tests. While imposing the Robin boundary conditions we take $\zeta$ to be 1 in the test examples. The numerical solutions are presented for Approach C and where applicable, results obtained from Approach C are compared against Approach A or B.  \REVIEW{Since Approach A constructs the flux-forcing function from the known solution to the problem, for a given indicator function $\chi$, Approach A is expected to perform better than or at least as well as Approach B and C. This expectation is also confirmed from the tests that follow next.} 

\subsection{Concentric circular annulus with spatially constant flux on the interface} \label{sec_circular_annulus_constant_flux}

We first consider the concentric circular annulus problem from Sakurai et al.~\cite{Sakurai2019} in which different inhomogeneous, but spatially constant,  Neumann boundary conditions are specified on the two interfaces defining the annulus. The inner radius of the annulus is $r_i = \pi/4$ and its outer radius is $r_o = 3\pi/4$. The center of the annulus is positioned at $(\pi,\pi)$. The circular annulus is embedded into a larger computational domain of extents $\Omega \in [0, 2\pi]^2$. The source term of the Poisson equation for this case is
\begin{equation}
f(r) = 16\cos(4r)+\frac{4\sin(4r)}{r},
\label{sakurai_annulus_src_term}
\end{equation}
in which $r = \sqrt{\left(x - \pi\right)^2 + \left(y - \pi\right)^2}$ and the Neumann boundary condition values on the two interfaces are taken to be
\begin{equation}
\left.\frac{\mathrm{d} q}{\mathrm{d} r}\right|_{r = \frac{\pi}{4}} = 3m \hspace{1 pc} \text{and} \hspace{1 pc} \left.\frac{\mathrm{d} q}{\mathrm{d} r}\right|_{r = \frac{3\pi}{4}} = m.
\end{equation}
The exact solution of this problem using the zero-mean condition $\int_{\Omegaf}rq(r)\; \text{d}r = 0$  reads as
\begin{equation}
\qexact(r) = \cos(4r) + \frac{3}{4} m \pi \log(r) - \frac{3}{32} m \pi \left(9 \log\left(\frac{3}{4}\pi\right) - \log\left(\frac{\pi}{4}\right) - 4\right).
\label{sakurai_2d_poisson_exact}
\end{equation}
\REVIEW{The mean of the numerical solution in the fluid region is subtracted as a post-processing step to impose the zero-mean condition numerically.} 

Relatively simple geometry and constant flux boundary conditions of this test problem allows for an analytical construction of $\vbeta$. Indeed, in~\cite{Sakurai2019,Thirumalaisamy2021}, the flux-forcing function $\vbeta$ was constructed analytically (Approach A) as  $\vbeta = \kappa\; \grad \qinexact = \d{\qinexact}{r} \e_r$,  in which $\e_r = \left(\frac{x-\pi}{r}, \frac{y-\pi}{r}\right)$ and $\d{\qinexact}{r}$ is
\begin{equation}
\d{\qinexact}{r} = \left \{  \begin{array}{ll}
      m \left(\frac{4r}{3\pi}\right)^2 \left(4\left(1 - \frac{r}{\pi}\right)\right)^3, & \text{if} \hspace{1pc} 0 \leq r \leq \pi, \\
      0, & \text{otherwise}. \\
\end{array} \right.
\label{eq_sakurai_beta}
\end{equation}
It is to be noted that $\d{\qinexact}{r}$ reduces to $3m$ at $r = \pi/4$ and $m$ at $r = 3\pi/4$, respectively. 

\begin{figure}[]
\centering
\subfigure[{\REVIEW{Spatial convergence rate for Approach A}}]{
\includegraphics[scale = 0.08]{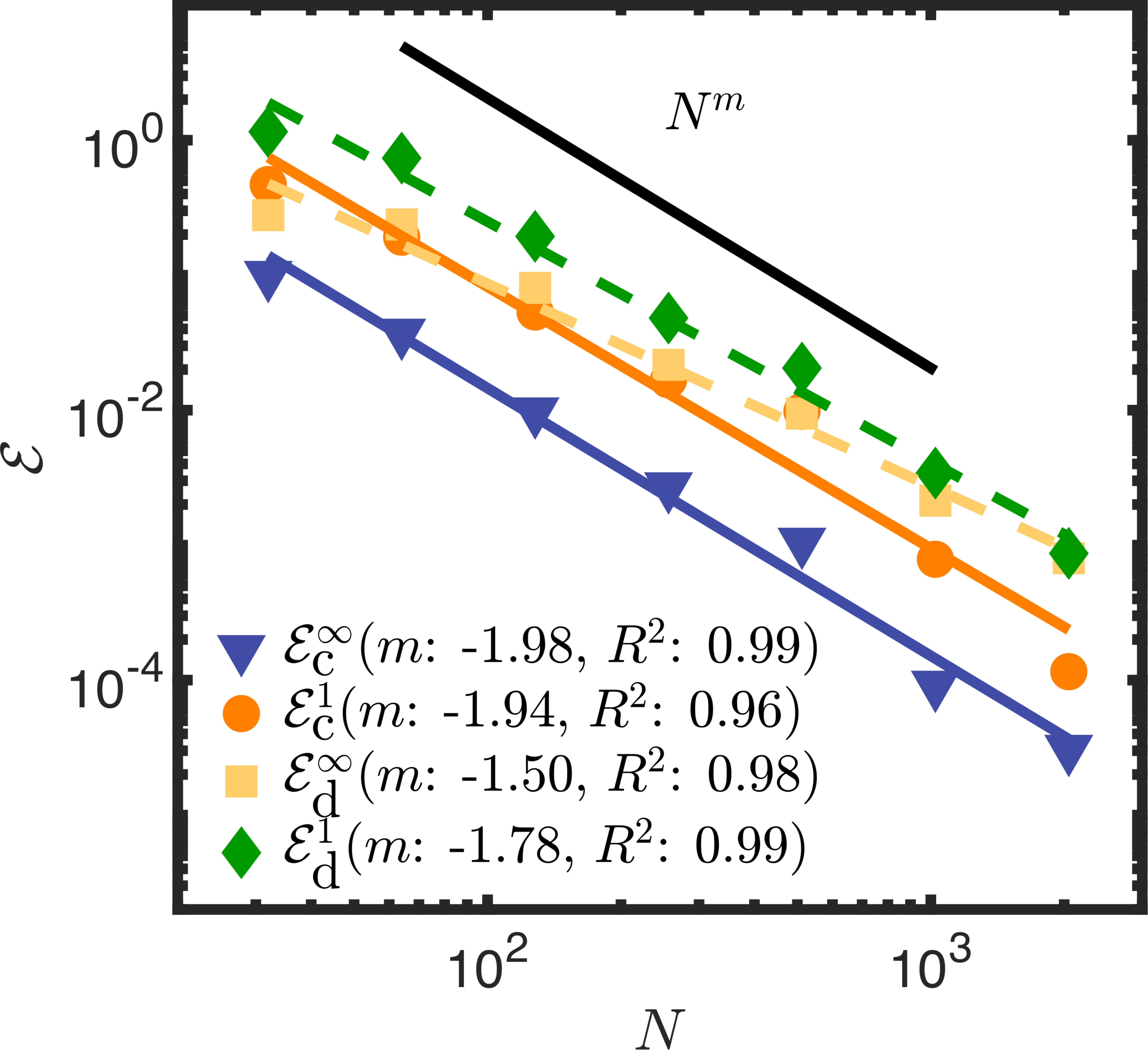}
\label{fig_circular_annulus_ooa_approach_A}
}
\subfigure[{\REVIEW{Spatial convergence rate for Approach B}}]{
\includegraphics[scale = 0.08]{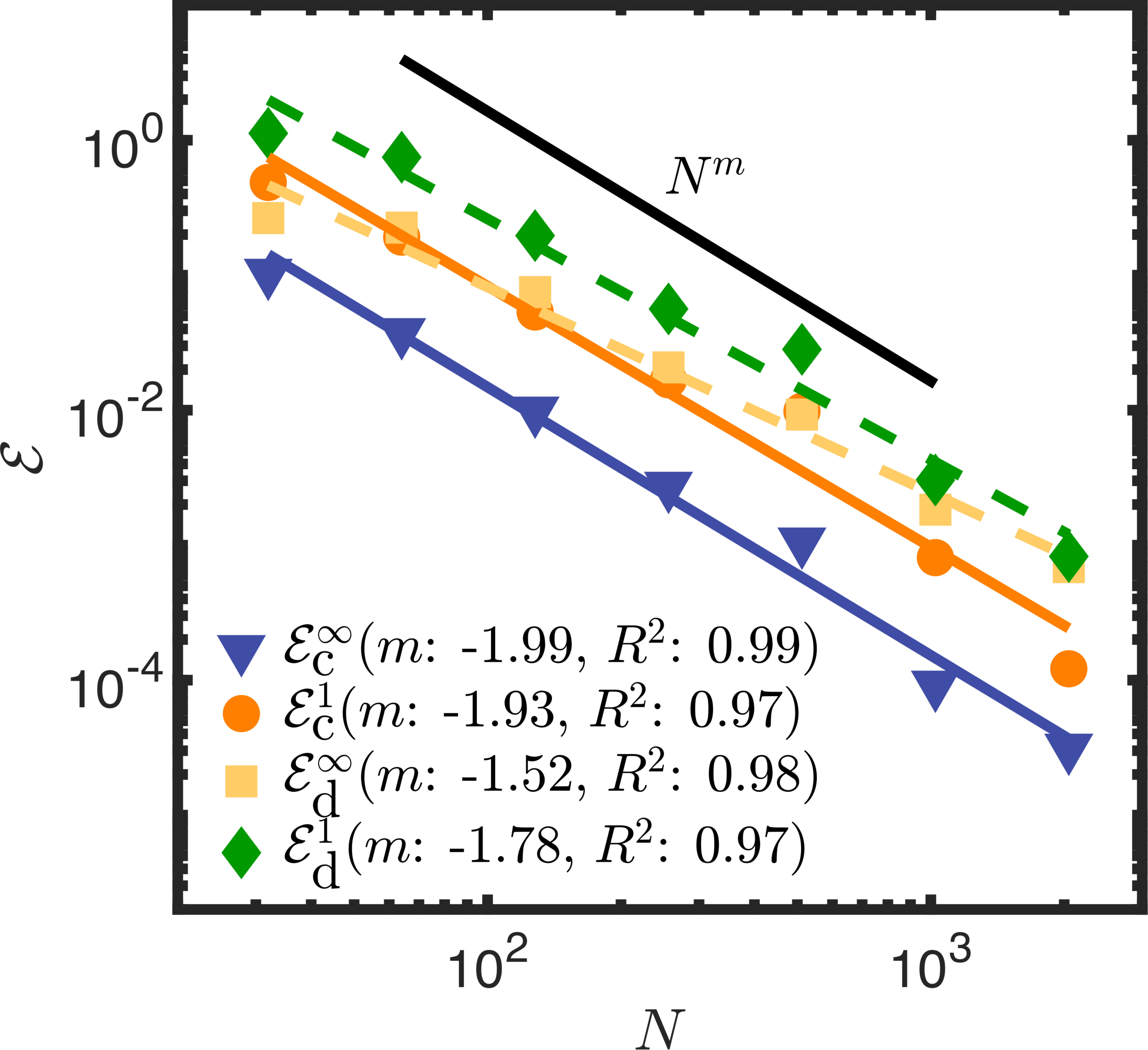}
\label{fig_circular_annulus_ooa_approach_B}
}
\subfigure[{Solution variation along $y-$direction}]{
\includegraphics[scale = 0.08]{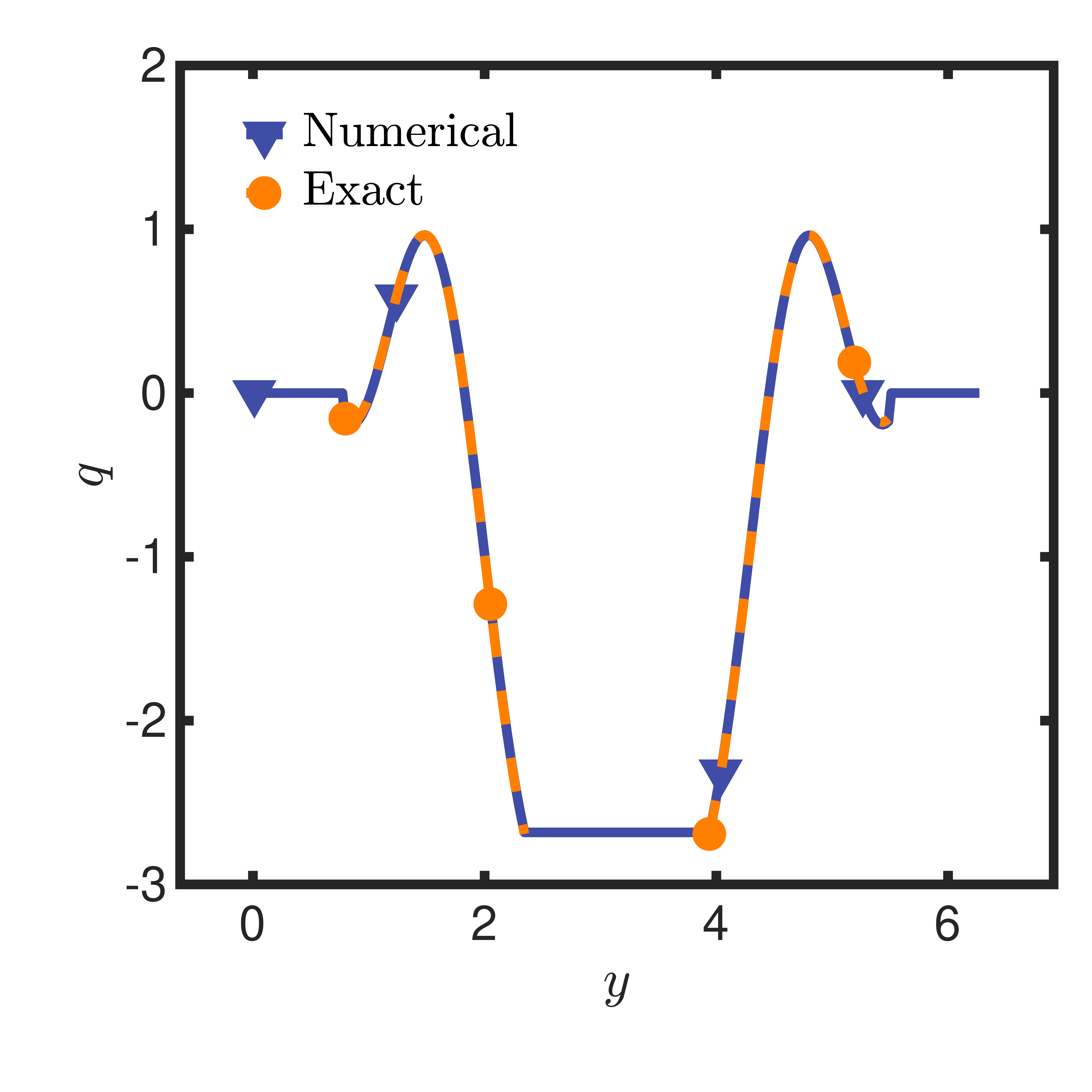}
\label{fig_circular_annulus_xcut}
}
\subfigure[{Solution variation along $x-$direction}]{
\includegraphics[scale = 0.08]{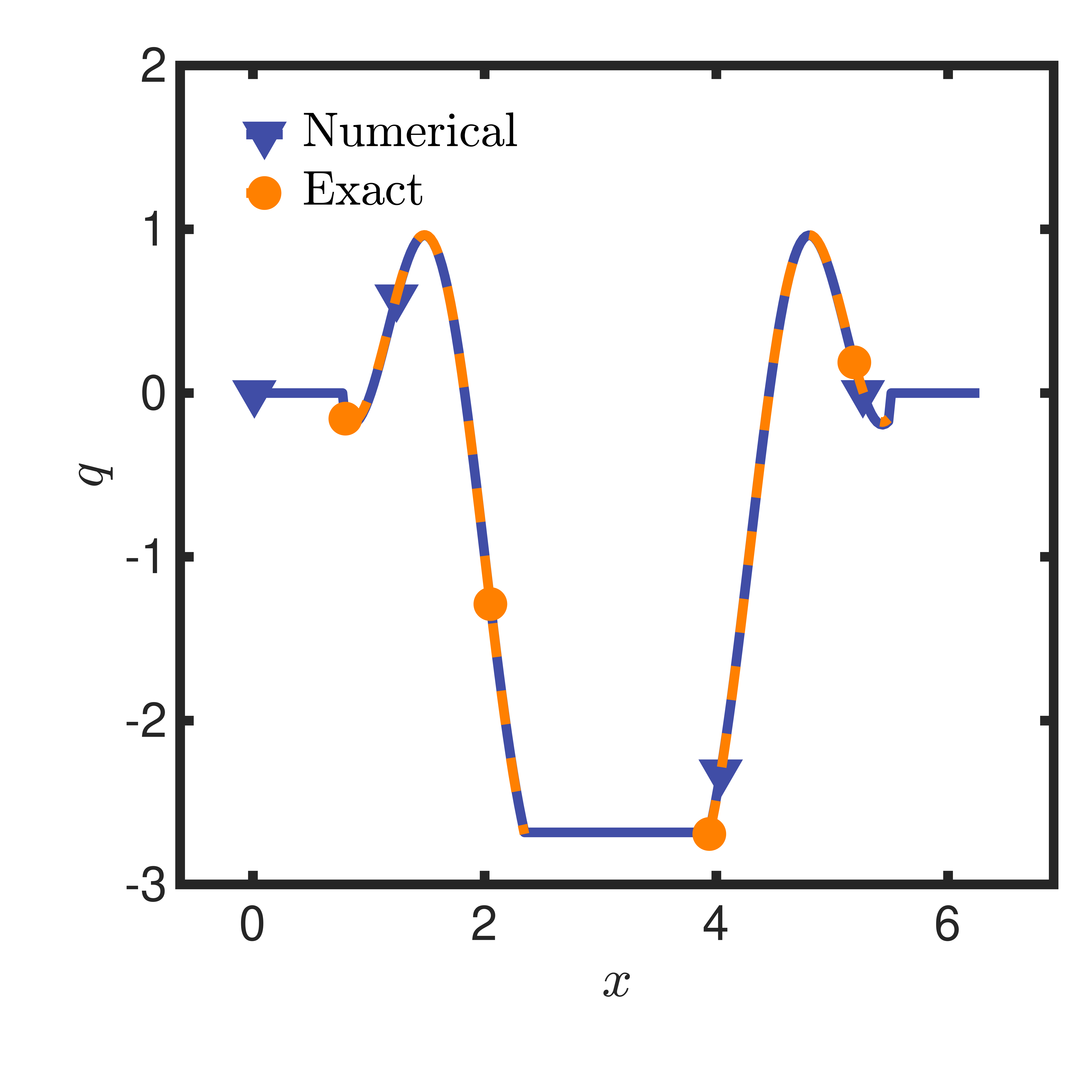}
\label{fig_circular_annulus_ycut}
}
 \caption{Concentric circular annulus with constant flux on the two interfaces using Approach A and B. Error norms \REVIEW{$\mathcal{E}^1$} and \REVIEW{$\mathcal{E}^\infty$} as a function of grid size $N$ using the continuous (solid lines with symbols) and discontinuous (dashed lines with symbols) indicator functions for~\subref{fig_circular_annulus_ooa_approach_A} \REVIEW{Approach A}; and  \subref{fig_circular_annulus_ooa_approach_B} \REVIEW{Approach B}. The penalization parameter $\eta$ is taken as $10^{-8}$ and  $m$ and $\kappa$ are taken as 1.~\subref{fig_circular_annulus_xcut} Variation of the numerical solution along $y-$direction at a fixed $x = 3.12$ location; and~\subref{fig_circular_annulus_ycut} variation of the numerical solution along $x-$direction at a fixed $y= 3.12$  location using $N = 256$ grid.}
\label{fig_circular_annulus}
\end{figure}

Since the flux boundary condition value is spatially constant on the interface, Approach B is also applicable for this test problem.  The $g$ value required for Approach B is $m$ and $-3m$ on the inner and outer interface, respectively. Figs.~\ref{fig_circular_annulus_ooa_approach_A} and~\ref{fig_circular_annulus_ooa_approach_B} show the order of accuracy of the solution as a function of mesh resolution for Approach A and B, respectively.  \REVIEW{As can be observed in Fig.~\ref{fig_circular_annulus_ooa_approach_A}, for Approach A,  $\mathcal{O}(h^{1.98})$ (respectively,  $\mathcal{O}(h^{1.94})$) convergence rate with an $R^2$ value of 0.99 (respectively, 0.96)  in $L^{\infty}$ (respectively, $L^1$ ) norm is achieved using the continuous indicator function. With the discontinuous $\chi$, $\mathcal{O}(h^{1.50})$ (respectively,  $\mathcal{O}(h^{1.78})$) convergence rate with an $R^2$ value of 0.98 (respectively, 0.99) in $L^{\infty}$ (respectively, $L^1$) norm is achieved. We note that the convergence rate using the continuous indicator function is better than the discontinuous function for Approach A.}  \REVIEW{Looking at Fig.~\ref{fig_circular_annulus_ooa_approach_B}, it is seen that Approach B exhibits a very similar convergence rate as Approach A, but in contrast to Approach A, Approach B is more versatile as it requires only $\grad \phi(\x)$ information, which can be constructed for any irregular boundary~\cite{Baerentzen2005}.} Figs.~\ref{fig_circular_annulus_xcut} and~\ref{fig_circular_annulus_ycut} compare the analytical and numerical solutions along $x-$ and $y-$direction, respectively. Numerical solutions using Approach B and discontinuous indicator function are presented. As can be observed in the figures, an excellent agreement is obtained between the analytical and numerical solutions.   

\begin{figure}[]
\centering
\subfigure[{\REVIEW{Spatial convergence rate for Approach C}}]{
\includegraphics[scale = 0.08]{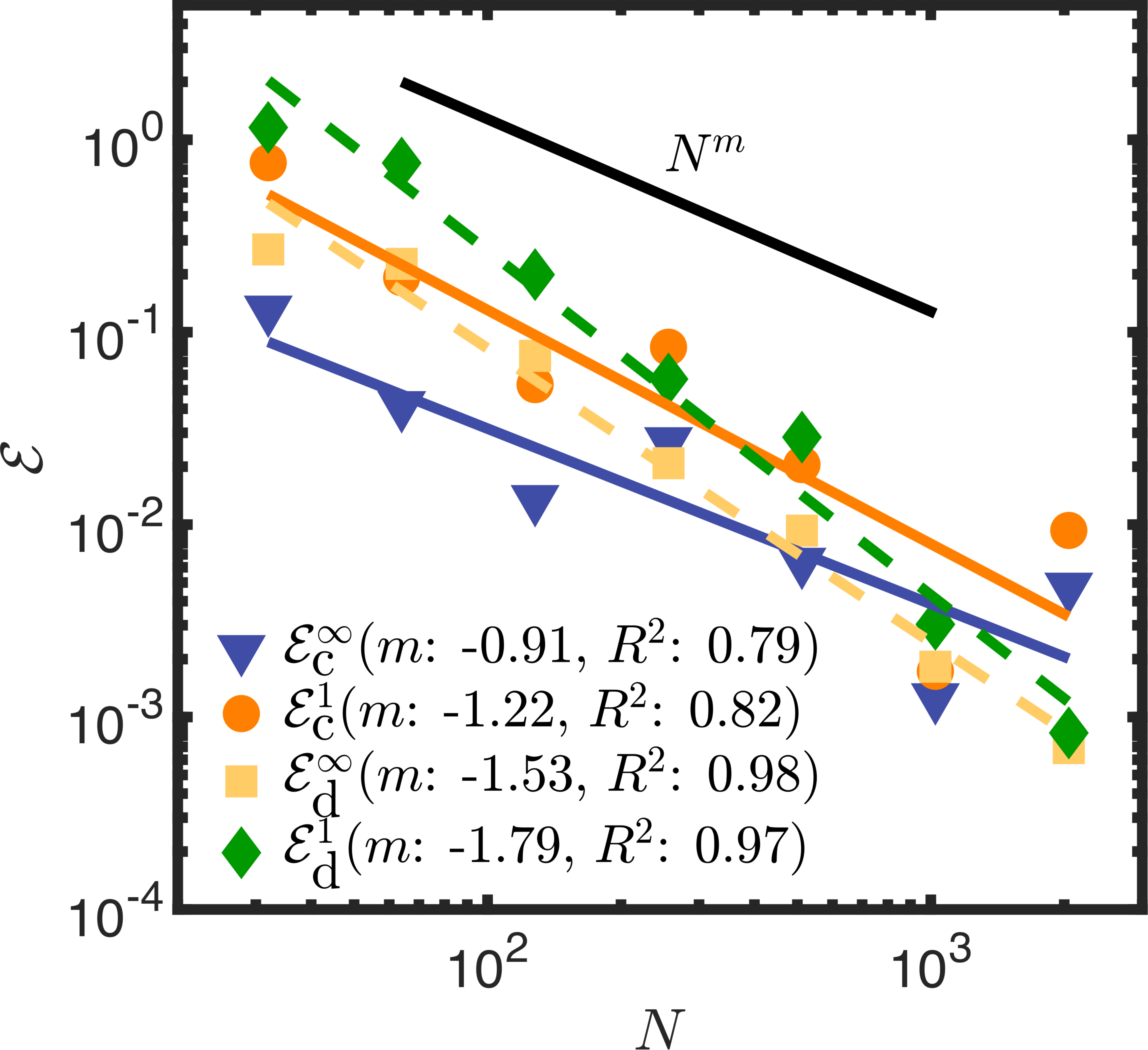}
\label{fig_circular_annulus_ooa_approach_C}
}
\subfigure[{\REVIEW{Effect of number of smear width and propagation cells}}]{
\includegraphics[scale = 0.08]{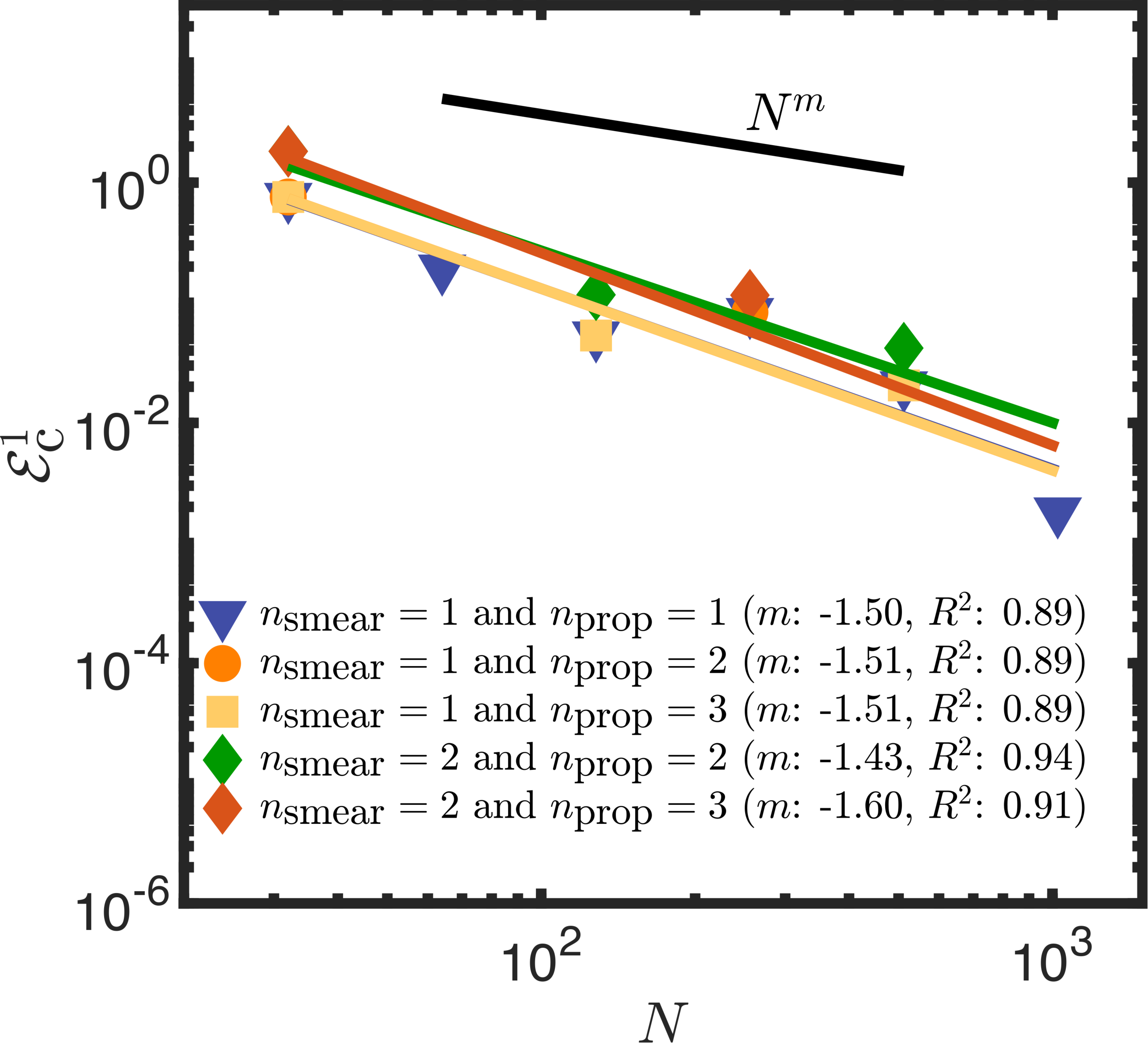}
\label{fig_circular_annulus_effect_of_prop_cells}
}
 \caption{Concentric circular annulus with constant flux on the two interfaces using Approach C. Error norms \REVIEW{$\mathcal{E}^1$} and \REVIEW{$\mathcal{E}^\infty$} as a function of grid size $N$ using the continuous (solid lines with symbols) and discontinuous (dashed lines with symbols) indicator functions.~\subref{fig_circular_annulus_ooa_approach_C} \REVIEW{Convergence rate using} \REVIEW{$\ncells = 1$ and $\nprop = 2$.}~\subref{fig_circular_annulus_effect_of_prop_cells} \REVIEW{Effect of $\ncells$ and $\nprop$ on the solution accuracy.}}
\label{fig_circular_annulus_approach_C}
\end{figure}

Next, we solve this problem using Approach C. \REVIEW{As can be seen in Fig.~\ref{fig_circular_annulus_ooa_approach_C}, Approach C also exhibits a very similar convergence rate as Approach A and B when the discontinuous indicator function is used, whereas the convergence rate is reduced when the continuous function is employed. Specifically, $\mathcal{O}(h^{0.91})$ (respectively, $\mathcal{O}(h^{1.22})$) convergence rate with an $R^2$ value of 0.79 (respectively, 0.82) in $L^{\infty}$ (respectively, $L^1$) norm is obtained using the continuous $\chi$.} However, Approach C is the most general one, since it can be used for imposing spatially varying flux values as demonstrated in later examples.  

The present example is also used to study the effect of the number of propagation cells $\nprop$ on the solution accuracy for the continuous indicator function. The results are shown in  Fig.~\ref{fig_circular_annulus_effect_of_prop_cells},  in which it can be observed that 2 cells on either side of the interface are sufficient for propagating $g$ values for a fixed number of $\ncells$ cells; the error norms are mostly affected by the $\ncells$ choice. Based on the results of this problem, we choose $\ncells = 1$ and $\nprop = 2$ for the continuous masking function, unless otherwise stated. For the discontinuous indicator function also, we use $\nprop = 2$  for the remainder of the problems (although $\nprop = 1$ is also sufficient).  
 
\subsection{Spatially varying flux values along complex interfaces} \label{sec_spatially_varying_flux}

\begin{figure}[]
  \centering
   \includegraphics[scale = 0.13]{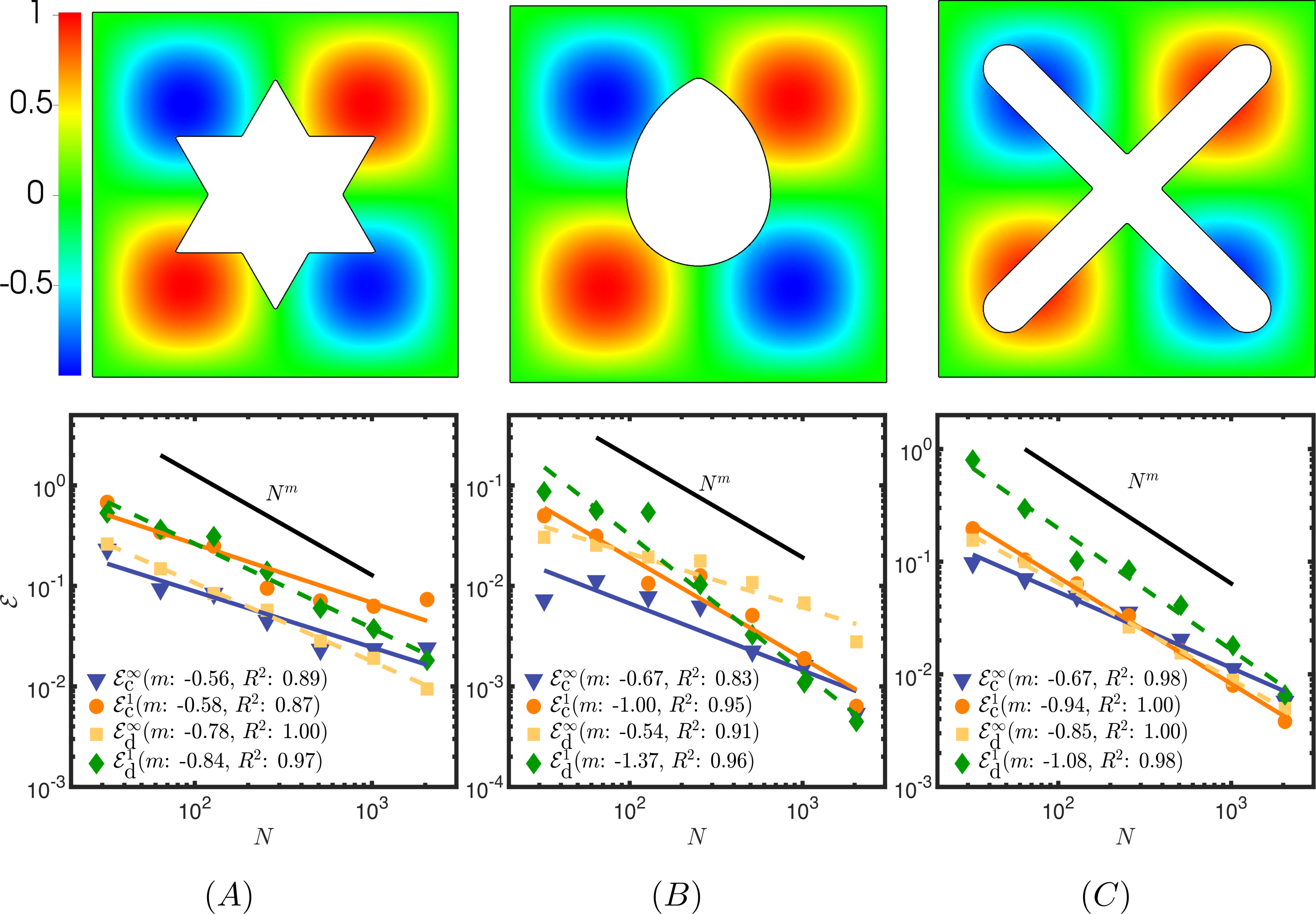}
   \caption{Numerical solution $q$ (top row) and the error norms \REVIEW{$\mathcal{E}^1$} and \REVIEW{$\mathcal{E}^\infty$} as a function of grid size $N$ (bottom row) for three complex shapes: (A) \REVIEW{Hexagram}; (B) \REVIEW{Egg}; and (C) \REVIEW{X-cross} using Approach C. The convergence rate results are shown for both continuous (solid line with symbols) and discontinuous (dashed line with symbols) indicator functions. The penalization parameter $\eta$ is taken as $10^{-8}$ and $\kappa$ is taken to be 1.}
     \label{fig_complex_annulus}
\end{figure}

In this section, we assess the accuracy of the numerical solution for spatially varying flux values using a manufactured solution of the form
\begin{equation}
\qexact(\x)= \sin(x)  \sin(y).
\label{eq_2dcomplex_analytical}
\end{equation}
Inhomogeneous Neumann boundary conditions $g(\x) = -\kappa\; \n \cdot \grad \qexact$ are imposed on the fluid-solid interface $\partial\Omegas$, whereas Dirichlet boundary conditions are imposed on the external boundaries of the computational domain, i.e., $\left.q\right|_{\partial \Omega(\x)} = \qexact(\x)$. Note that $g(\x)$ varies spatially, and therefore, Approach B is not applicable for this test. Eq.~\eqref{eq_2dcomplex_analytical} is plugged into the non-penalized Poisson Eq.~\eqref{eq_non_penalized_poisson} to generate the required source term $f(\x)$. We consider three geometrically complex solid domains: a hexagram, an egg, and a x-cross; these geometries are embedded in a larger Cartesian domain of extents $\Omega \in [0, 2\pi]^2$ and the numerical solutions are computed in the corresponding fluid domains $\Omegaf = \Omega \setminus \Omegas$.

Our prior work~\cite{Thirumalaisamy2021} considered Approach A for constructing the flux-forcing function for this problem, wherein $\qexact$ was used to define $\vbeta = \kappa\; \grad \qexact$. Although not feasible in practice (solution is unknown), Approach A results in second-order convergence rate of the numerical solution for this problem; see Appendix~\ref{sec_approach_a_results}. Next, we solve the same problem using Approach C. Fig.~\ref{fig_complex_annulus} presents the numerical solution and its convergence rate as a function of grid resolution. \REVIEW{For the hexagram case, at least $\mathcal{O}(h^{0.56})$ accuracy is achieved using the continuous indicator function whereas at least $\mathcal{O}(h^{0.78})$ accuracy is achieved using the discontinuous indicator function. Similarly, for the egg case, at least $\mathcal{O}(h^{0.67})$ accuracy is achieved using the continuous $\chi$ and at least $\mathcal{O}(h^{0.54})$ accuracy is obtained using the discontinuous indicator function. Lastly, for the x-cross geometry,  Approach C exhibits at least $\mathcal{O}(h^{0.67})$ accuracy using the continuous indicator function and at least $\mathcal{O}(h^{0.85})$ accuracy with the discontinuous one. As noted in the previous section also, Approach C with the discontinuous indicator function is able to achieve a better convergence rate than with the discontinuous one.} We remark that for Approach C, the reduction in accuracy (when compared to Approach A) is attributed to the codimension-0 extension of the spatially varying $g$ function in the neighborhood of the interface. Nevertheless, Approach C is able to impose spatially varying flux values on a complex interface (sharp corners, etc.) and the solution accuracy is also reasonable. \REVIEW{Later in Sec.~\ref{sec_rounded_hexagram} we demonstrate that smoothing of geometric features like sharp corners improves the accuracy of Approach C further.}

\subsection{Constant and spatially varying flux on three-dimensional interfaces} \label{sec_three_dimensional}
In this section, we consider two complex geometries in three spatial dimensions: a sphere and a torus. These geometries are embedded into a larger computational domain of extents $\Omega \in [0, 2\pi]^3$. 

For the spherical geometry, the  manufactured solution is taken to be
\begin{equation}
\qexact(r)= r^2 + c,
\label{eq_qsphere_analytical}
\end{equation}
in which $r = \sqrt{\left(x - \pi\right)^2 + \left(y - \pi\right)^2 + \left(z - \pi\right)^2}$ and $c$ is a constant.  Eq.~\eqref{eq_qsphere_analytical} when plugged into the non-penalized Poisson Eq.~\eqref{eq_non_penalized_poisson} yields a constant source term $f(\x) = -6$. The radius of the sphere is taken to be $R = 3/2$. Two cases are considered for the spherical geometry: fluid inside the sphere and fluid outside it. For the first case, a constant flux boundary condition $g = -\kappa\; \partial \qexact/\partial n = -2R$ is imposed on the spherical surface and a homogeneous Dirichlet boundary condition is imposed on the external domain boundary $\partial \Omega$.  Since the solution of this Poisson problem is defined up to an additive constant $c$, we use the zero-mean condition $\int_0^R\; 4\pi r^2 \; \qexact \; \Dr = 0$ to determine the constant $c = \frac{-3 R^2}{5}$. For the second case, in which the fluid is considered between the spherical interface and the computational domain boundary, constant flux boundary condition $g = -\kappa \partial \qexact/\partial n = -2R$ is imposed on the spherical interface and inhomogeneous Dirichlet boundary conditions $\left.q\right|_{\partial \Omega(\x)} = \qexact(\x)$ are imposed on the external boundary. The constant $c$ is taken to be zero for this case. 

\begin{figure}[]
  \centering
   \includegraphics[scale = 0.13]{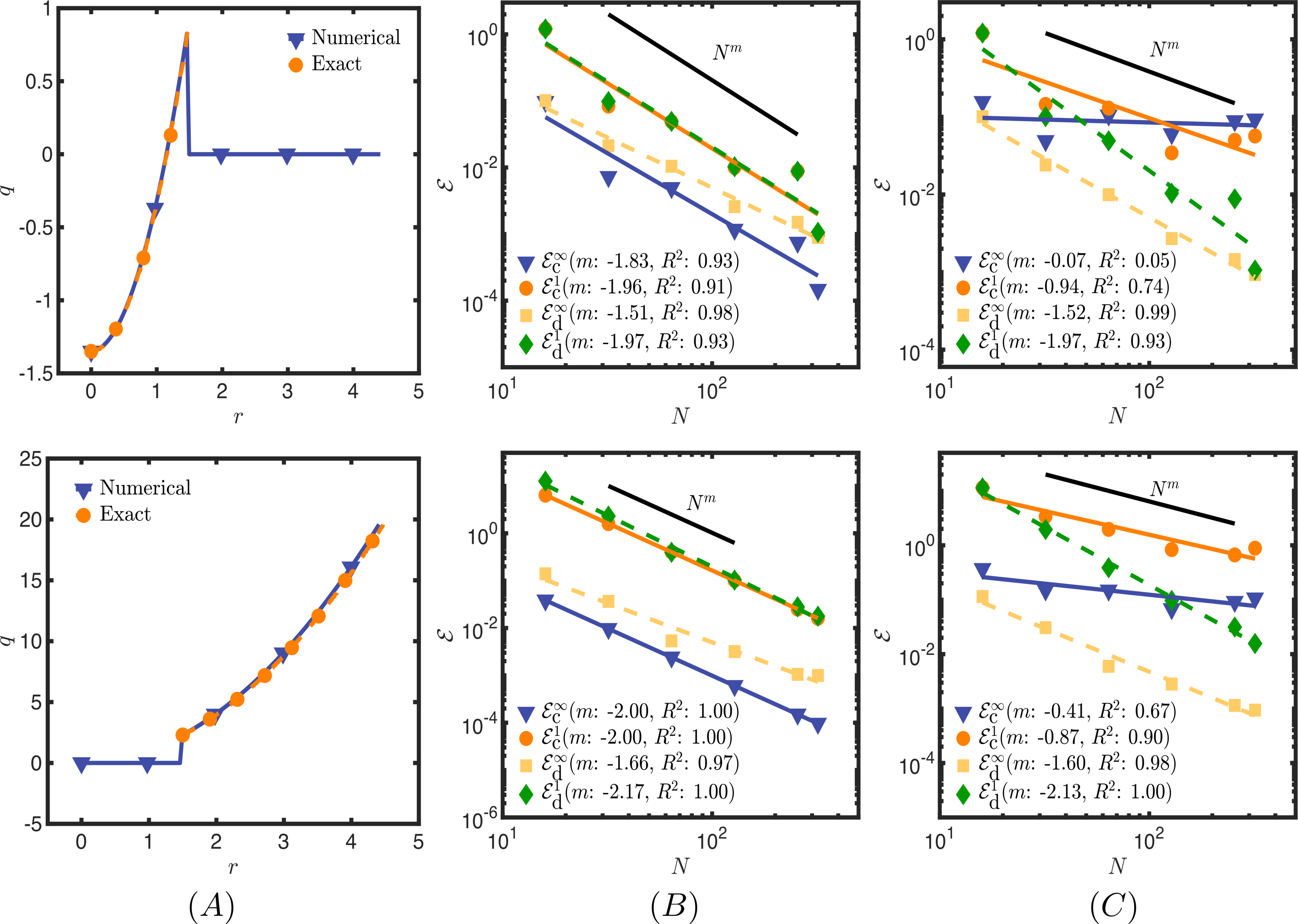}
   \caption{Spherical interface with constant flux boundary condition using Approach B and C. Top row corresponds to fluid inside the sphere case and the bottom row corresponds to fluid outside the sphere case. Also shown are the error norms, \REVIEW{$\mathcal{E}^1$} and \REVIEW{$\mathcal{E}^\infty$} as a function of grid size $N$ using continuous (solid lines with symbols) and discontinuous (dashed lines with symbols) indicator functions. (A) Comparison of numerical solution at $z = \pi$ using $N = 256$ grid. Error norms using (B) \REVIEW{Approach B}; and (C)  \REVIEW{Approach C}. The penalization parameter $\eta$ is taken as $10^{-8}$, $\kappa$ as 1 and the radius of the sphere is 3/2.}
     \label{fig_sphere}
\end{figure}

The comparison between the numerical and analytical solutions, as well as the spatial convergence rate of $\mathcal{E}^1$ and $\mathcal{E}^\infty$ error norms using Approach B and C are shown in Fig.~\ref{fig_sphere}. As can be observed in the figure, the numerical solution is in excellent agreement with the exact solution. Largely \REVIEW{$\mathcal{O}(h^2)$} convergence rate is obtained for this example using Approach B with both continuous and discontinuous indicator functions. \REVIEW{Approach C also yields the same order of accuracy with the discontinuous indicator function as Approach B. However, the convergence rate using the continuous indicator function is between 0 and 1 for Approach C. Clearly, the discontinuous indicator function performs better than the continuous one for Approach C.} Second-order convergence rate is also obtained with Approach A, when $\vbeta = \kappa \; \grad \qexact$ is used for the spherical geometry (data not shown for brevity). \REVIEW{Better performance of Approach A compared to Approach B and C is expected, as mentioned in the beginning of Sec.~\ref{sec_results_and_discussion}. }

\begin{figure}[]
\centering
\subfigure[{Toroidal interface}]{
\includegraphics[scale = 0.15]{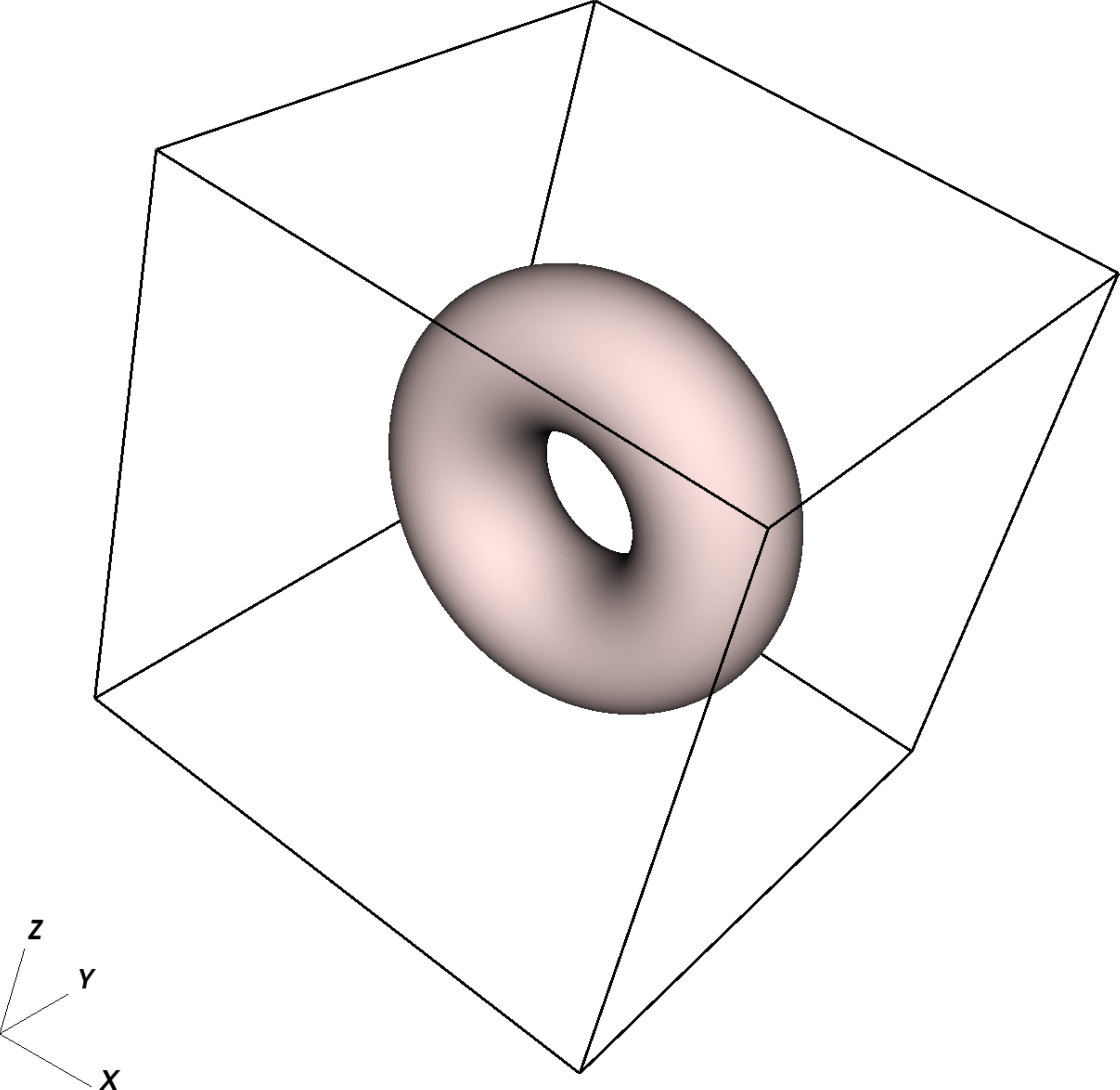}
\label{fig_torus_zero_contour}
}
\subfigure[{Numerical solution at $y = \pi$}]{
\includegraphics[scale = 0.15]{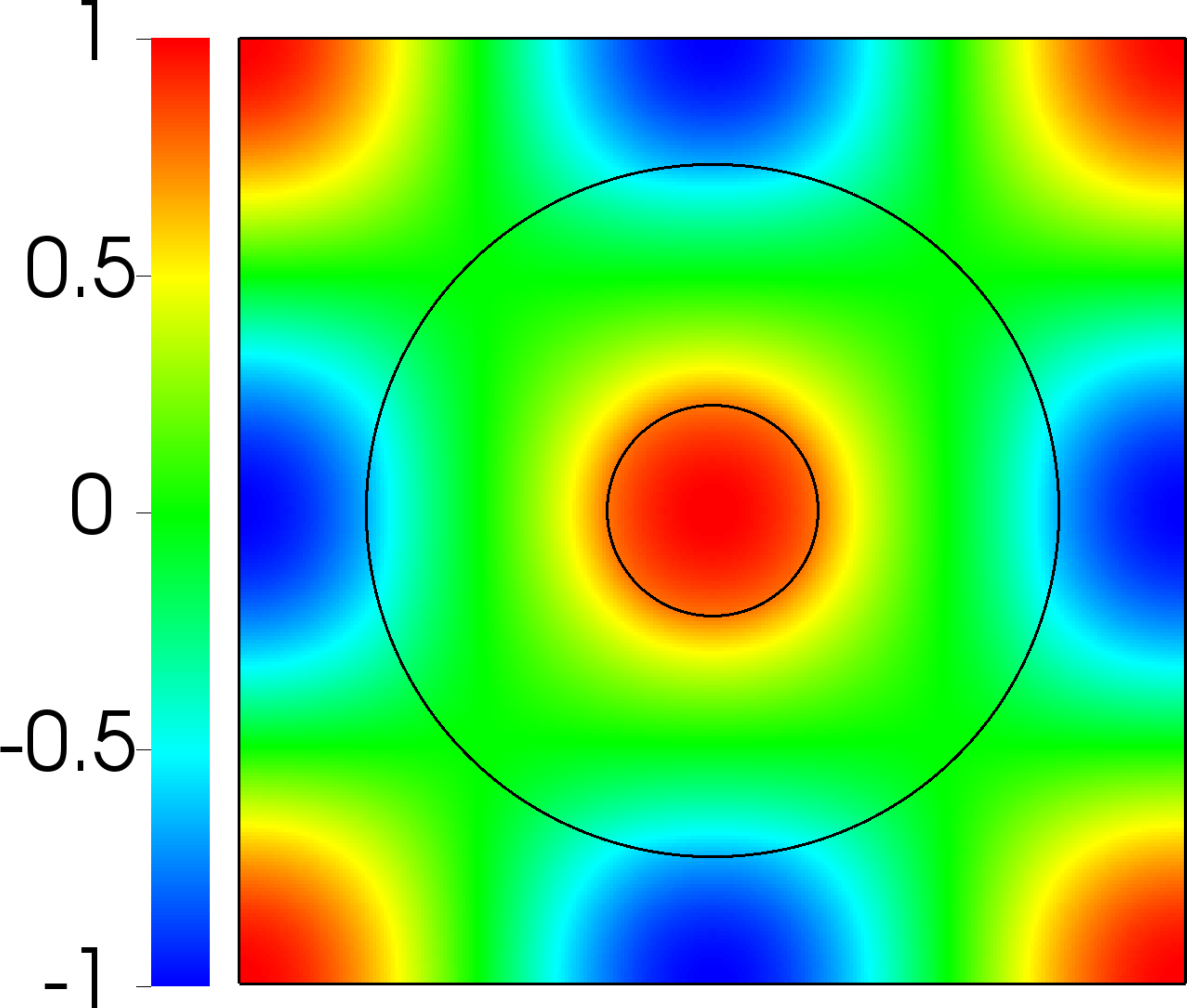}
\label{fig_torus_y_plane}
}
\subfigure[{\REVIEW{Spatial convergence rate for Approach C}}]{
\includegraphics[scale = 0.08]{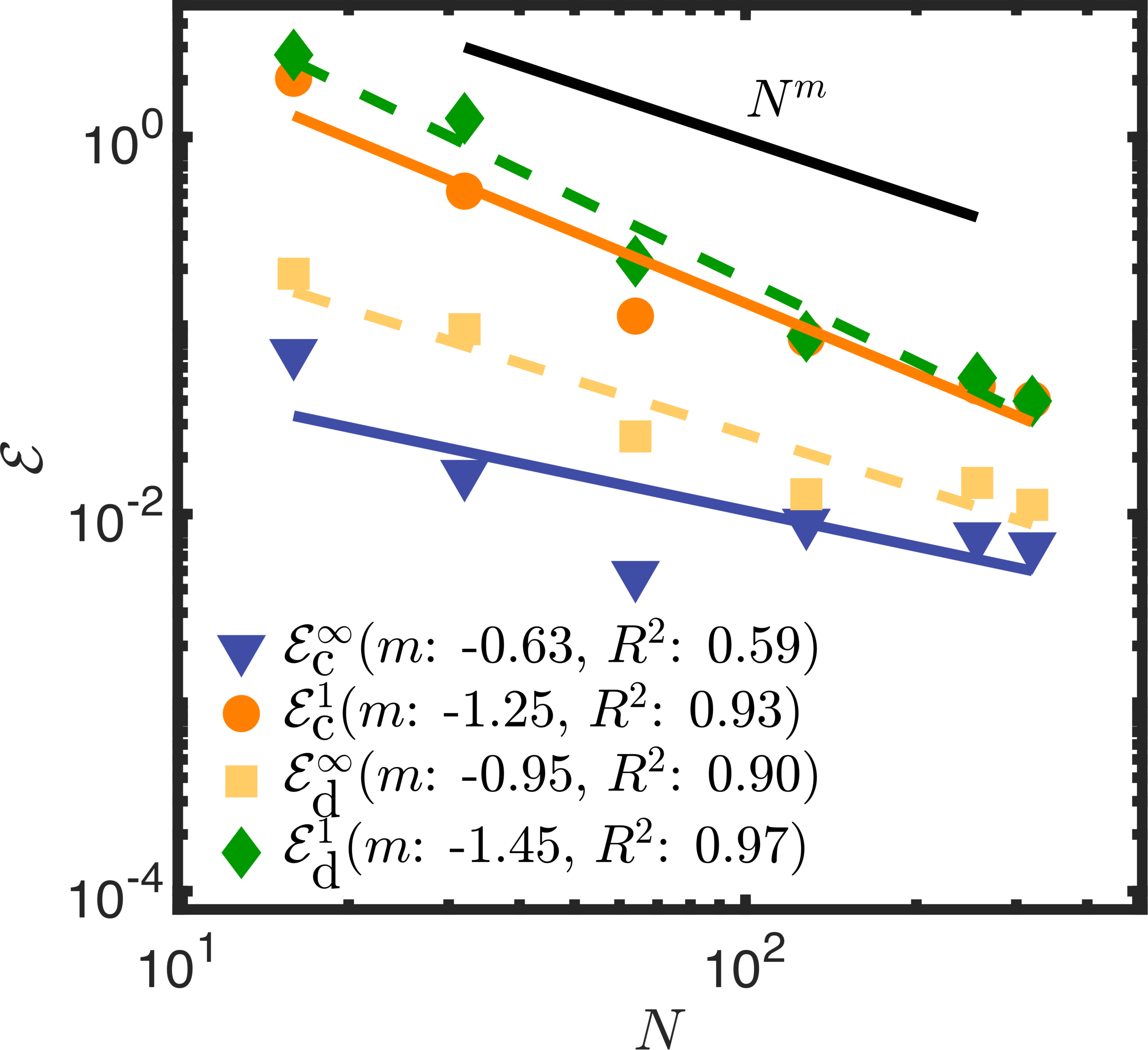}
\label{fig_torus_ooa_approach_C}
}
 \caption{Toroidal interface with spatially varying flux boundary conditions using Approach C.~\subref{fig_torus_zero_contour} Zero-contour of the solid torus;~\subref{fig_torus_y_plane} numerical solution at $y = \pi$ using $N = 256$ grid; and~\subref{fig_torus_ooa_approach_C} the \REVIEW{error norms \REVIEW{$\mathcal{E}^1$} and \REVIEW{$\mathcal{E}^\infty$} as a function of grid size $N$ using the continuous (solid lines with symbols) and discontinuous (dashed lines with symbols) indicator functions.} The penalization parameter $\eta$ is taken as $10^{-8}$ and $\kappa$ is taken to be 1.}
\label{fig_torus}
\end{figure}

For the next three-dimensional test example, a solid torus is embedded in a computational domain of extents $\Omega \in [0, 2\pi]^3$ as shown in Fig.~\ref{fig_torus_zero_contour} and the fluid region is taken outside of the torus. We consider a manufactured solution of the form
\begin{equation}
\qexact(\x)= -\cos(x)  \cos(y) \cos(z).
\label{eq_torus_analytical}
\end{equation}
Eq.~\eqref{eq_torus_analytical} is plugged into the non-penalized Poisson Eq.~\eqref{eq_non_penalized_poisson} to generate the required source term $f(\x)$. On the toroidal interface, spatially varying inhomogeneous Neumann boundary conditions are imposed, whereas on the external domain boundary inhomogeneous Dirichlet boundary conditions using the exact solution are imposed. We solve this test problem using Approach C and the results are shown in Fig.~\ref{fig_torus}.  \REVIEW{As can be observed in Fig.~\ref{fig_torus}, at least $\mathcal{O}(h^{1.25})$ is achieved with the discontinuous indicator function and at least $\mathcal{O}(h^{0.63})$ is achieved with the continuous indicator function.}

\subsection{Spatially constant Robin boundary condition on two-dimensional interfaces} \label{sec_spatially_constant_robin}
We consider the concentric circular annulus problem of Sec.~\ref{sec_circular_annulus_constant_flux} with the same exact solution $\qexact(\r)$ and source term $f(r)$, as written in Eqs.~\eqref{sakurai_2d_poisson_exact} and~\eqref{sakurai_annulus_src_term}, respectively. Plugging the exact solution into the Robin boundary condition Eq.~\ref{eqn_robin} yields a spatially constant $g$ value for the inner and outer interface, respectively. 

\begin{figure}[]
\centering
\subfigure[{Numerical solution}]{
\includegraphics[scale = 0.08]{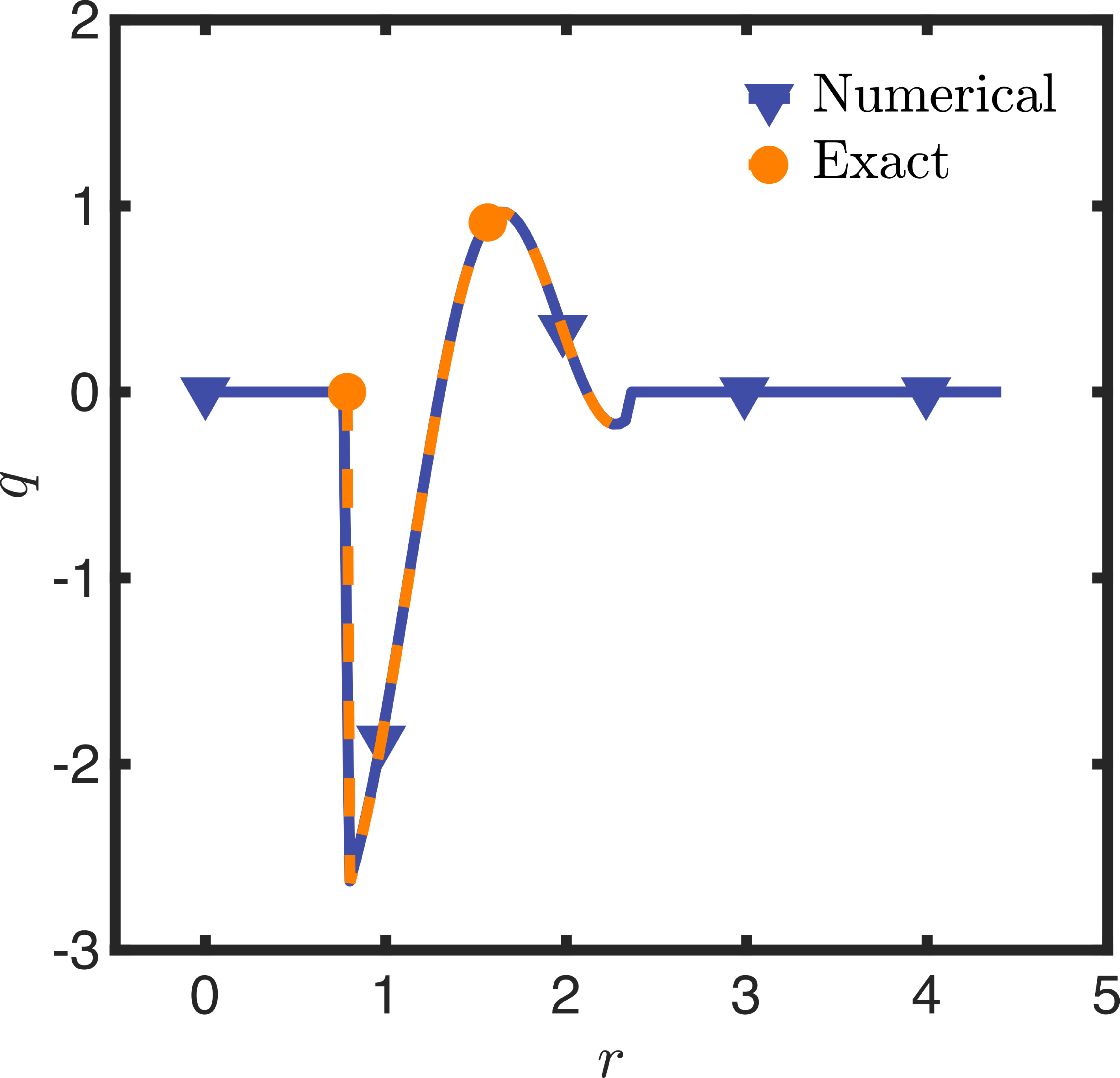}
\label{fig_robin_sakurai_case_solution}
}

\subfigure[{\REVIEW{Order of convergence using Approach B}}]{
\includegraphics[scale = 0.08]{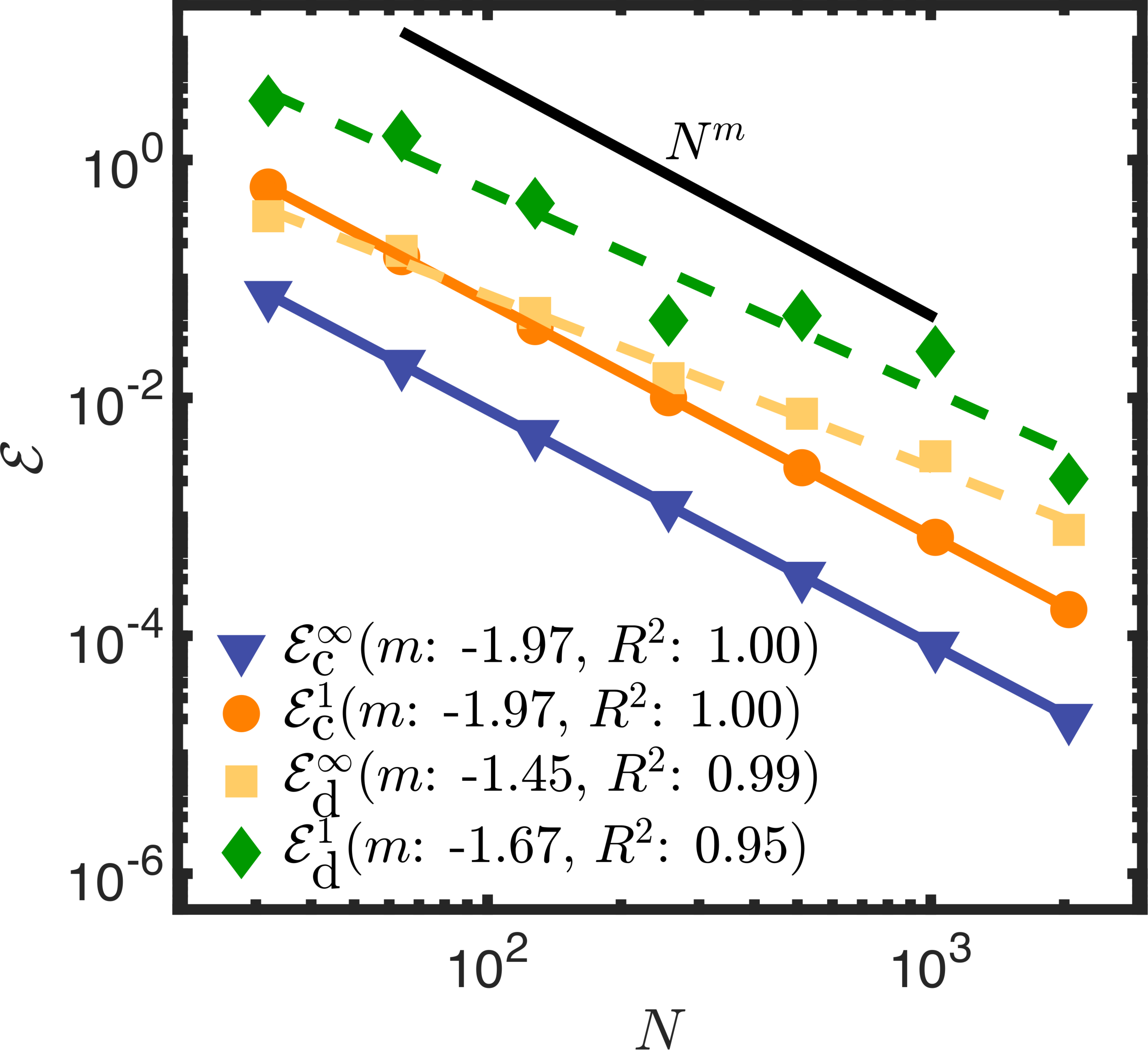}
\label{fig_robin_sakurai_case_approach_B}
}
\subfigure[{\REVIEW{Order of convergence using Approach C}}]{
\includegraphics[scale = 0.08]{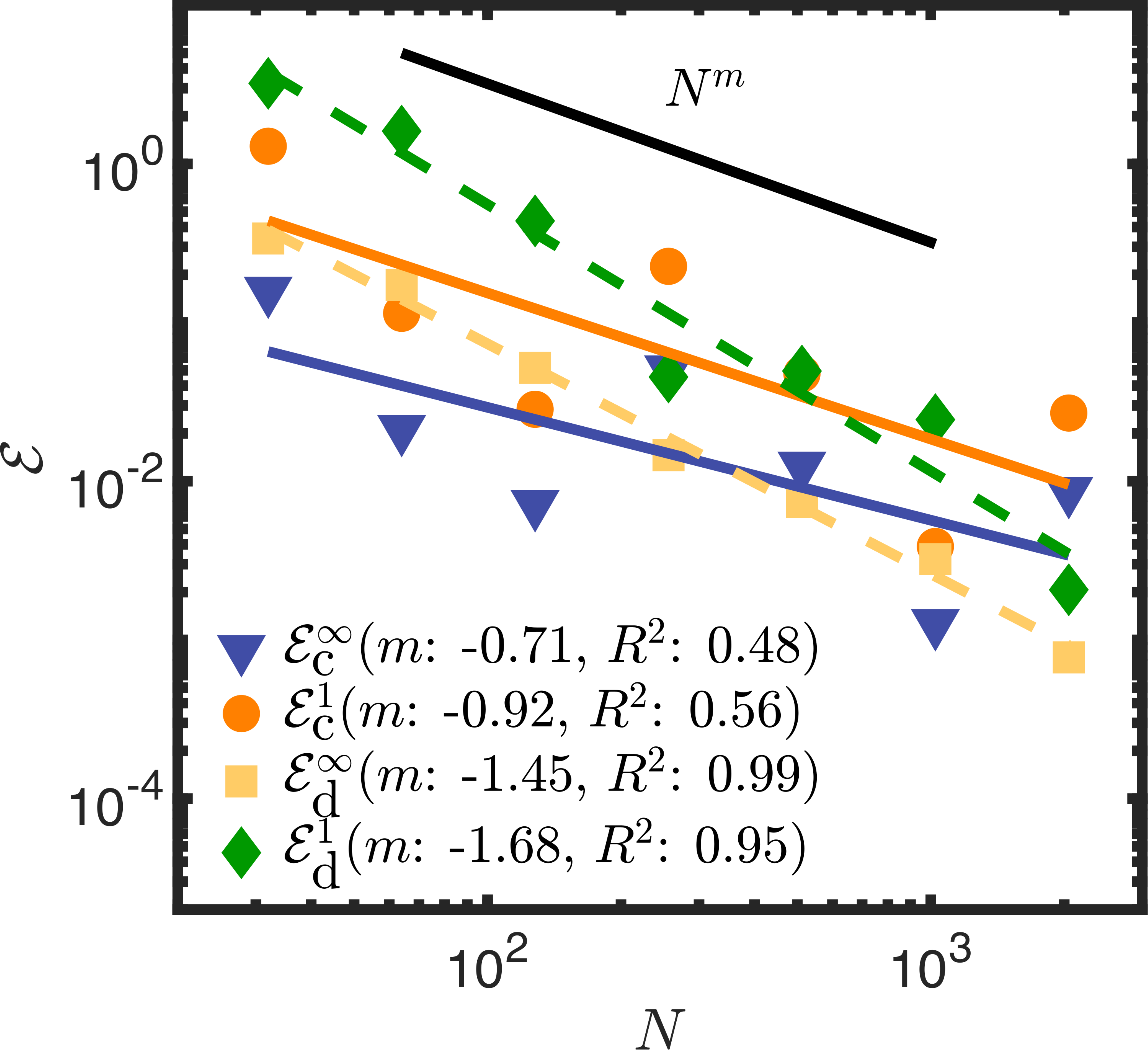}
\label{fig_robin_sakurai_case_approach_C}
}
 \caption{ Concentric annulus with spatially constant Robin boundary conditions using Approach B and C.~\subref{fig_robin_sakurai_case_solution} Numerical solution using $N = 256$ grid. Error norms \REVIEW{$\mathcal{E}^1$} and \REVIEW{$\mathcal{E}^\infty$} as a function of grid size $N$ using the continuous (solid lines with symbols) and discontinuous (dashed lines with symbols) indicator functions for~\subref{fig_robin_sakurai_case_approach_B} \REVIEW{Approach B}; and~\subref{fig_robin_sakurai_case_approach_C}  \REVIEW{Approach C}. The penalization parameter $\eta$ is taken as $10^{-8}$. The values of $\kappa$ and $\zeta$ are taken to be 1.}
\label{fig_robin_sakurai_case}
\end{figure}

The VP Poisson Eq.~\ref{eqn_vp_robin} is solved using Approach B and C for this problem. The numerical solution compared against Eq.~\ref{sakurai_2d_poisson_exact} using Approach C is shown in Fig.~\ref{fig_robin_sakurai_case_solution}; an excellent agreement is obtained. We also present the convergence rate for Approach B and  C in Fig.~\ref{fig_robin_sakurai_case}.  \REVIEW{As can be seen in Fig.~\ref{fig_robin_sakurai_case}, the convergence rates obtained by using the discontinuous indicator function for Approach B and C are quite close to what we had obtained in Sec.~\ref{sec_circular_annulus_constant_flux}. For the continuous indicator function, we obtain approximately second-order accuracy with Approach B and at least $\mathcal{O}(h^{0.71})$ accuracy with Approach C. As also observed in Sec.~\ref{sec_circular_annulus_constant_flux}, the discontinuous indicator function performs better than the continuous one with Approach C; the reverse is true for Approach B.}

\subsection{Spatially varying Robin boundary condition on a complex two-dimensional interface}\label{sec_spatially_varying_robin}

\begin{figure}[]
\centering
\subfigure[{Numerical solution}]{
\includegraphics[scale = 0.15]{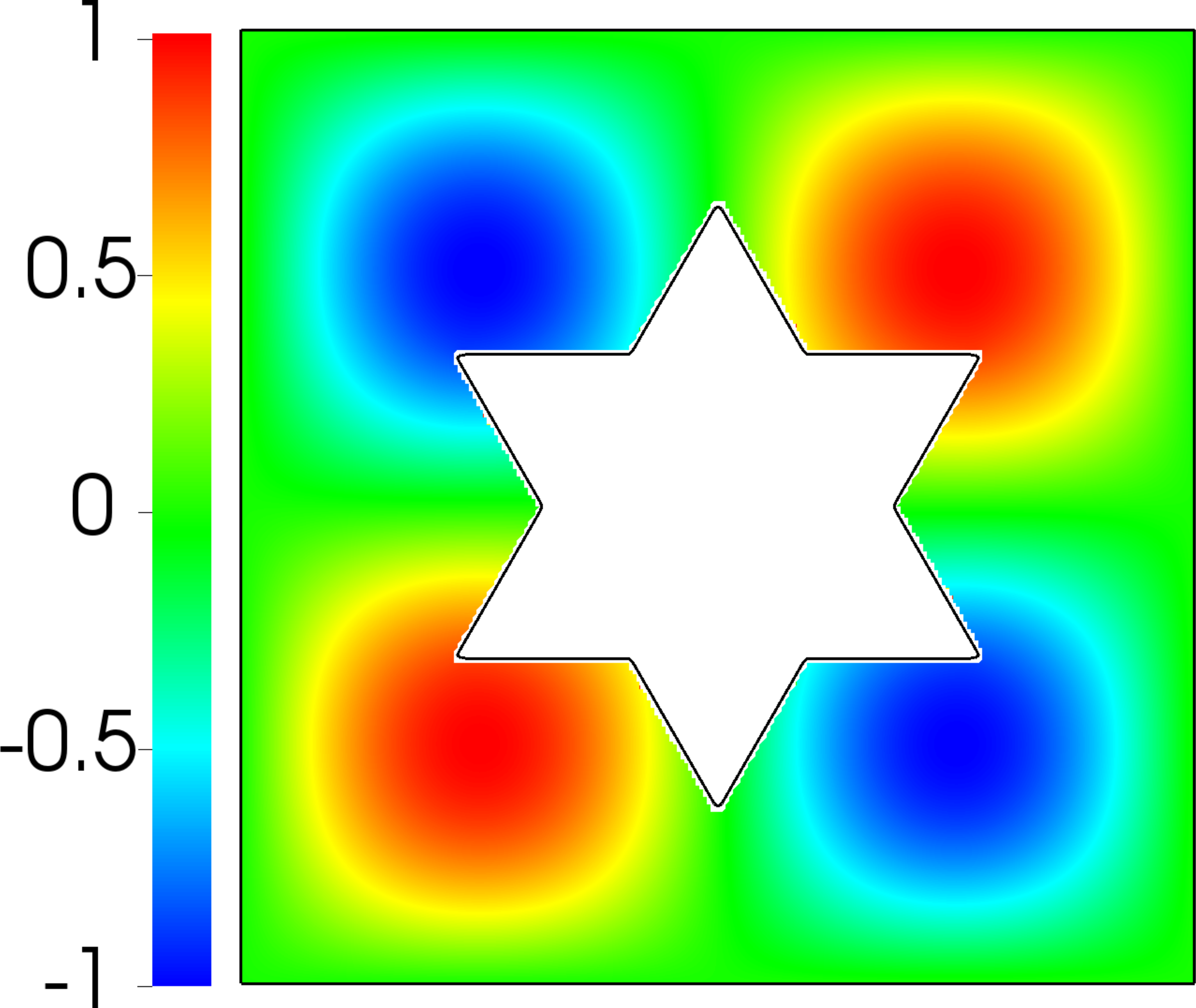}
\label{fig_robin_hexagram_solution}
}
\subfigure[{\REVIEW{Order of convergence}}]{
\includegraphics[scale = 0.08]{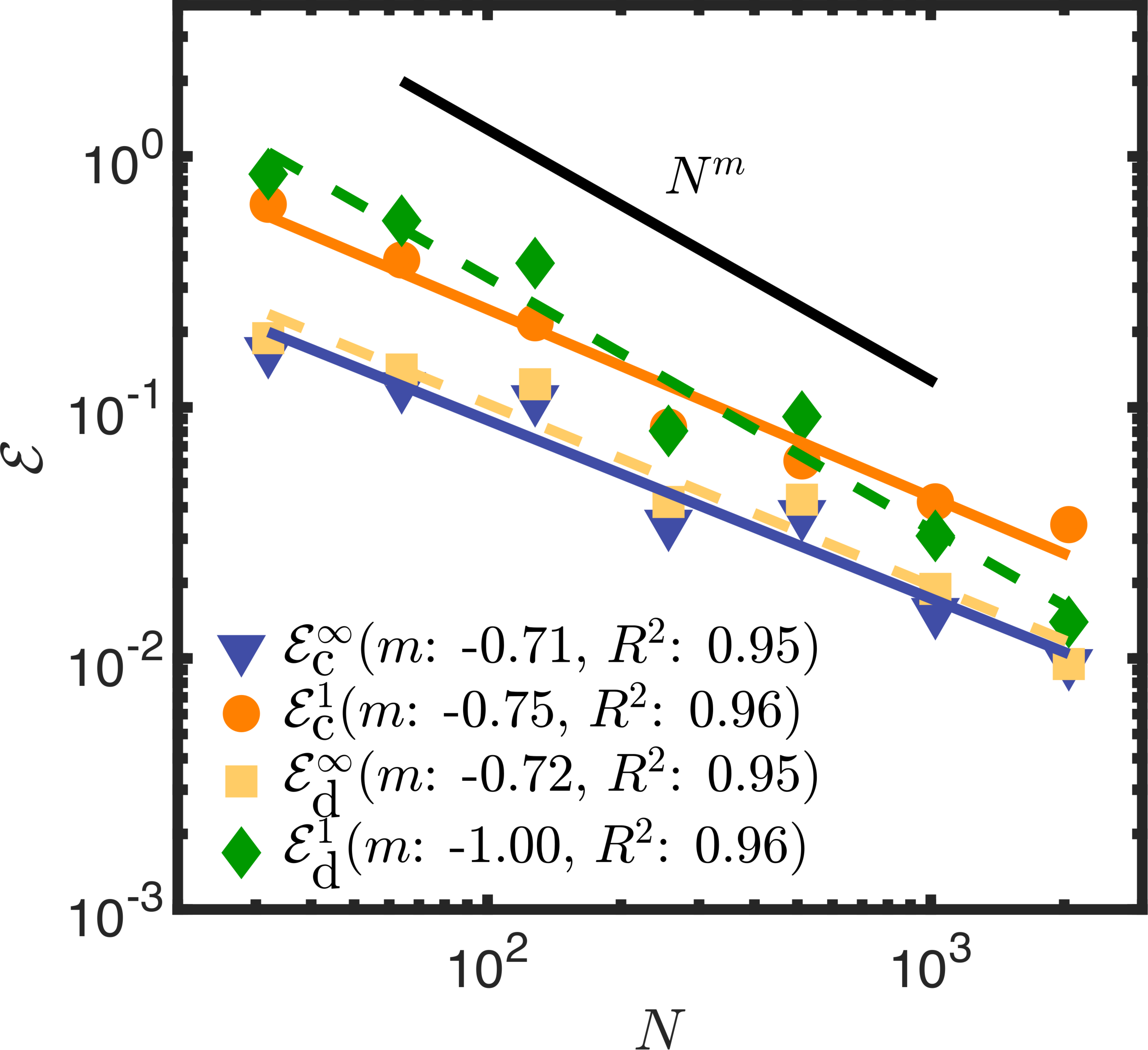}
\label{fig_robin_hexagram_ooc}
}
 \caption{Hexagram domain with spatially varying Robin boundary conditions using Approach C.~\subref{fig_robin_hexagram_solution} Numerical solution using $N = 256$ grid; and~\subref{fig_robin_hexagram_ooc} the \REVIEW{error norms \REVIEW{$\mathcal{E}^1$} and \REVIEW{$\mathcal{E}^\infty$} as a function of grid size $N$ using the continuous (solid lines with symbols) and discontinuous (dashed lines with symbols) indicator functions}. The penalization parameter $\eta$ is taken as $10^{-8}$. The values of $\kappa$ and $\zeta$ are taken to be 1.}
\label{fig_robin_hexagram}
\end{figure}

In this section, we assess the accuracy of Approach C for spatially varying Robin boundary conditions on a complex two-dimensional interface. A hexagram geometry is embedded into a computational domain of extents $\Omega \in [0, 2\pi]^2$, as considered in Sec.~\ref{sec_spatially_varying_flux}. The fluid is occupied between the computational domain boundary $\partial \Omega$ and the fluid-solid interface $\partial \Omega_s$. The same manufactured solution as written in Eq.~\ref{eq_2dcomplex_analytical} is considered here; this solution yields spatially varying $g$ values when plugged into the Robin boundary condition Eq.~\eqref{eqn_robin}. We solve the VP Poisson Eq.~\ref{eqn_vp_robin} using Approach C. The numerical solution and the spatial convergence rate of the error norms are presented in Fig.~\ref{fig_robin_hexagram}. \REVIEW{As observed in the figure, at least $\mathcal{O}(h^{0.71})$ accuracy is achieved with the continuous indicator function and at least $\mathcal{O}(h^{0.72})$ accuracy is achieved with the discontinuous indicator function. } 

\REVIEW{\subsection{Effect of smoothing geometric features on the convergence rate of Approach C} \label{sec_rounded_hexagram}
In this section, we study the effect of sharp geometric features, such as corners on the convergence rate of Approach C. We consider a slightly modified version of the hexagram interface, which was considered earlier in Secs.~\ref{sec_spatially_varying_flux} and~\ref{sec_spatially_varying_robin}. Instead of the sharp corners, in this case, we embed a hexagram having smooth exterior corners (see Fig.~\ref{fig_rounded_hexagram_solution}) in a larger Cartesian domain of extents $\Omega \in [0, 2\pi]^2$ and solve the Neumann/Robin problem using Approach C.  The order of accuracy results for spatially varying Neumann and Robin boundary conditions using the discontinuous indicator function are presented  in Fig.~\ref{fig_rounded_hexagram}. For the sharp corner geometry case with spatially varying Neumann boundary conditions (as shown in Fig~\ref{fig_complex_annulus}(A)), the convergence rates were $\mathcal{O}(h^{0.78})$ and $\mathcal{O}(h^{0.84})$ in $L^\infty$ and $L^1$ norm, respectively. In contrast, with smooth corners, the convergence rates are $\mathcal{O}(h^{0.95})$ and $\mathcal{O}(h^{1.08})$ in $L^\infty$ and $L^1$ norm, respectively. A similar trend is obtained when spatially varying Robin boundary conditions are considered: For the sharp geometry case (as shown in  Fig.~\ref{fig_robin_hexagram}), the convergence rates were $\mathcal{O}(h^{0.71})$ and $\mathcal{O}(h^{1.00})$ in $L^\infty$ and $L^1$ norm, respectively. With smooth corners, the convergence rates are $\mathcal{O}(h^{1.00})$ and $\mathcal{O}(h^{1.26})$ in $L^\infty$ and $L^1$ norm, respectively. This test  demonstrates that the convergence rate of Approach C also depends upon local geometric features.

\begin{figure}[]
\centering
\subfigure[{\REVIEW{Smoothed hexagram geometry}}]{
\includegraphics[scale = 0.137]{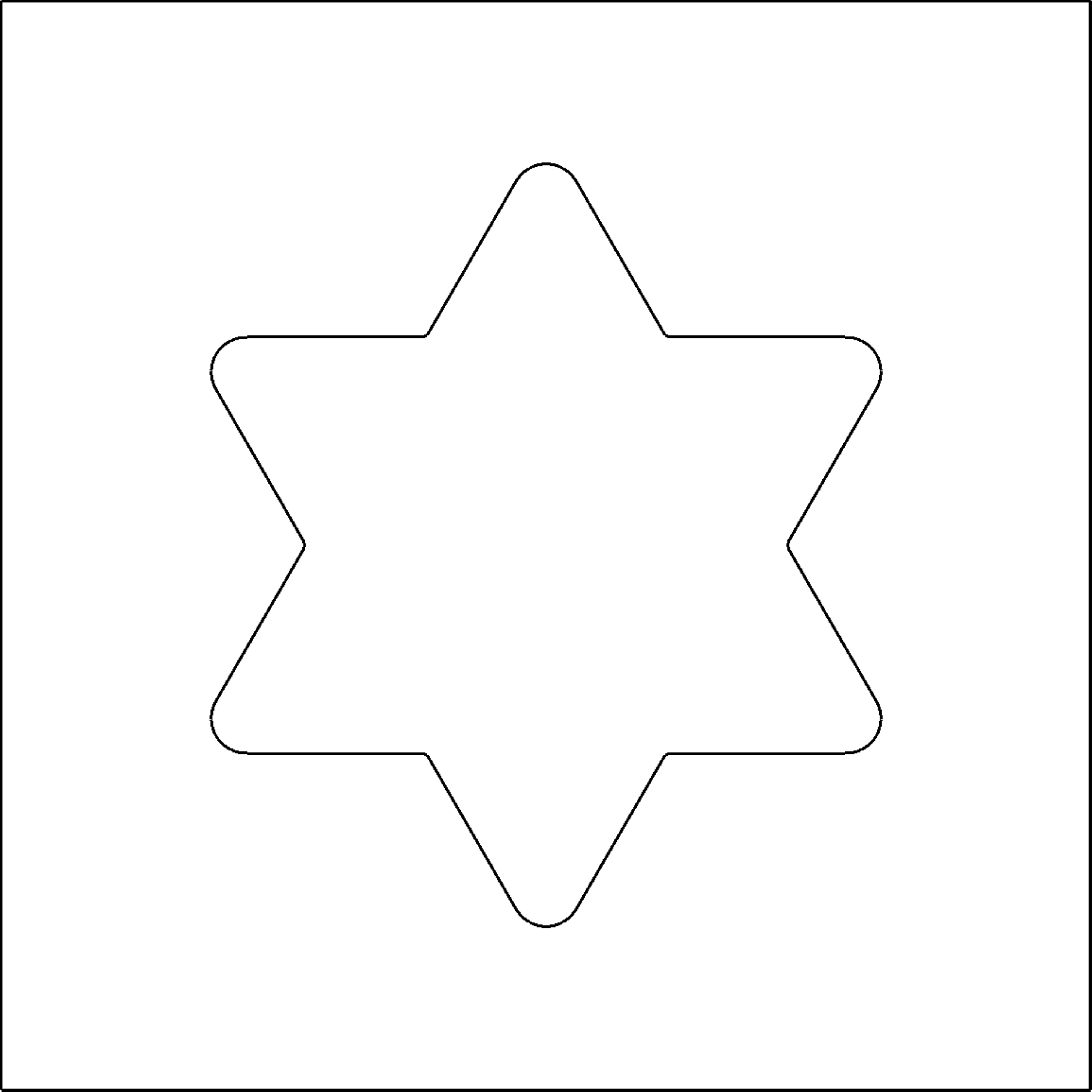}
\label{fig_rounded_hexagram_solution}
}
\subfigure[{\REVIEW{Order of convergence for spatially varying Neumann boundary condition}}]{
\includegraphics[scale = 0.08]{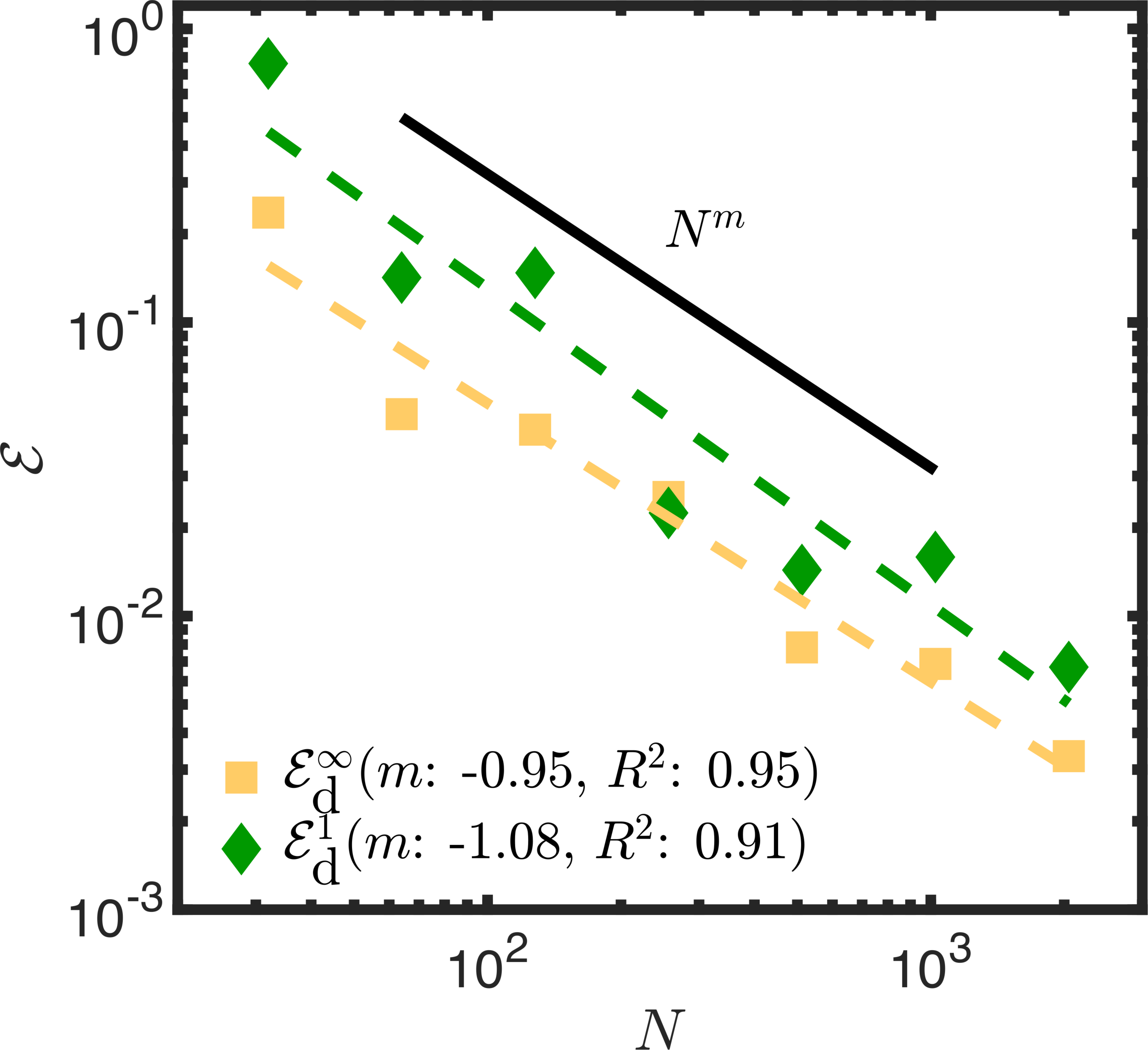}
\label{fig_rounded_hexagram_neumann_ooc}
}
\subfigure[{\REVIEW{Order of convergence for spatially varying Robin boundary condition}}]{
\includegraphics[scale = 0.08]{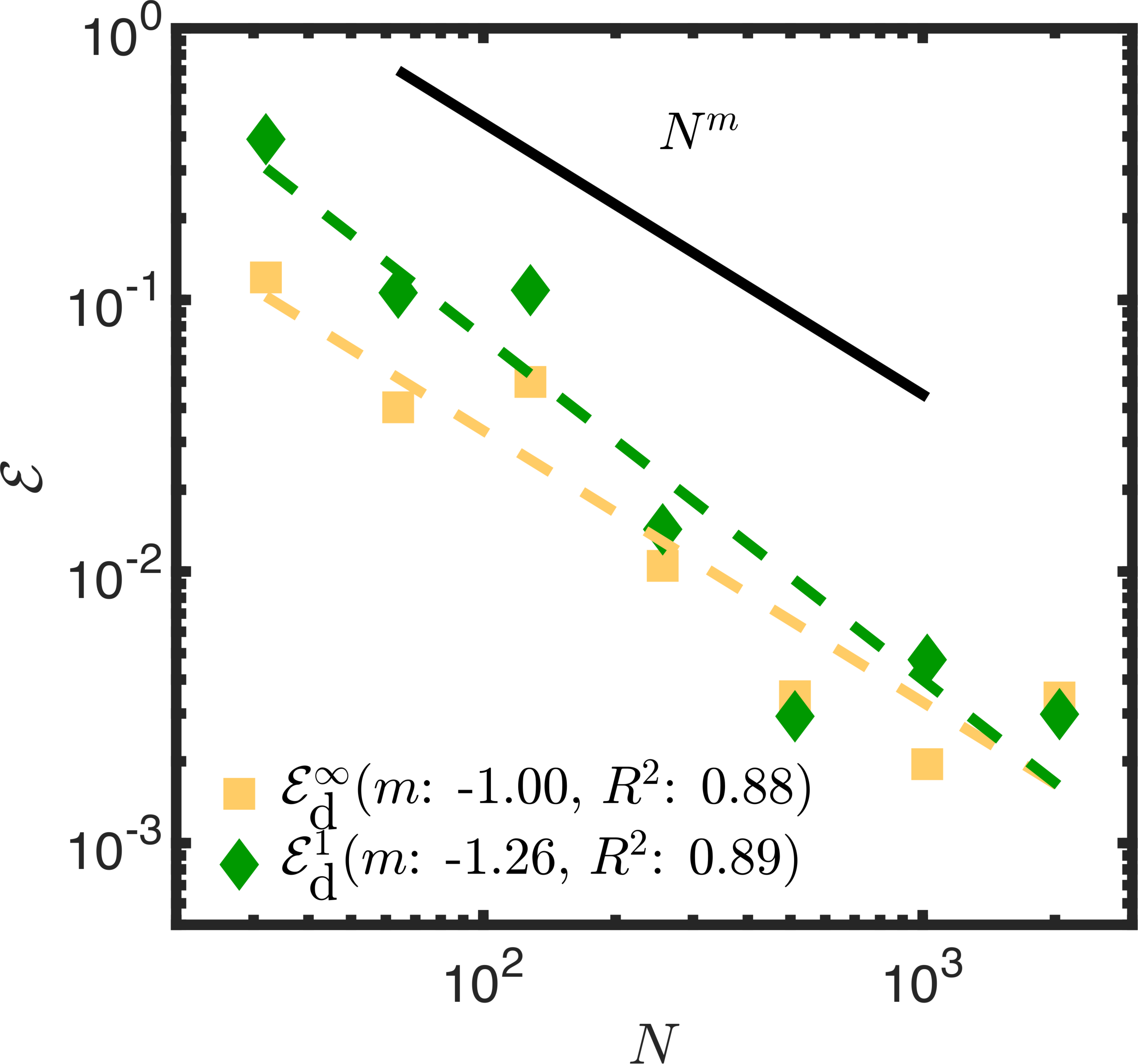}
\label{fig_rounded_hexagram_robin_ooc}
}
 \caption{\REVIEW{Hexagram domain with smooth exterior corners.~\subref{fig_rounded_hexagram_solution} Zero-contour of the smoothed hexagram interface. Error norms \REVIEW{$\mathcal{E}^1$} and \REVIEW{$\mathcal{E}^\infty$} as a function of grid size $N$ using the discontinuous indicator function with Approach C.~\subref{fig_rounded_hexagram_neumann_ooc} Spatially varying Neumann boundary conditions; and~\subref{fig_rounded_hexagram_robin_ooc} spatially varying Robin boundary conditions. The penalization parameter $\eta$ is taken as $10^{-8}$. The values of $\kappa$ and $\zeta$ are taken to be 1.}}
\label{fig_rounded_hexagram}
\end{figure}
}

\subsection{Application to free convection problem}\label{sec_free_conv}

\begin{figure}[]
  \centering
   \includegraphics[scale = 0.13]{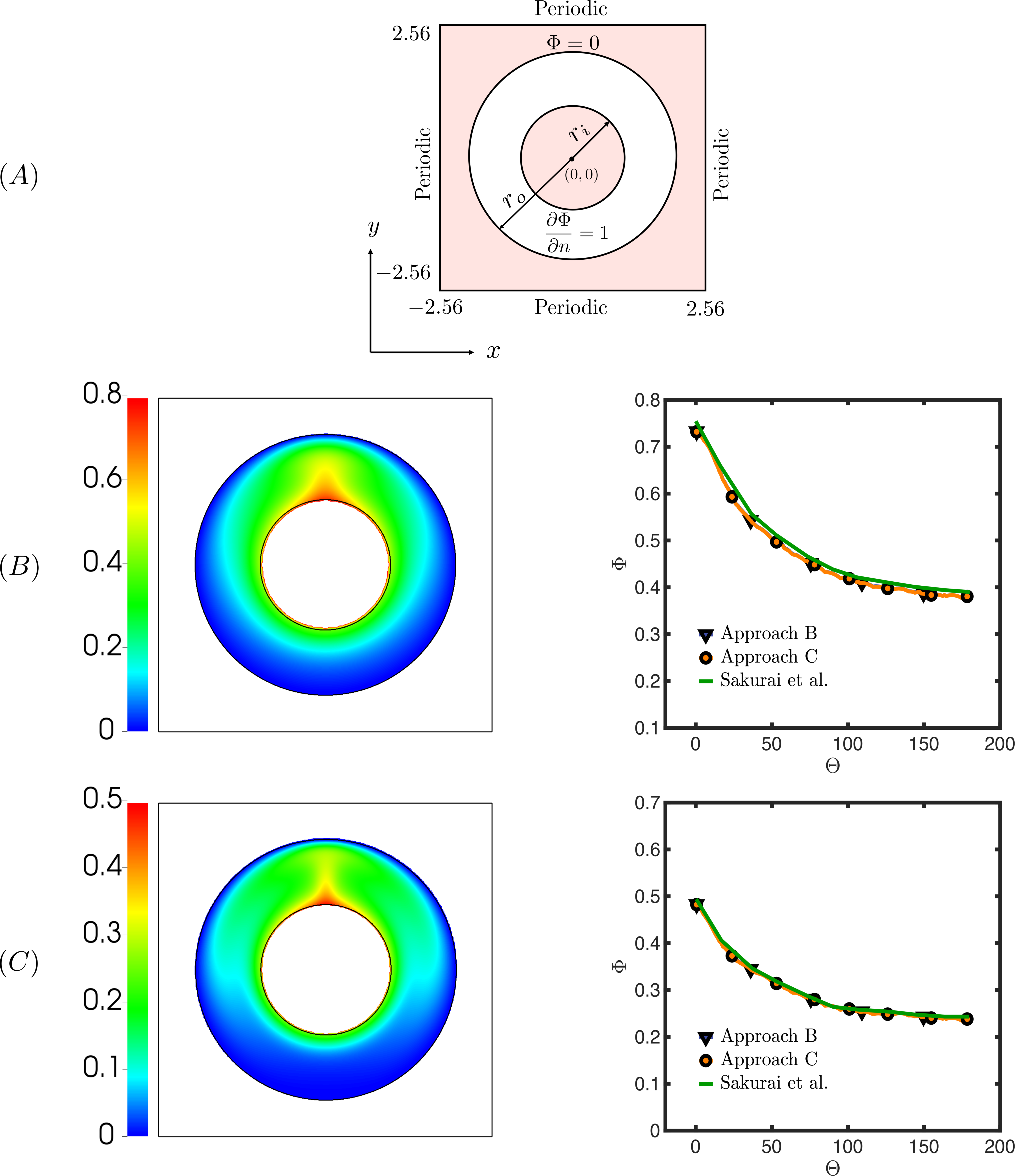}
   \caption{2D free convection problem at two Rayleigh numbers: Ra = 5700 and Ra = $5 \times 10^4$.(A) Schematic of the problem. Steady state temperature field and temperature distribution on the surface of the left half of the inner cylinder at (B) Ra $= 5700$ and (C) Ra = $5\times 10^4$.
}
     \label{fig_free_conv}
\end{figure}

Finally, we consider steady natural convection in a concentric annulus. A constant heat flux $Q$ is imposed on the inner cylinder of radius $r_i$ and a fixed temperature $T_o$ is maintained on the outer cylinder of radius $r_o$. The concentric annulus is embedded into a larger computational domain of extents $\Omega \in [-2.56,2.56]^2$,  as shown in Fig.~\ref{fig_free_conv}. This example was studied in Yoo~\cite{Yoo2003} using a body-fitted grid approach and more recently, it has been used to validate the IB/FCT relying on time-splitting approach to handle the flux boundary conditions on embedded interfaces~\cite{Ren2013, Wang2016, Guo2019}; see Introduction Sec.~\ref{sec_intro} for a brief discussion on IB/FCT.  
 
We solve the volume penalized advection-diffusion equation for the temperature field coupled to the volume penalized incompressible Navier-Stokes equations (Eqs.~\eqref{eqn_momentum}-\eqref{eqn_adv_diff}) in a non-dimensional form. The non-dimensional quantities are defined as: dimensionless temperature $\Phi^* =  k(T - T_0)/(QL)$, velocity $\u^* = \u L/\alpha$, time $t^* = t\alpha/L^2$, and position $\x^* = \x/L$. Here, $k$ is the thermal conductivity, $\alpha$ is the thermal diffusivity $\alpha = k/(\rho c_p)$, $\rho$ is the density, $c_p$ is the heat capacity at constant pressure, and $L = r_o - r_i$ is the annulus thickness. For this case we consider $r_i^* = r_i/L = 1$ and $r_o^* = r_o/L = 2$. Dropping the $^*$ superscript from the non-dimensional quantities, the system of non-dimensional equations reads as
\begin{align}
\D{\u}{t} + \div (\u\u) &= -\grad p + \text{Pr} \grad^2 \u  + [1 -  (\chi^{\rm d} + \chi^{\rm n})] \text{RaPr} \Phi \V{e}_y - \frac{\chi^{\rm d} + \chi^{\rm n}}{\eta_{\rm d}} \u, \label{eqn_ndim_mom} \\
\div \u &= 0, \label{eqn_ndim_cont} \\
\D{\Phi}{t} + \left(1 - \chi^{\rm n} \right) \left(\u\cdot\grad \Phi\right) &= \div \left[\left\{\left(1 - \chi^{\rm n} \right) + \eta_{\rm n} \chi^{\rm n} \right\} \grad\Phi \right] + \div \left(\chi^{\rm n} \vbeta\right) - \chi^{\rm n} \div \vbeta - \frac{\chi^{\rm d}}{\eta_{\rm d}} \Phi. \label{eqn_ndim_adv_diff}
\end{align} 
Here, $\V{e}_y = (0,1)$ is a unit vector in the $y-$direction and the continuous indicator functions $\chi^{\rm d}$ and $\chi^{\rm n}$ are defined to be
\begin{align}
\label{eqn_chi_freeconv}
\chi^{\rm d}&=\begin{cases}
       1,  & \phi^{\rm d}(\x) < - \nsmear \; h,\\
       1 - \frac{1}{2}\left(1 + \frac{1}{\nsmear h} \phi^{\rm d}(\x)+ \frac{1}{\pi} \sin\left(\frac{\pi}{ \nsmear h} \phi^{\rm d}(\x)\right)\right) ,  & |\phi^{\rm d}(\x)| \le \nsmear\; h,\\
        0,  & \textrm{otherwise},
\end{cases} \\ \chi^{\rm n} &=\begin{cases}
       1,  & \phi^{\rm n}(\x) < - \nsmear \; h,\\
       1 - \frac{1}{2}\left(1 + \frac{1}{\nsmear h} \phi^{\rm n}(\x)+ \frac{1}{\pi} \sin\left(\frac{\pi}{ \nsmear h} \phi^{\rm n}(\x)\right)\right) ,  & |\phi^{\rm n}(\x)| \le \nsmear\; h,\\
        0,  & \textrm{otherwise}.
\end{cases}
\end{align}
In the above, $\phi^{\rm d}(\x)$ and $\phi^{\rm n}(\x)$ are the signed distance functions for the Dirichlet (outer cylinder) and Neumann (inner cylinder) boundary, respectively.  
The Rayleigh number Ra = $G \gamma Q L^4/(k \alpha \nu)$ and the Prandtl number Pr = $\nu/\alpha$ are the two main non-dimensional parameters that characterize buoyancy-driven flows; these parameters are seen in the right-hand side of the non-dimensional momentum  Eq.~\eqref{eqn_ndim_mom}. Here, $\gamma$ is the coefficient of thermal expansion, $G$ is the gravitational constant, and $\nu$ is the kinematic viscosity.  The flux-forcing function $\vbeta$ imposing the constant flux boundary condition on the surface of the inner cylinder, $\n_i \cdot \grad \Phi = 1$, is constructed using Approach B and C for this problem. Here, $\n_i$ is the unit outward normal vector of the inner cylinder. On the outer cylinder homogeneous Dirichlet boundary condition, $\Phi = 0$, is imposed through the last term of Eq.~\eqref{eqn_ndim_adv_diff}. Periodic boundary conditions are used on the external domain boundaries. 

Two Rayleigh numbers Ra = $5700$ and Ra = $5 \times 10^4$ are considered for this problem. The same Prandtl number Pr = 0.71 is used for the two cases. The computational domain $\Omega$ is discretized by a uniform Cartesian grid of size $256 \times 256$. The penalization parameters $\eta_{\rm d}$ and $\eta_{\rm n}$ are taken to be $10^{-8}$. We treat the convective and the advective terms of Eqs.~\eqref{eqn_ndim_mom} and~\eqref{eqn_ndim_adv_diff} explicitly, whereas the rest of the terms are treated implicitly. The implicit treatment of volume penalization terms in Eqs.~\eqref{eqn_ndim_mom} and~\eqref{eqn_ndim_adv_diff} allows us to use a relatively large time step sizes of  $\Delta t = 10^{-4}$ and $\Delta t = 5\times10^{-5}$  for Ra $= 5700$ and Ra $= 5 \times 10^4$ cases, respectively. In contrast, Sakurai et al.~\cite{Sakurai2019} used a time step size of $\Delta t = 10^{-6}$ for these two cases as they employed an explicit Euler time marching scheme.   More details on the second-order accurate spatial discretization and time-stepping scheme employed in the fluid solver can be found in our prior works~\cite{Nangia2019MF,Nangia2019WSI}.  

To compare our results with those reported in~\cite{Sakurai2019} we plot the steady-state temperature distribution on the left half of the inner cylinder for both Ra cases in Fig.~\ref{fig_free_conv}. In the figure, the polar angle $\Theta = 0^{\circ}$ starts from the top position $(x,y) = (0,1)$ of the inner cylinder and ends at its  bottom position $(x,y) = (0,-1)$, where $\Theta = 180^{\circ}$.  As observed in the figure, the numerical results obtained using both Approach B and C are in excellent agreement with those reported in~\cite{Sakurai2019} who used Approach A for constructing $\vbeta$. Sakurai et al. compared their numerical results with Yoo~\cite{Yoo2003} and Ren et al.~\cite{Ren2013}; comparison with Yoo and Ren et al. is therefore omitted in Fig.~\ref{fig_free_conv} in the interest of clarity. We also present the steady-state temperature field in the whole annular domain at the two Rayleigh numbers in Fig.~\ref{fig_free_conv}.


%% file: Conclusions.tex
In this work, we proposed a numerical technique for constructing flux-forcing functions for the flux-based VP method introduced by Sakurai et al. We also extended the flux-based VP approach to include Robin boundary conditions. Our method of flux-forcing functions is more general than the analytical approach (denoted Approach A in this work) of Sakurai et al. and requires only a signed distance function to construct the flux-forcing function. Two numerical-based approaches were presented for constructing flux-forcing functions: Approach B for imposing spatially constant and Approach C for imposing spatially varying (as well spatially constant) Neumann/Robin boundary conditions. Within Approach C we extended the (spatially varying) codimension-1 $g$ function to the neighborhood of the interface using a simple propagation strategy. We considered several two- and three-dimensional Poisson problems in complex domains to assess the accuracy of the numerical solutions. Results were presented for both continuous and discontinuous indicator functions. \REVIEW{For Approach B, largely $\mathcal{O}(h^2)$ accuracy is observed using the continuous indicator function. Between $\mathcal{O}(h^1)$ and $\mathcal{O}(h^2)$ convergence rate is observed for Approach B with the discontinuous indicator function and a similar convergence rate is observed for Approach C with the discontinuous indicator function when it is used for solving the constant Neumann/Robin boundary condition problem. For spatially varying boundary conditions, Approach C using the discontinuous indicator function exhibits close to $\mathcal{O}(h^1)$  convergence rate; the accuracy of the method is further improved by smoothing the sharp geometric features.  However, Approach C using the continuous indicator function exhibits a convergence rate between $\mathcal{O}(h^0)$ and $\mathcal{O}(h^1)$. Based on our results, we recommend using the discontinuous indicator function with Approach C to  impose spatially varying Neumann/Robin boundary conditions and using the continuous indicator function with Approach B to impose spatially constant Neumann/Robin boundary conditions.} Finally, we used Approach B and C to study the flux-driven thermal convection problem in a concentric annulus and compared our results against the literature. An excellent agreement was obtained. We also provided formal derivation of the flux-based volume penalized Poisson equations in strong form for both Neumann and Robin problems. The formulation shows that an explicit construction of the delta function is not necessary for the flux-based VP method, which makes it different from other diffuse domain equations presented in the literature.

%

%% file: Appendix.tex
\section{Derivation of the flux-based volume penalized Poisson equation: the Neumann problem} \label{sec_neumann_derivation}

In this section, we derive Eq.~\eqref{eqn_vp_poisson} by following the diffuse domain equation derivation provided in Li et al.~\cite{Li2009}. Similar derivation appeared in Rami\`ere et al.~\cite{Ramiere2007a}. To begin, multiply Eq.~\eqref{eq_non_penalized_poisson} by a test function $\psi$ and integrate it over the fluid domain $\Omegaf$ to obtain
\begin{equation}
\int_{\partial \Omegaf} \psi g \; \dS + \int_{\Omegaf} \kappa\; \grad q \cdot \grad \psi \; \dV = \int_{\Omegaf} \psi f \; \dV. \label{eqn_weak_fluid}
\end{equation}
In the above equation we used the vector identity 
\begin{equation} 
\grad \cdot (a \V b) = (\grad a) \cdot \V b + a\; (\grad \cdot \V{b}),
\label{eqn_vector_idn}
\end{equation} 
with the scalar field $a = \psi$ and vector field $\V{b} = \kappa \; \grad q$, along with the Neumann boundary condition on the fluid-solid interface $\partial \Omegaf$ as written in Eq.~\eqref{eq_neuman_bc}. Note that Eq.~\eqref{eqn_weak_fluid} is the weak form of the Poisson Eq.~\eqref{eq_non_penalized_poisson} defined in the fluid domain $\Omegaf$. Next, extend the integration domain from region $\Omegaf$ to $\Omega$ in the integral Eq.~\eqref{eqn_weak_fluid} by introducing the indicator function $\chi$ ($\chi = 0$ in $\Omegaf$ and $\chi = 1$ in $\Omegas$) and a surface delta function $\deltaf$ to obtain
\begin{equation}
\int_{\Omega} \deltaf\; \psi g \; \dV + \int_{\Omega} (1 - \chi) \kappa\; \grad q \cdot \grad \psi \; \dV = \int_{\Omega} (1- \chi) \psi f \; \dV. \label{eqn_weak_full}
\end{equation}
Again invoking the vector identity defined in Eq.~\eqref{eqn_vector_idn}, but with $\V{b} = (1 - \chi) \kappa \; \grad q$ this time, the second integrand in the left-hand side of the above equation can be written as  
\begin{align*}
(1 - \chi) \kappa\; \grad q \cdot \grad \psi &= \grad \cdot [\psi (1 - \chi) \kappa \; \grad q ]  -  \psi\; \grad \cdot [(1 - \chi) \kappa \; \grad q]. 
\end{align*}
This allows us to simplify the second integral in the left-hand side of Eq.~\eqref{eqn_weak_full} as
\begin{align*} 
 \int_{\Omega} (1 - \chi) \kappa\; \grad q \cdot \grad \psi \; \dV  &= \underbrace{ \int_{\partial \Omega} \psi (1 - \chi) \kappa \; (\grad q  \cdot \n_{\partial \Omega}) \; \dS}_{\text{ = 0 as } \chi  \text{ = 1 on } \partial \Omega }  - \int_{\Omega} \psi\; \grad \cdot [(1 - \chi) \kappa \; \grad q] \; \dV.
\end{align*}
Therefore, the weak form of the Poison equation in the extended domain can be written as 
\begin{equation}
\int_{\Omega} \psi \left( -\grad \cdot [(1 - \chi) \kappa\; \grad q ]+ \deltaf\; g - (1- \chi)  f  \right) \; \dV = 0.  \label{eqn_weak_full_simplified}
\end{equation}
Since $\psi$ is an arbitrary test function, the collective term multiplying $\psi$ in Eq.~\eqref{eqn_weak_full_simplified} should evaluate to zero at each point in the domain. This gives the strong form of the extended domain Poisson equation as
\begin{equation}
  -\grad \cdot [(1 - \chi) \kappa\; \grad q ]+ \deltaf\; g = (1- \chi)  f .  \label{eqn_strong_full}
\end{equation}

Next, we show that the flux-based volume penalized Poisson Eq.~\eqref{eqn_vp_poisson} can be obtained from Eq.~\eqref{eqn_strong_full} using a specific definition of the surface delta function $\deltaf$. First, simplify the forcing term $f_{\rm b}$ in the right hand-side of Eq.~\eqref{eqn_vp_poisson} to 
\begin{equation} 
\div \left(\chi \vbeta\right) - \chi \div \vbeta = \vbeta \cdot \grad \chi.
\end{equation}
Next, noticing that $\grad \chi = \deltaf\; \n$, the forcing term of VP poisson equation becomes 
\begin{equation} 
 \vbeta \cdot \grad \chi = \deltaf (\vbeta \cdot \n) = -\deltaf\; g.  \label{eqn_fb} 
\end{equation}
Substituting $\deltaf\; g$ term from Eq.~\eqref{eqn_fb} into the extended domain Poisson equation~\eqref{eqn_strong_full} \REVIEW{eliminates the explicit representation of the delta function} and the extended domain equation reads as 
\begin{equation}
  -\grad \cdot [(1 - \chi) \kappa\; \grad q ] = (1- \chi)  f  + \div \left(\chi \vbeta\right) - \chi \div \vbeta.  \label{eqn_strong_sakurai}
\end{equation}
The flux-based VP Poisson equation is obtained from Eq.~\eqref{eqn_strong_sakurai} by introducing a small amount of diffusion in the solid domain which is controlled by the penalization parameter $\eta$. For an easy reference, the VP equation is re-written below
\begin{equation*}
  -\grad \cdot [ \{\kappa(1 - \chi) + \eta \chi\} \grad q ] = (1- \chi)  f  + \div \left(\chi \vbeta\right) - \chi \div \vbeta.  \end{equation*}

\section{Derivation of the flux-based volume penalized Poisson equation: the Robin problem} \label{sec_robin_derivation}

The flux-based volume penalization method can be easily extended to include Robin boundary conditions of the type
\begin{equation}
\zeta\; q + \kappa\; \n \cdot \grad q = -g
\label{eqn_robin_app}
\end{equation}
on the irregular boundary $\partial \Omegaf$ (or $\partial \Omegas$). First, it can be easily verified that the weak form of the Poisson equation defined in the fluid domain and satisfying Robin boundary conditions written in Eq.~\eqref{eqn_robin_app} is  
\begin{equation}
\int_{\partial \Omegaf} \psi (\zeta q + g) \; \dS + \int_{\Omegaf} \kappa\; \grad q \cdot \grad \psi \; \dV = \int_{\Omegaf} \psi f \; \dV. \label{eqn_weak_fluid_robin}
\end{equation}
Next, following the procedure to reformulate the PDE on the entire domain as described in Appendix~\ref{sec_neumann_derivation}, the strong form of the Poisson equation reads as
\begin{equation}
  -\grad \cdot [(1 - \chi) \kappa\; \grad q ]+ \deltaf\; (\zeta q + g) = (1- \chi)  f .  
  \label{eqn_strong_full_robina}
\end{equation}
Defining a flux function $\vbeta$ that satisfies the property of $\vbeta \cdot \n = -g$ on $\partial \Omegaf$,  the above equation can be written as
\begin{align}
& -\grad \cdot [(1 - \chi) \kappa\; \grad q ]+ \grad \chi \cdot (\zeta q \; \n - \vbeta) = (1- \chi)  f.  
\label{eqn_strong_full_robinb}
\end{align}
The flux-based VP Poisson equation satisfying the Robin boundary conditions is obtained from Eq.~\eqref{eqn_strong_full_robinb} by adding a small amount of diffusion in the solid domain
\begin{equation}
(\zeta\; \grad \chi \cdot \n)  q -\grad \cdot [\left\{ \kappa \left(1 - \chi\right) + \eta \chi \right\} \; \grad q] = (1- \chi) f + \div (\chi \vbeta) - \chi \div \vbeta.
\label{eqn_strong_full_robinc}
\end{equation} 

In order to avoid computing the gradient of a possible discontinuous indicator function $\chi$, the above equation is re-written as
\begin{equation}
\zeta [ \div (\chi \n) - \chi \div \n ] q -\grad \cdot [\left\{ \kappa \left(1 - \chi\right) + \eta \chi \right\} \; \grad q] = (1- \chi) f + \div (\chi \vbeta) - \chi \div \vbeta.
\label{eqn_strong_full_robind}
\end{equation} 
The normal vector appearing in the first term of  Eq.~\eqref{eqn_strong_full_robind} can be computed numerically using the signed distance function as $\n = - \grad \phi$.
  
\section{Order of accuracy results using Approach A with spatially varying flux values}\label{sec_approach_a_results}

\begin{figure}[]
\centering
\subfigure[{\REVIEW{Convergence rate for the egg geometry}}]{
\includegraphics[scale = 0.08]{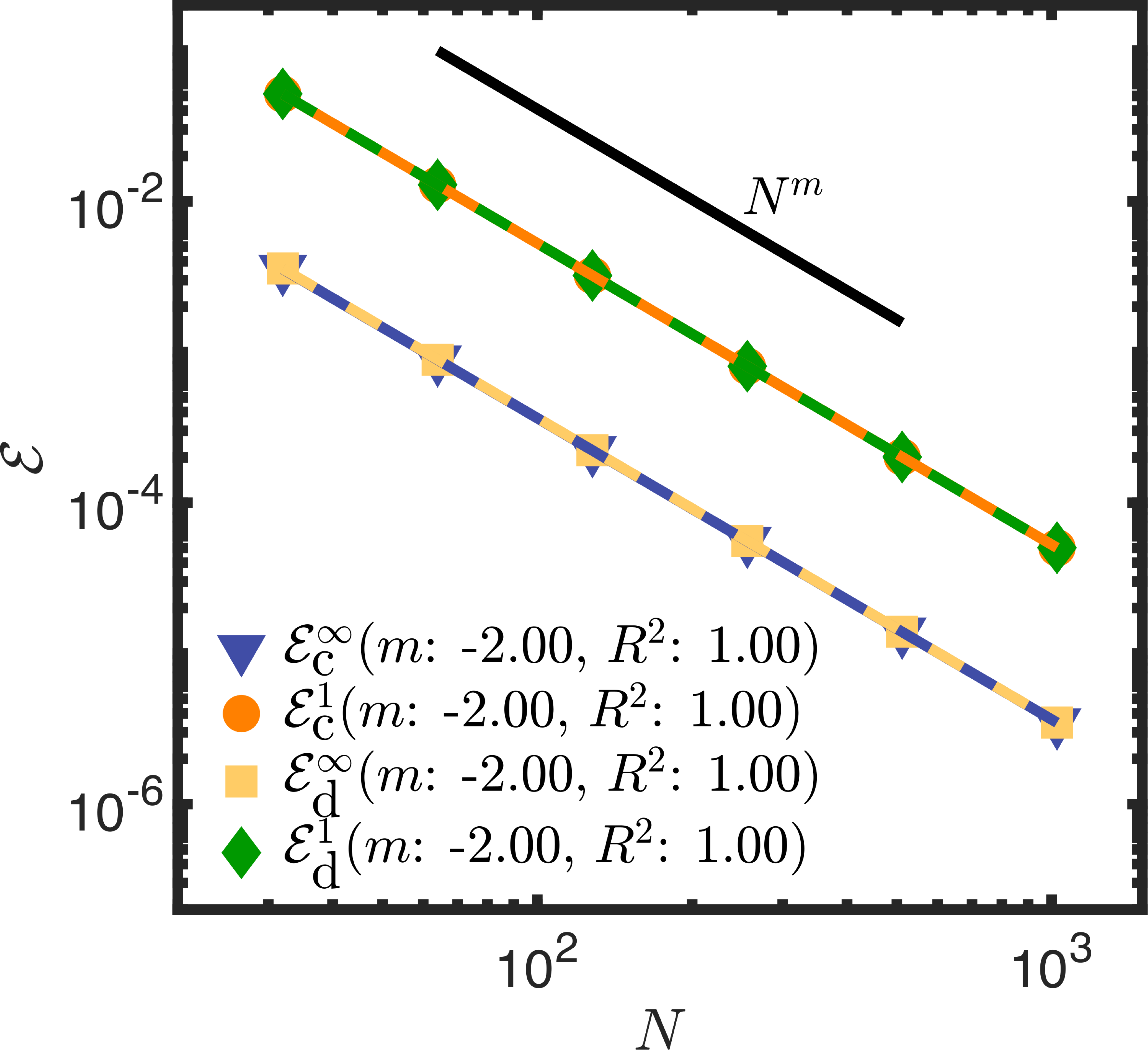}
\label{fig_egg_ooa_approach_A}
}
\subfigure[{\REVIEW{Convergence rate for the torus geometry}}]{
\includegraphics[scale = 0.08]{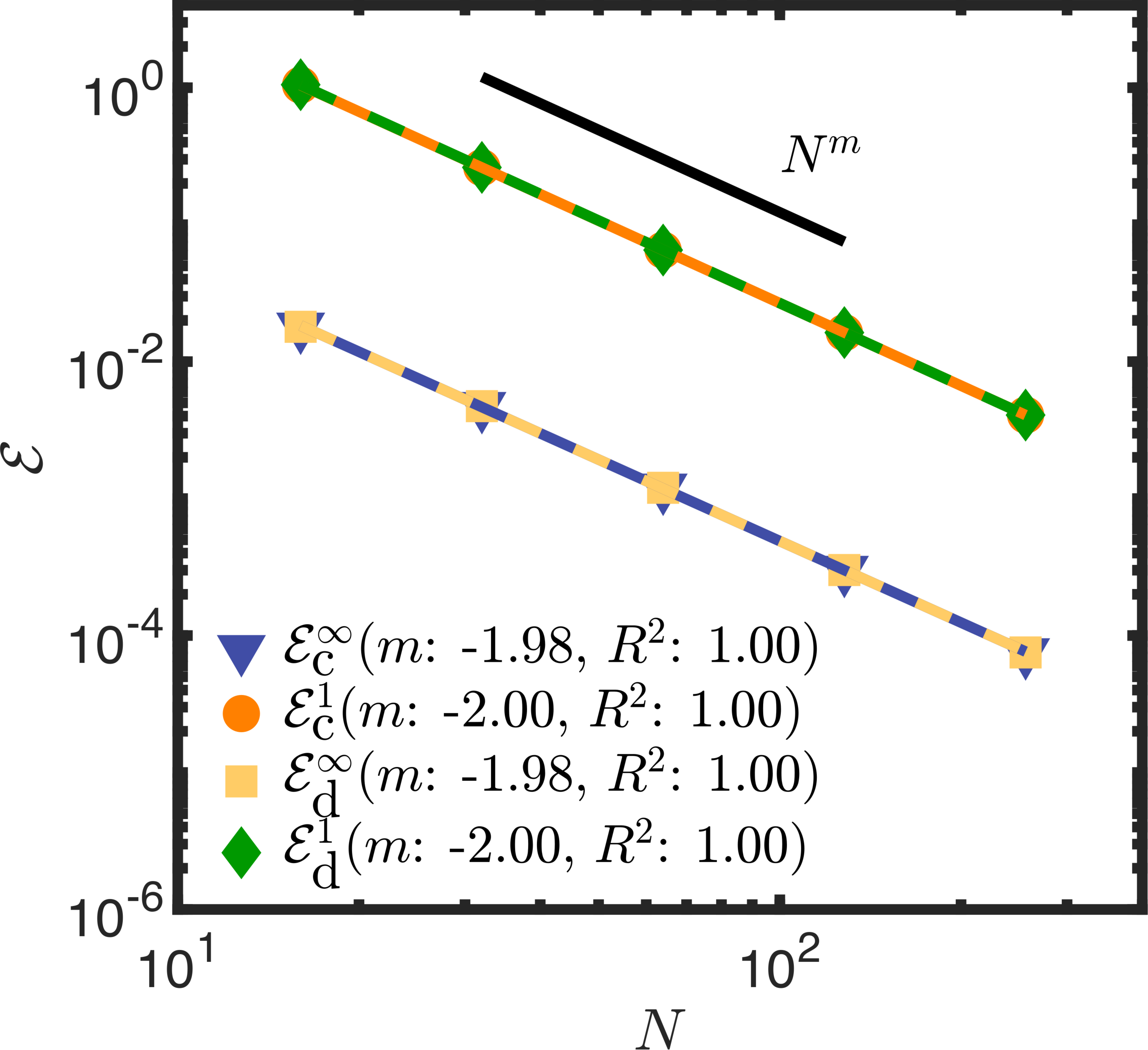}
\label{fig_torus_ooa_approach_A}
}
\caption{Convergence rate of Approach A with spatially varying flux boundary conditions. Error norms \REVIEW{$\mathcal{E}^{1}$} and \REVIEW{$\mathcal{E}^\infty$} are shown as a function of grid size $N$ using the continuous (solid lines with symbols) and discontinuous (dashed lines with symbols) indicator functions for~\subref{fig_egg_ooa_approach_A} the \REVIEW{egg}; and~\subref{fig_torus_ooa_approach_A} \REVIEW{torus domains}. }
\label{fig_approach_A}
\end{figure}

In this section, we present the spatial convergence rate of the error norms using Approach A (analytical construction of  $\vbeta$) for the egg and torus domains. The test problem remains the same as defined in Secs.~\ref{sec_spatially_varying_flux} and \ref{sec_three_dimensional} for the egg and torus domain, respectively.  The flux-forcing function using Approach A is $\vbeta = \kappa \grad \qexact$. As shown in Fig~\ref{fig_approach_A}, we observe second-order convergence rate for both problems using Approach A.

\REVIEW{
\section{Comparison of Approach C with Approach D} \label{sec_other_ext_results}

 Here, we demonstrate the efficacy of Approach C in comparison to Approach D for numerically constructing the flux-forcing functions. To implement Approach D,  we discretize the interface into a set of discrete Lagrangian/marker points with position $\X \equiv \left(X,Y\right)$ and spread the two components (in 2D) of $\V{\beta}(\X) \equiv \left( \beta_X, \beta_Y\right)$ to the nearby $x$- and $y$-faces of the Cartesian grid cells, respectively. We consider two kernel functions for the spreading operator: a one-point top hat function and a six-point Gaussian bell-like spline function. The one-dimensional (in the $x$-direction) form of the top hat function can be defined in terms of $r = (x - X)/ h$ and it reads as 
 \begin{align}
f_\textrm{top hat} &=
\begin{cases}
       1,  & |r| \le 0.5,\\
        0,  & \textrm{otherwise}.
\end{cases}       \label{eqn_top_hat}
\end{align}
Here, $x$ is the face-center location of the x-face and $h$ is the grid cell size. Similarly, the one-dimensional form of the six-point spline function reads as
 \begin{align}
f_\textrm{spline} &=
\begin{cases}
       \frac{1}{60} \left( -5\sigma^5 + 90\sigma^4 - 630\sigma^3 + 2130\sigma^2 - 3465\sigma + 2193 \right),  & 0 \le |r| < 1,\\
        \frac{1}{120} \left( 5\sigma^5 - 120\sigma^4 + 1140\sigma^3 - 5340\sigma^2 + 12270\sigma - 10974\right),  & 1 \le |r| < 2,\\
        \frac{1}{120} \left( -\sigma^5 + 30\sigma^4 - 360\sigma^3 + 2160\sigma^2 - 6480\sigma + 7776\right),  & 2 \le |r| < 3,\\
        0,  & 3 \le |r|.
\end{cases}       \label{eqn_top_hat}
\end{align}
in which $\sigma = |r| + 3$. The graphical representation of these functions is shown in Fig.~\ref{fig_extension_functional_form}.
  \begin{figure}[]
\centering
\subfigure[{\REVIEW{One-point top hat function}}]{
\includegraphics[scale = 0.08]{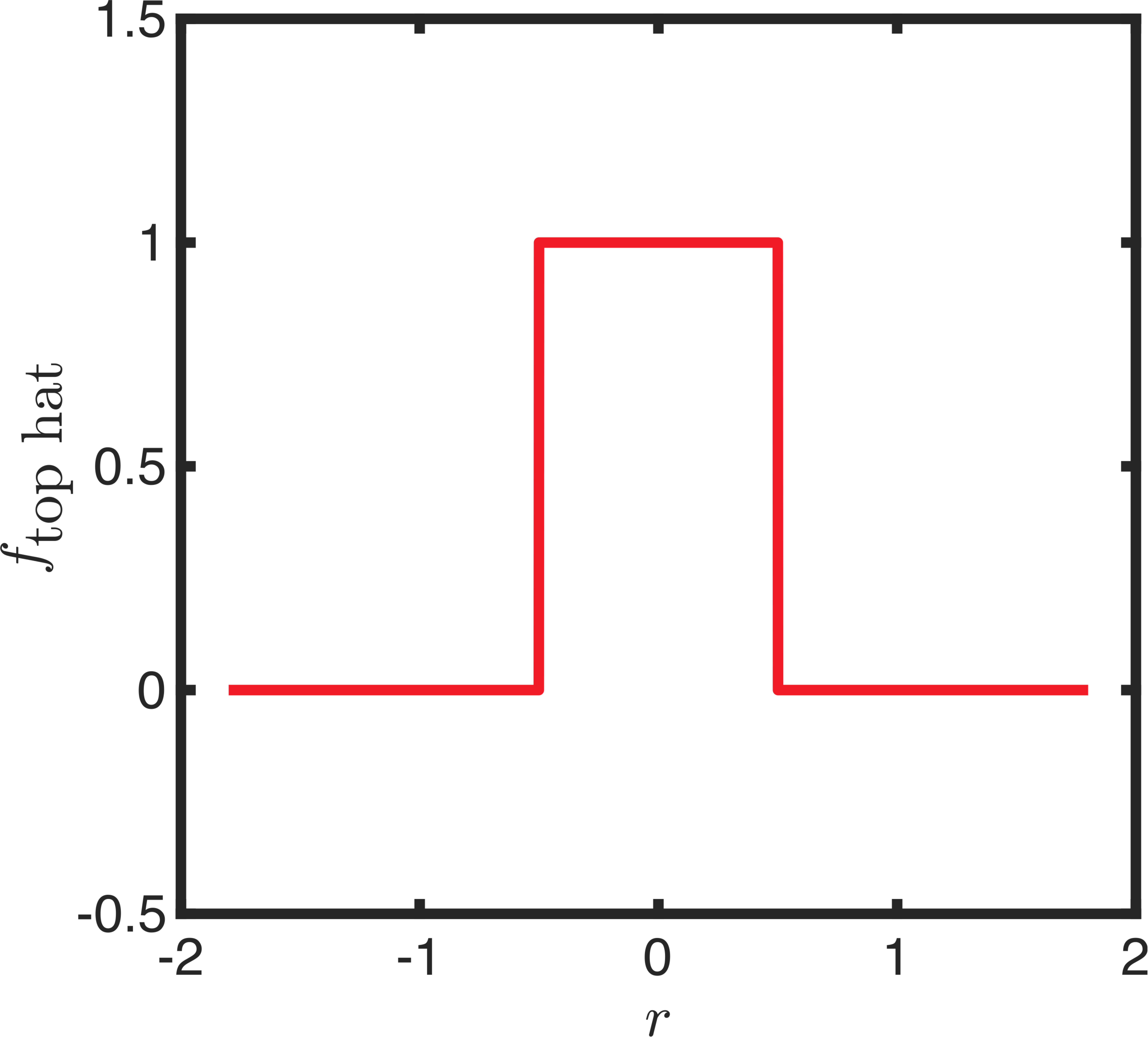}
\label{fig_top_hat_function}
}
\subfigure[{\REVIEW{Six-point spline function}}]{
\includegraphics[scale = 0.08]{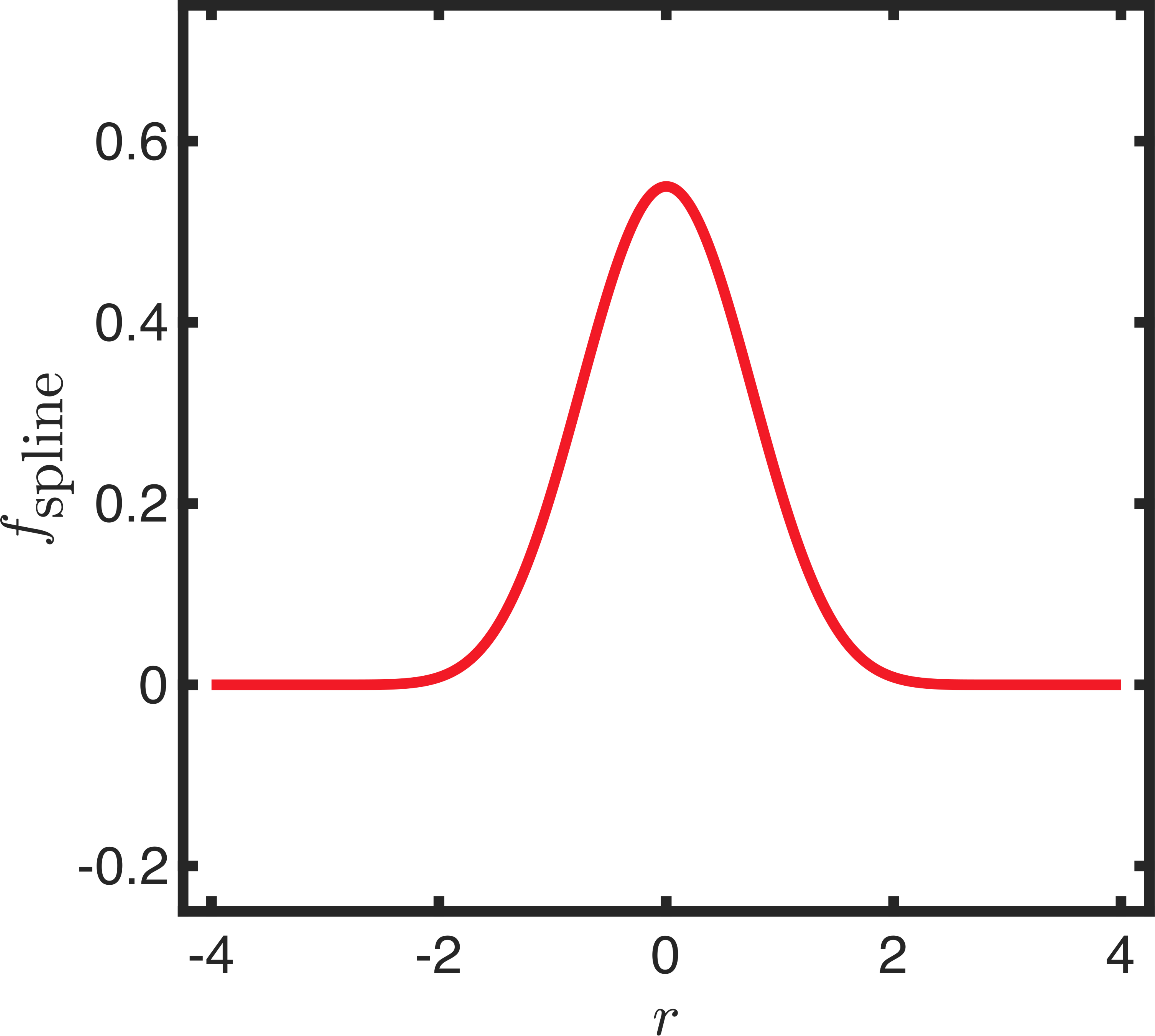}
\label{fig_spline_solution}
}
 \caption{ \REVIEW{One-dimensional representation of the~\subref{fig_top_hat_function} one-point top hat function; and~\subref{fig_circle_top_hat_solution} six-point spline function.}}
\label{fig_extension_functional_form}
\end{figure}
In dimensions higher than one, a tensor-product form of the one-dimensional functions is used. We refer the readers to Peskin~\cite{Peskin02} for more details on the spreading operator. The distance between the maker points is kept approximately equal to $h$, although increasing or decreasing the distance up to a factor of two did not affect the overall accuracy of the scheme (data not shown).  

We consider a test problem similar to the one defined in Sec.~\ref{sec_spatially_varying_flux},  in which a circle of radius 3/2 is embedded into a larger computational domain of extents $\Omega \in [0, 2\pi]^2$.  Inhomogeneous Neumann boundary conditions with $g(\x) = -\kappa\; \n \cdot \grad \qexact$ are imposed on the fluid-solid interface $\partial\Omegas$, whereas Dirichlet boundary conditions are imposed on the external boundaries of the computational domain, i.e., $\left.q\right|_{\partial \Omega(\x)} = \qexact(\x)$.  We present the contours of the numerical (red) and exact (blue) solutions and the convergence rate of the numerical schemes based on Approach C and D in Fig.~\ref{fig_Approach_C_comparison_with_other_extension_functions}. As observed in the figure, excellent agreement is obtained between the numerical and exact solutions in the fluid domain (considered to be outside the cylinder) with Approach C. However,  when Approach D is used considering either the top hat or the spline function, there is a large disagreement between the numerical and exact solutions. Furthermore, Approach C exhibits approximately first-order accuracy using the continuous indicator function, whereas Approach D exhibits zeroth-order accuracy. A similar discrepancy is observed using the discontinuous indicator function with Approach D, although it is slightly less severe than the continuous case. In contrast,  the order of accuracy of Approach C improves further when the discontinuous $\chi$ is used.  Data for the discontinuous function is not shown in the interest of brevity. Based on the results of this section we do not recommend Approach D to impose the spatially varying Neumann/Robin boundary conditions. 

 
\begin{figure}[]
\centering
\subfigure[{\REVIEW{Numerical (red) vs. analytical (blue) solution using Approach C}}]{
\includegraphics[scale = 0.137]{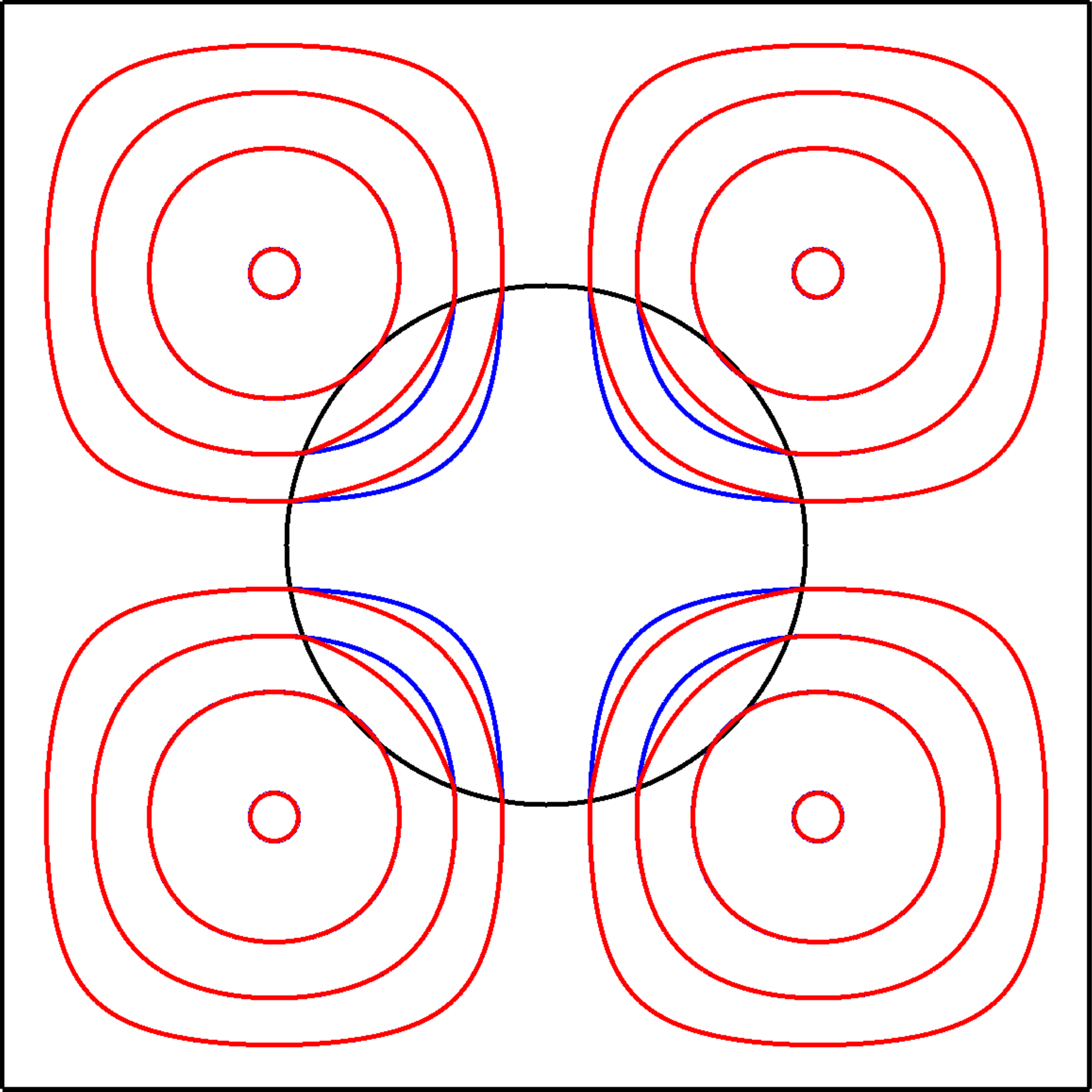}
\label{fig_circle_Approach_C_solution}
}\hspace{1 pc}
\subfigure[{\REVIEW{Numerical (red) vs. analytical (blue) solution using Approach D (top hat)} }]{
\includegraphics[scale = 0.137]{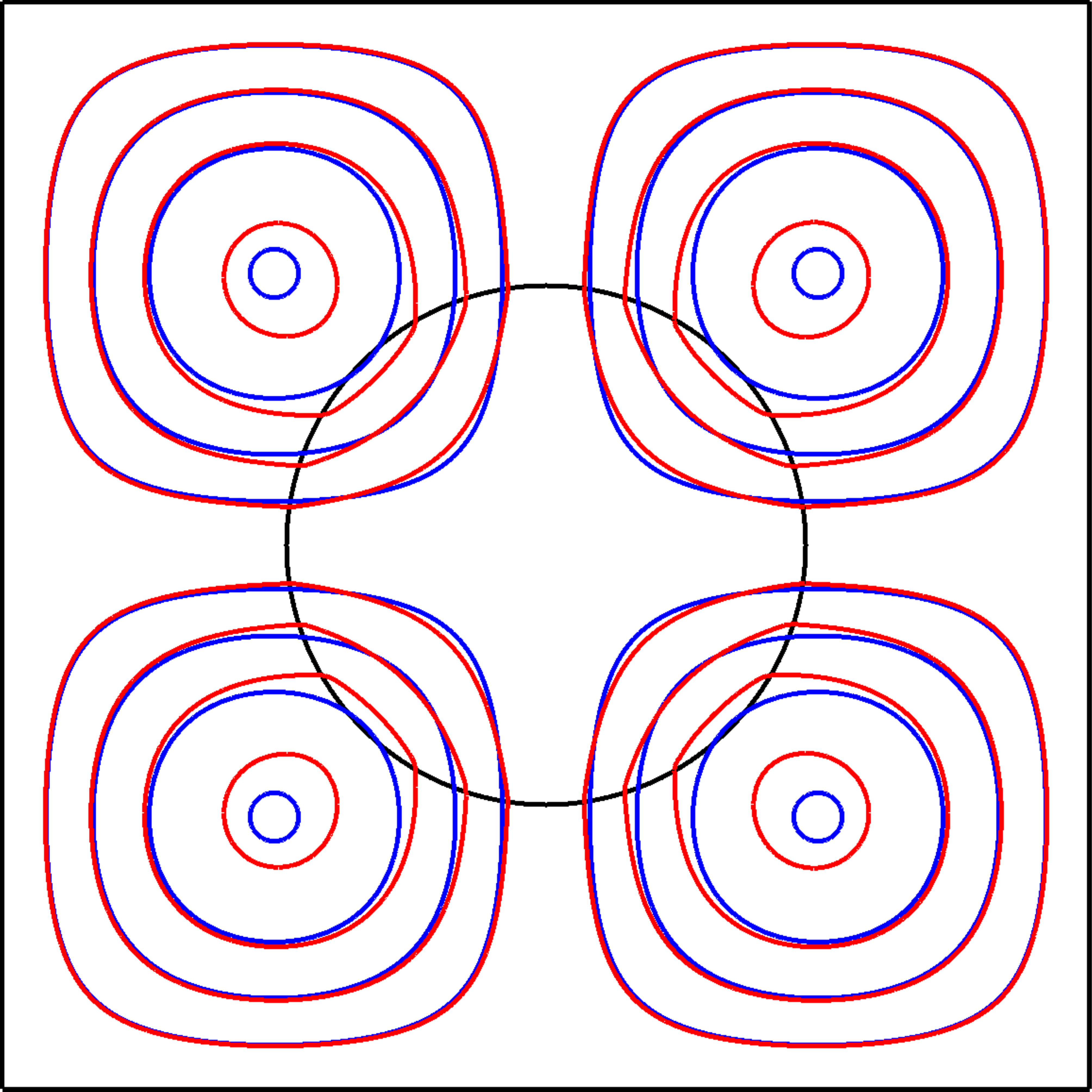}
\label{fig_circle_top_hat_solution}
}

\subfigure[{\REVIEW{Numerical (red) vs. analytical (blue) solution using Approach D (spline)}}]{
\includegraphics[scale = 0.137]{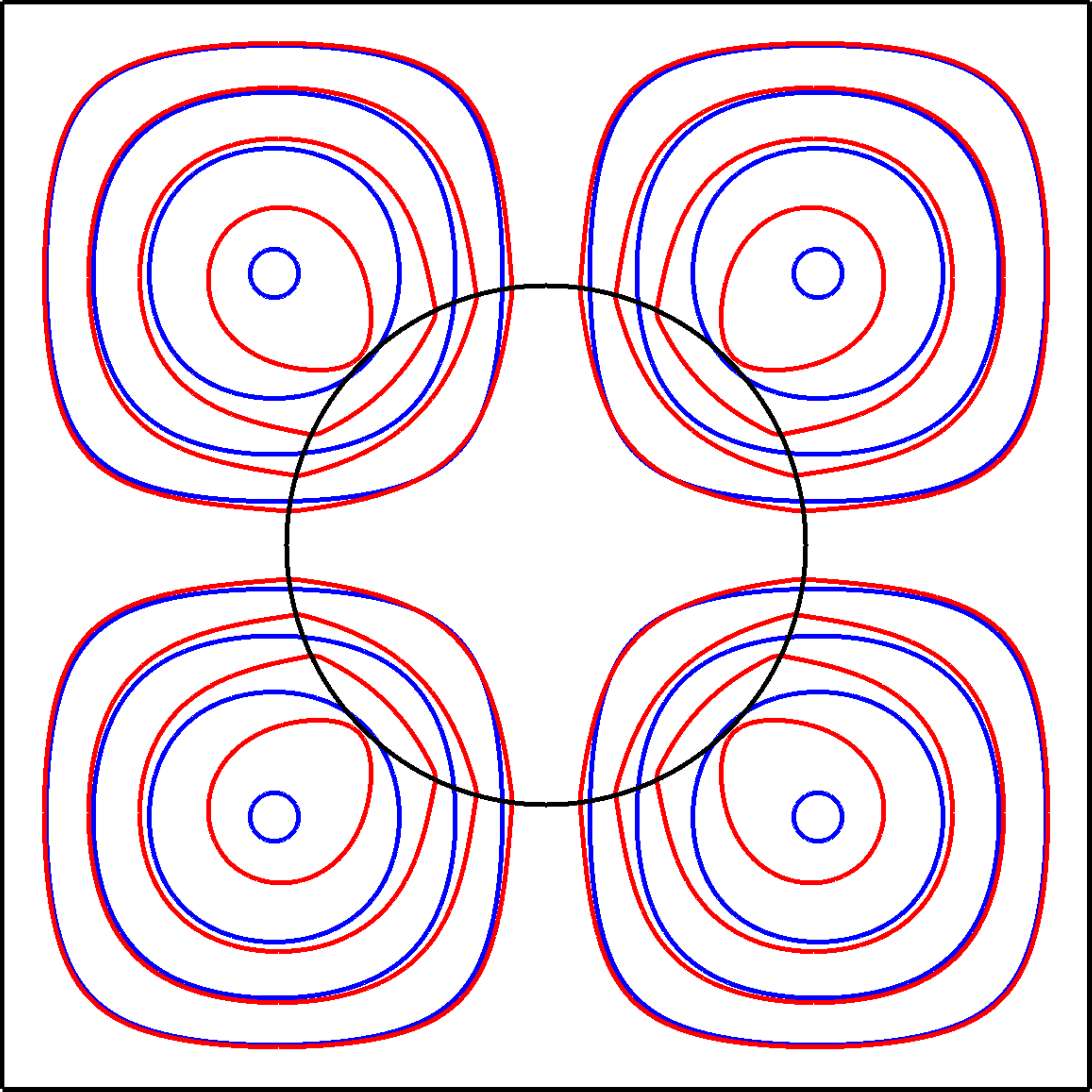}
\label{fig_circle_spline_solution}
}
\subfigure[{\REVIEW{Order of convergence using the continuous indicator function}}]{
\includegraphics[scale = 0.08]{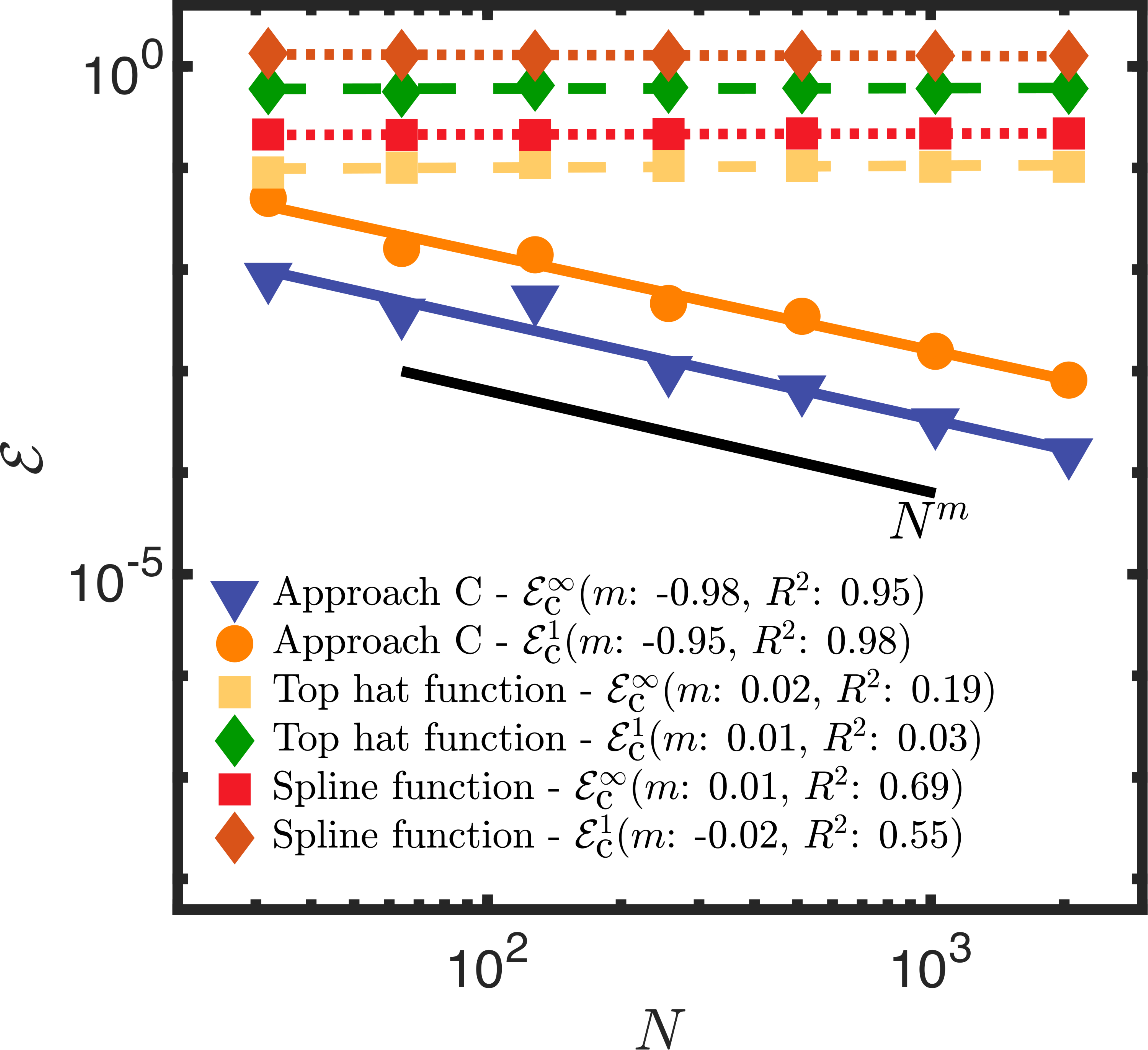}
\label{fig_circle_OOC_cont}
}
 \caption{ \REVIEW{Circular domain with spatially varying Neumann boundary conditions. Contours of the numerical (red) and analytical (blue) solutions at $N = 256$ grid using~\subref{fig_circle_Approach_C_solution}  Approach C and using Approach D with~\subref{fig_circle_top_hat_solution} one-point top hat function and~\subref{fig_circle_spline_solution} six-point spline function.~\subref{fig_circle_OOC_cont} Error norms $\mathcal{E}^{1} $ and $\mathcal{E}^\infty$ as a function of grid size $N$ using the continuous indicator function with Approach C (solid line with symbols), and with Approach D using the top hat (dashed line with symbols) and spline  (dotted line with symbols) kernel functions. The penalization parameter $\eta$ is taken as $10^{-8}$. The values of $\kappa$ is taken to be 1.}}
\label{fig_Approach_C_comparison_with_other_extension_functions}
\end{figure}
}
\REVIEW{
\section{Effect of the penalization parameter} \label{sec_eta_effect}

In this section, we study the effect of the penalization parameter $\eta$ on the order of accuracy of the flux-based VP method. The test problem described in Sec.~\ref{sec_spatially_varying_flux} for the hexagram interface is considered here. We solve this problem using four different $\eta$ values: $\eta = \{10^{-2}, 10^{-4}, 10^{-8}, 10^{-12}\}$. As noted in Fig.~\ref{fig_effect_of_eta}, except for the largest value of $\eta = 10^{-2}$, the convergence rate remains the same for the rest of the $\eta$ values. Based on the results of this section we chose $\eta = 10^{-8}$ for all our test cases.
\begin{figure}[]
\centering
\subfigure[{\REVIEW{Order of convergence using the continuous indicator function}}]{
\includegraphics[scale = 0.08]{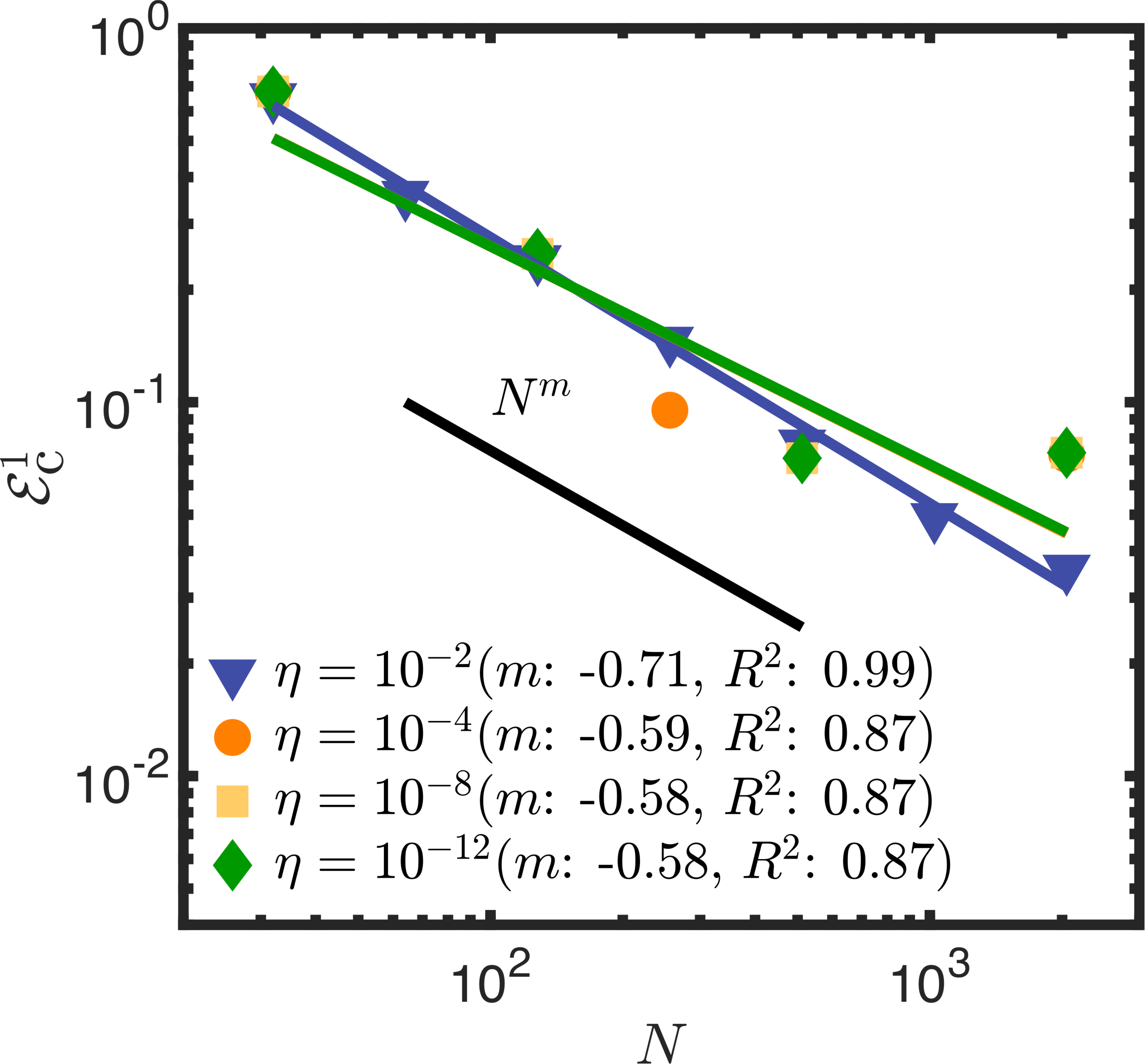}
\label{fig_effect_of_eta_L1_cont}
}
\subfigure[{\REVIEW{Order of convergence using the discontinuous indicator function}}]{
\includegraphics[scale = 0.08]{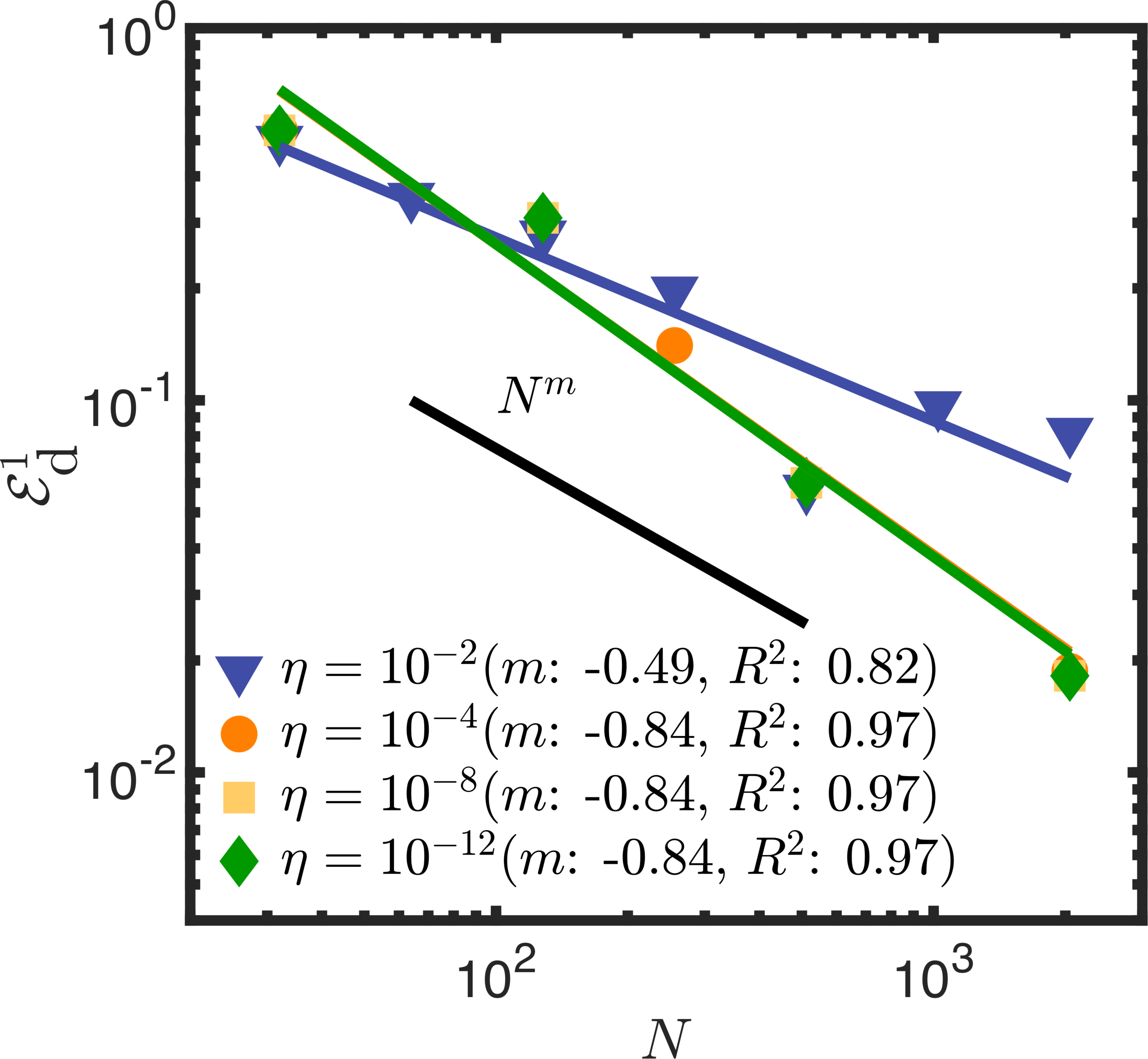}
\label{fig_effect_of_eta_L1_discont}
}
 \caption{ \REVIEW{Effect of the penalization parameter $\eta$ on the order of accuracy of the flux-based VP mehod. Spatially varying Neumann boundary conditions are imposed on the hexagram interface using Approach C.~\subref{fig_effect_of_eta_L1_cont} Error norms $\mathcal{E}^{1} $  as a function of grid size $N$ using the continuous indicator function~\subref{fig_effect_of_eta_L1_discont} Error norms $\mathcal{E}^1$  as a function of grid size $N$ using the discontinuous indicator function.} }
\label{fig_effect_of_eta}
\end{figure}
}
\REVIEW{
\section{Error norm and curve-fitting of the data}\label{sec_convergence_data}
}

\begin{table}[h!]
\centering
\caption{\REVIEW{Error norm data for the concentric circular annulus case using Approach A considered in Sec.~\ref{sec_circular_annulus_constant_flux}}}
\label{tab_xcross}
\begin{tabular}{cllllllll}
\hline
\multicolumn{1}{l}{}  & \multirow{2}{*}{N} & \multirow{2}{*}{h} & \multicolumn{3}{l}{Continuous indicator function}       & \multicolumn{3}{l}{Discontinuous indicator function}    \\ \cline{4-9} 
\multicolumn{1}{l}{}  &                    &                    & Error    & Order & Fit ($m,R^2$)        & Error    & Order & Fit ($m,R^2$)        \\ \hline
\multirow{7}{*}{$\mathcal{E}^1$} &32 & 1.96$\times10^{-1}$ & 4.6978$\times10^{-1}$        & --      & \multirow{7}{*}{1.94,0.96}     & 1.1552          & --         & \multirow{7}{*}{1.78, 0.99}    \\
                 &64     & 9.82$\times10^{-2}$  & 1.9149$\times10^{-1}$        & 1.29     &  & 7.3710$\times10^{-1}$             &  0.65     \\
               &128       & 4.91$\times10^{-2}$  & 5.3379$\times10^{-2}$        & 1.84    &   & 1.9393$\times10^{-1}$          & 1.93       \\
               &256       & 2.45$\times10^{-2}$  & 1.6787$\times10^{-2}$         & 1.67    &    & 4.8113$\times10^{-2}$          & 2.01         \\
                &512      & 1.22$\times10^{-2}$  & 9.8071$\times10^{-3}$       & 0.78     &   & 2.0514$\times10^{-2}$          & 1.23         \\
               &1024       & 6.13$\times10^{-3}$ & 7.8649$\times10^{-4}$        & 3.64   &    & 3.4221$\times10^{-3}$          & 2.58         \\
               &2048  & 3.06$\times10^{-3}$ & 1.1559$\times10^{-4}$        & 2.77   &     & 8.6893$\times10^{-4}$          & 1.98 \\  \hline
\multirow{7}{*}{$\mathcal{E}^\infty$}       &32        & 1.96$\times10^{-1}$ & 1.0015$\times10^{-2}$        & -- & \multirow{7}{*}{1.98, 0.99}          & 2.7865$\times10^{-1}$          & --         & \multirow{7}{*}{1.50, 0.98}    \\
                  &64     & 9.82$\times10^{-2}$  & 3.7365$\times10^{-2}$        & 1.42       &   & 2.3764$\times10^{-2}$          & 0.23           \\
                  &128    & 4.91$\times10^{-2}$  & 9.7656$\times10^{-3}$        & 1.94       &   & 8.0442$\times10^{-2}$          & 1.56           \\
                  &256     & 2.45$\times10^{-2}$  & 2.7216$\times10^{-3}$        & 1.84      &   & 2.1673$\times10^{-2}$          & 1.89          \\
                  &512    & 1.22$\times10^{-2}$  & 1.0573$\times10^{-3}$       & 1.36     &  & 9.4256$\times10^{-3}$          & 1.20           \\
                 &1024     & 6.13$\times10^{-3}$ & 9.2188$\times10^{-5}$        & 3.52    &      & 2.1399$\times10^{-3}$          & 2.14          \\ 
                 &2048  & 3.06$\times10^{-3}$ & 3.0885$\times10^{-5}$        & 1.58   &     & 7.9984$\times10^{-4}$          & 1.42  \\ \hline
\end{tabular}
\end{table}

\begin{table}[]
\centering
\caption{\REVIEW{Error norm data for the concentric circular annulus case using Approach B considered in Sec.~\ref{sec_circular_annulus_constant_flux}}}
\label{tab_xcross}
\begin{tabular}{cllllllll}
\hline
\multicolumn{1}{l}{}  & \multirow{2}{*}{N} & \multirow{2}{*}{h} & \multicolumn{3}{l}{Continuous indicator function}       & \multicolumn{3}{l}{Discontinuous indicator function}    \\ \cline{4-9} 
\multicolumn{1}{l}{}  &                    &                    & Error    & Order & Fit ($m,R^2$)       & Error    & Order & Fit ($m,R^2$)        \\ \hline
\multirow{7}{*}{$\mathcal{E}^1$} &32 & 1.96$\times10^{-1}$ & 4.8654$\times10^{-1}$        & --     & \multirow{7}{*}{1.93, 0.97}      & 1.1256          & --     & \multirow{7}{*}{1.78, 0.97}       \\
                      &64	& 9.82$\times10^{-2}$  & 1.9149$\times10^{-1}$        & 1.35      &  & 7.4819$\times10^{-1}$             &  0.59       \\
                     &128	& 4.91$\times10^{-2}$  & 5.3379$\times10^{-2}$        & 1.84     &  & 1.9643$\times10^{-1}$          & 1.93         \\
                      &256	& 2.45$\times10^{-2}$  & 1.6787$\times10^{-2}$         & 1.67    &    & 5.6235$\times10^{-2}$          & 1.80         \\
                      &512	& 1.22$\times10^{-2}$  & 9.8071$\times10^{-3}$       & 0.78      &  & 2.8225$\times10^{-2}$          & 0.99         \\
                      &1024	& 6.13$\times10^{-3}$ & 8.0832$\times10^{-4}$        & 3.60    &   & 3.0398$\times10^{-3}$          & 3.21        \\ 
                      &2048  & 3.06$\times10^{-3}$ & 1.2137$\times10^{-4}$        & 2.74   &     & 8.2524$\times10^{-4}$          & 1.88 \\ \hline
\multirow{7}{*}{$\mathcal{E}^\infty$}           &32	& 1.96$\times10^{-1}$ & 1.0388$\times10^{-1}$	        & --     & \multirow{7}{*}{1.99, 0.99}     & 2.6491$\times10^{-1}$          & --        & \multirow{7}{*}{1.52, 0.98}    \\
                      &64	& 9.82$\times10^{-2}$  & 3.7747$\times10^{-2}$        & 1.46        &   & 2.2513$\times10^{-1}$          & 0.23           \\
                      &128	& 4.91$\times10^{-2}$  & 9.8976$\times10^{-3}$        & 1.93        &   & 7.4914$\times10^{-2}$          & 1.59           \\
                      &256	& 2.45$\times10^{-2}$  & 2.7456$\times10^{-3}$        & 1.85        &  & 2.0745$\times10^{-2}$          & 1.85          \\
                     &512	& 1.22$\times10^{-2}$  & 1.0606$\times10^{-3}$       & 1.37      &  & 9.2775$\times10^{-3}$          & 1.16           \\
                     &1024	& 6.13$\times10^{-3}$ & 9.2313$\times10^{-5}$        & 3.52     &    & 1.8323$\times10^{-3}$          & 2.34          \\ 
                     &2048  & 3.06$\times10^{-3}$ & 3.1502$\times10^{-5}$        & 1.55   &     & 6.8778$\times10^{-4}$          & 1.41 \\ \hline
\end{tabular}
\end{table}

\begin{table}[]
\centering
\caption{\REVIEW{Error norm data for the concentric circular annulus case using Approach C considered in Sec.~\ref{sec_circular_annulus_constant_flux}}}
\label{tab_xcross}
\begin{tabular}{cllllllll}
\hline
\multicolumn{1}{l}{}  & \multirow{2}{*}{N} & \multirow{2}{*}{h} & \multicolumn{3}{l}{Continuous indicator function}       & \multicolumn{3}{l}{Discontinuous indicator function}    \\ \cline{4-9} 
\multicolumn{1}{l}{}  &                    &                    & Error    & Order & Fit ($m,R^2$)        & Error    & Order & Fit ($m,R^2$)        \\ \hline
\multirow{7}{*}{$\mathcal{E}^1$} &32 & 1.96$\times10^{-1}$ & 7.5945$\times10^{-1}$        & --     & \multirow{7}{*}{1.22, 0.82}       & 1.1582          & --       & \multirow{7}{*}{1.79, 0.97}      \\
                     &64  & 9.82$\times10^{-2}$  & 1.9216$\times10^{-1}$        & 1.98      &  & 7.5692$\times10^{-1}$             &  0.61       \\
                     &128  & 4.91$\times10^{-2}$  & 5.3385$\times10^{-2}$        & 1.85    &   & 1.9907$\times10^{-1}$          & 1.93         \\
                     &256 & 2.45$\times10^{-2}$  & 8.3493$\times10^{-2}$         & -0.64   &     & 5.7075$\times10^{-2}$          & 1.80         \\
                    &512  & 1.22$\times10^{-2}$  & 2.0512$\times10^{-2}$       & 2.03    &    & 2.8476$\times10^{-2}$          & 1.00         \\
                    &1024  & 6.13$\times10^{-3}$ & 1.7283$\times10^{-3}$        & 3.57    &   & 3.0329$\times10^{-3}$          & 3.23         \\ 
                    &2048  & 3.06$\times10^{-3}$ & 9.3256$\times10^{-3}$        & -2.43   &     & 8.2703$\times10^{-4}$          & 1.87 \\ \hline
\multirow{7}{*}{$\mathcal{E}^\infty$}    &32           & 1.96$\times10^{-1}$ & 1.3074$\times10^{-1}$	        & --       & \multirow{7}{*}{0.91, 0.79}    & 2.6907$\times10^{-1}$          & --          & \multirow{7}{*}{1.53, 0.98}  \\
                     &64  & 9.82$\times10^{-2}$  & 4.1687$\times10^{-2}$        & 1.65       &   & 2.2643$\times10^{-1}$          & 0.25           \\
                      &128 & 4.91$\times10^{-2}$  & 1.3548$\times10^{-2}$        & 1.62      &    & 7.5307$\times10^{-2}$          & 1.59           \\
                     &256 & 2.45$\times10^{-2}$  & 2.7071$\times10^{-2}$        & -1.00     &    & 2.0848$\times10^{-2}$          & 1.85          \\
                    &512  & 1.22$\times10^{-2}$  & 6.3965$\times10^{-3}$       & 2.08     &  & 9.3043$\times10^{-3}$          & 1.16          \\
                    &1024  & 6.13$\times10^{-3}$ & 1.2647$\times10^{-3}$        &2.34     &    & 1.8262$\times10^{-3}$          & 2.35          \\ 
                    &2048  & 3.06$\times10^{-3}$ & 4.7482$\times10^{-3}$        & -1.91   &     & 6.8950$\times10^{-4}$          & 1.41 \\ \hline
\end{tabular}
\end{table}

\begin{table}[]
\centering
\caption{\REVIEW{Error norm data for the hexagram case considered in Sec.~\ref{sec_spatially_varying_flux}}}
\label{tab_hexagram}
\begin{tabular}{cllllllll}
\hline
\multicolumn{1}{l}{}  & \multirow{2}{*}{N} & \multirow{2}{*}{h} & \multicolumn{3}{l}{Continuous indicator function}       & \multicolumn{3}{l}{Discontinuous indicator function}    \\ \cline{4-9} 
\multicolumn{1}{l}{}  &                    &                    & Error    & Order & Fit ($m,R^2$)        & Error    & Order & Fit ($m,R^2$)       \\ \hline
\multirow{7}{*}{$\mathcal{E}^1$} &32 & 1.96$\times10^{-1}$ & 6.7780$\times10^{-1}$        & --        & \multirow{7}{*}{0.58, 0.87}   & 5.3142$\times10^{-1}$          & --       & \multirow{7}{*}{0.84, 0.97}     \\
                   &64   & 9.82$\times10^{-2}$  & 3.4265$\times10^{-1}$        & 0.98     &   & 3.6469$\times10^{-1}$          & 0.54         \\
                   &128   & 4.91$\times10^{-2}$  & 2.4912$\times10^{-1}$        & 0.46    &    & 3.0923$\times10^{-1}$          & 0.24         \\
                   &256   & 2.45$\times10^{-2}$  & 9.4835$\times10^{-2}$         & 1.39   &     & 1.3984$\times10^{-1}$          & 1.14         \\
                   &512   & 1.22$\times10^{-2}$  & 7.0835$\times10^{-2}$        & 0.42    &    & 5.9982$\times10^{-2}$          & 1.22         \\
                   &1024   & 6.13$\times10^{-3}$ & 7.6675$\times10^{-2}$        & -0.11   &     & 3.5921$\times10^{-2}$          & 0.74        \\
                    &2048  & 3.06$\times10^{-3}$ & 7.3231$\times10^{-2}$        & 0.07   &     & 1.8182$\times10^{-2}$          & 0.98         \\ \hline
\multirow{7}{*}{$\mathcal{E}^\infty$}       &32        & 1.96$\times10^{-1}$ & 2.3021$\times10^{-1}$        & --  &\multirow{7}{*}{0.56, 0.89}           & 2.6085$\times10^{-1}$          & --      & \multirow{7}{*}{0.78, 1.00}      \\
                  &64    & 9.82$\times10^{-2}$  & 9.2871$\times10^{-2}$        & 1.31      &    & 1.4826$\times10^{-1}$          & 0.82           \\
                  &128    & 4.91$\times10^{-2}$  & 8.3922$\times10^{-2}$        & 0.15    &      & 8.4091$\times10^{-2}$          & 0.82           \\
                 &256     & 2.45$\times10^{-2}$  & 4.5081$\times10^{-2}$        & 0.90     &     & 5.7254$\times10^{-2}$          & 0.55         \\
                 &512     & 1.22$\times10^{-2}$  & 2.3251$\times10^{-2}$        & 0.96    &     & 2.8368$\times10^{-2}$          & 1.01           \\
                 &1024     & 6.13$\times10^{-3}$ & 2.7484$\times10^{-2}$        & -0.24    &      & 1.9064$\times10^{-2}$          & 0.57          \\
                  &2048    & 3.06$\times10^{-3}$ & 2.3980$\times10^{-2}$        & 0.20    &      & 9.4295$\times10^{-3}$          & 1.02           \\ \hline
\end{tabular}
\end{table}

\begin{table}[]
\centering
\caption{\REVIEW{Error norm data for the egg case considered in Sec.~\ref{sec_spatially_varying_flux}}}
\label{tab_egg}
\begin{tabular}{cllllllll}
\hline
\multicolumn{1}{l}{}  & \multirow{2}{*}{N} & \multirow{2}{*}{h} & \multicolumn{3}{l}{Continuous indicator function}       & \multicolumn{3}{l}{Discontinuous indicator function}    \\ \cline{4-9} 
\multicolumn{1}{l}{}  &                    &                    & Error    & Order & Fit ($m,R^2$)       & Error    & Order & Fit ($m,R^2$)       \\ \hline
\multirow{7}{*}{$\mathcal{E}^1$} &32 & 1.96$\times10^{-1}$ & 4.9929$\times10^{-2}$        & --       & \multirow{7}{*}{1.00, 0.95}    & 8.6565$\times10^{-1}$          & --          & \multirow{7}{*}{1.37, 0.96}  \\
                  &64    & 9.82$\times10^{-2}$  & 3.1338$\times10^{-2}$        & 0.67     &   & 5.6384$\times10^{-1}$            & 0.62         \\
                  &128     & 4.91$\times10^{-2}$  & 1.0582$\times10^{-2}$        & 1.57   &     & 5.3923$\times10^{-1}$          & 0.06         \\
                  &256    & 2.45$\times10^{-2}$  & 1.2636$\times10^{-2}$         & -0.26   &     & 1.0306$\times10^{-2}$          & 2.39         \\
                  &512     & 1.22$\times10^{-2}$  & 5.0855$\times10^{-3}$        & 1.31    &    & 3.2609$\times10^{-3}$          & 1.66         \\
                  &1024    & 6.13$\times10^{-3}$ & 1.8937$\times10^{-3}$        & 1.43    &   & 1.0862$\times10^{-3}$          & 1.59         \\
                  &2048    & 3.06$\times10^{-3}$ & 6.3391$\times10^{-4}$        & 1.58    &    & 4.4750$\times10^{-4}$          & 1.28         \\ \hline
\multirow{7}{*}{$\mathcal{E}^\infty$}    &32           & 1.96$\times10^{-1}$ & 7.2046$\times10^{-3}$        & --           & \multirow{7}{*}{0.67, 0.83} & 3.0454$\times10^{-2}$          & --         & \multirow{7}{*}{0.54, 0.91}   \\
                 &64     & 9.82$\times10^{-2}$  & 1.1199$\times10^{-2}$        & -0.64      &    & 2.5406$\times10^{-2}$          & 0.26           \\
                 &128      & 4.91$\times10^{-2}$  & 7.7484$\times10^{-3}$        & 0.53    &      & 1.9357$\times10^{-2}$          & 0.39           \\
                 &256     & 2.45$\times10^{-2}$  & 6.2520$\times10^{-3}$        & 0.31     &    & 1.7718$\times10^{-2}$          & 0.13          \\
                 &512      & 1.22$\times10^{-2}$  & 2.2338$\times10^{-3}$        & 1.48    &    & 1.0803$\times10^{-2}$          & 0.71           \\
                 &1024     & 6.13$\times10^{-3}$ & 1.6059$\times10^{-3}$        & 0.48    &      & 6.7722$\times10^{-3}$          & 0.67          \\
                 &2048     & 3.06$\times10^{-3}$ & 5.3624$\times10^{-4}$        & 1.58    &      & 2.7800$\times10^{-3}$          & 1.28           \\ \hline
\end{tabular}
\end{table}

\begin{table}[]
\centering
\caption{\REVIEW{Error norm data for the x-cross case considered in Sec.~\ref{sec_spatially_varying_flux}}}
\label{tab_xcross}
\begin{tabular}{cllllllll}
\hline
\multicolumn{1}{l}{}  & \multirow{2}{*}{N} & \multirow{2}{*}{h} & \multicolumn{3}{l}{Continuous indicator function}       & \multicolumn{3}{l}{Discontinuous indicator function}    \\ \cline{4-9} 
\multicolumn{1}{l}{}  &                    &                    & Error    & Order & Fit ($m,R^2$)        & Error    & Order & Fit ($m,R^2$)        \\ \hline
\multirow{7}{*}{$\mathcal{E}^1$} &32 & 1.96$\times10^{-1}$ & 1.9617$\times10^{-1}$        & --      & \multirow{7}{*}{0.94, 1.00}     & 7.9954$\times10^{-1}$          & --     & \multirow{7}{*}{1.08, 0.98}       \\
                    &64  & 9.82$\times10^{-2}$  & 1.0360$\times10^{-1}$        & 0.92    &    & 2.9599$\times10^{-1}$            & 1.43         \\
                   &128    & 4.91$\times10^{-2}$  & 6.3540$\times10^{-2}$        & 0.71  &     & 1.0168$\times10^{-1}$          & 1.54         \\
                    &256   & 2.45$\times10^{-2}$  & 3.3087$\times10^{-2}$         & 0.94  &      & 8.4574$\times10^{-2}$          & 0.27         \\
                    &512  & 1.22$\times10^{-2}$  & 1.6187$\times10^{-2}$       & 1.03   &     & 4.1020$\times10^{-2}$          & 1.04         \\
                   &1024   & 6.13$\times10^{-3}$ & 8.0014$\times10^{-3}$        & 1.01  &     & 1.7976$\times10^{-2}$          & 1.19         \\
                  &2048    & 3.06$\times10^{-3}$ & 3.8168$\times10^{-3}$        & 1.07   &     & 6.5705$\times10^{-3}$          & 1.45         \\ \hline
\multirow{7}{*}{$\mathcal{E}^\infty$}        &32       & 1.96$\times10^{-1}$ & 9.7360$\times10^{-2}$        & --   & \multirow{7}{*}{0.67, 0.98}        & 1.5537$\times10^{-1}$          & --        & \multirow{7}{*}{0.85, 1.00}    \\
                  &64    & 9.82$\times10^{-2}$  & 7.0146$\times10^{-2}$        & 0.42     &     & 9.9996$\times10^{-2}$          & 0.64           \\
                  &128     & 4.91$\times10^{-2}$  & 4.9175$\times10^{-3}$        & 0.51   &       & 6.0211$\times10^{-2}$          & 0.73           \\
                  &256     & 2.45$\times10^{-2}$  & 3.5404$\times10^{-3}$        & 0.47    &     & 2.6342$\times10^{-2}$          & 1.19          \\
                 &512     & 1.22$\times10^{-2}$  & 2.0327$\times10^{-2}$        & 0.80     &   & 1.5461$\times10^{-2}$          & 0.77           \\
                 &1024     & 6.13$\times10^{-3}$ & 1.1164$\times10^{-2}$        & 0.86    &      & 8.8542$\times10^{-3}$          & 0.80          \\
                  &2048    & 3.06$\times10^{-3}$ & 5.7269$\times10^{-3}$        & 0.96    &     & 4.9815$\times10^{-3}$          & 0.83           \\ \hline
\end{tabular}
\end{table}

\begin{table}[]
\centering
\caption{\REVIEW{Error norm data for fluid inside the sphere case using Approach B considered in Sec.~\ref{sec_three_dimensional}}}
\label{tab_sphere_fo_apparoach_B}
\begin{tabular}{cllllllll}
\hline
\multicolumn{1}{l}{}  & \multirow{2}{*}{N} & \multirow{2}{*}{h} & \multicolumn{3}{l}{Continuous indicator function}       & \multicolumn{3}{l}{Discontinuous indicator function}    \\ \cline{4-9} 
\multicolumn{1}{l}{}  &                    &                    & Error    & Order & Fit ($m,R^2$)        & Error    & Order & Fit ($m,R^2$)        \\ \hline
\multirow{6}{*}{$\mathcal{E}^1$} &16 & 3.92$\times10^{-1}$ &  1.2034       & --    & \multirow{6}{*}{1.96, 0.91}       &1.2034           & --   &\multirow{6}{*}{1.97, 0.93}          \\
                  &32     & 1.96$\times10^{-1}$  & 8.3704$\times10^{-1}$        & 3.85     &  & 9.8016$\times10^{-2}$            & 3.62    &     \\
                   &64    & 9.82$\times10^{-2}$  & 4.8865$\times10^{-1}$        & 0.78     &  & 4.9230$\times10^{-2}$          & 0.99         \\
                  &128     & 4.91$\times10^{-2}$  & 1.0000$\times10^{-2}$         & 2.29   &     & 1.0273$\times10^{-2}$          & 2.26         \\
                    &256   & 2.45$\times10^{-2}$  & 8.7930$\times10^{-3}$       & 0.19   &     & 8.7922$\times10^{-3}$          & 0.22        \\
                    &320      & 1.96$\times10^{-2}$ & 9.7462$\times10^{-4}$        & 9.86    &   & 1.0631$\times10^{-3}$          & 9.47  \\ \hline
\multirow{6}{*}{$\mathcal{E}^\infty$}     &16          & 3.92$\times10^{-1}$ &9.7340$\times10^{-2}$        & --   &\multirow{6}{*}{1.83, 0.93}          & 1.0115$\times10^{-1}$          & --     & \multirow{6}{*}{1.51, 0.98}       \\
                 &32      & 1.96$\times10^{-1}$   & 7.2565$\times10^{-3}$        & 3.75      &    & 2.1147$\times10^{-2}$          & 2.26           \\
                  &64     & 9.82$\times10^{-2}$  & 4.9534$\times10^{-3}$        & 0.55      &    & 1.0438$\times10^{-2}$          & 1.02           \\
                  &128     & 4.91$\times10^{-2}$   & 1.1726$\times10^{-3}$        & 2.08     &    & 2.5863$\times10^{-3}$          & 2.01          \\
                   &256    &  2.45$\times10^{-2}$  & 7.4807$\times10^{-4}$        & 0.65    &    & 1.5037$\times10^{-3}$          & 0.78          \\ 
                   &320      & 1.96$\times10^{-2}$ & 1.4668$\times10^{-4}$        & 7.30    &   & 8.9006$\times10^{-4}$          & 2.35  \\ \hline
\end{tabular}
\end{table}

\begin{table}[]
\centering
\caption{\REVIEW{Error norm data for fluid inside the sphere case using Approach C considered in Sec.~\ref{sec_three_dimensional}}}
\label{tab_sphere_fo_apparoach_B}
\begin{tabular}{cllllllll}
\hline
\multicolumn{1}{l}{}  & \multirow{2}{*}{N} & \multirow{2}{*}{h} & \multicolumn{3}{l}{Continuous indicator function}       & \multicolumn{3}{l}{Discontinuous indicator function}    \\ \cline{4-9} 
\multicolumn{1}{l}{}  &                    &                    & Error    & Order &  Fit ($m,R^2$)        & Error    & Order &  Fit ($m,R^2$)        \\ \hline
\multirow{6}{*}{$\mathcal{E}^1$} &16 & 3.92$\times10^{-1}$ & 1.2034       & --    & \multirow{6}{*}{0.94, 0.74}       & 1.2034         & --   &\multirow{6}{*}{1.97, 0.93}          \\
                  &32     & 1.96$\times10^{-1}$  & 1.4493$\times10^{-1}$        & 3.05     &  & 1.0073$\times10^{-1}$            & 3.58    &     \\
                   &64    & 9.82$\times10^{-2}$  & 1.2962$\times10^{-1}$        & 0.16     &  & 4.9176$\times10^{-2}$          & 1.03         \\
                  &128     & 4.91$\times10^{-2}$  & 3.4356$\times10^{-2}$         & 1.92   &     & 1.0362$\times10^{-2}$          & 2.25         \\
                    &256   & 2.45$\times10^{-2}$  & 4.9902$\times10^{-2}$       & -0.54   &     & 8.7919$\times10^{-3}$          & 0.24         \\ 
                    &320      & 1.96$\times10^{-2}$ & 5.6696$\times10^{-2}$        & -0.57    &   & 1.0625$\times10^{-3}$          & 9.47  \\ \hline
\multirow{6}{*}{$\mathcal{E}^\infty$}     &16          & 3.92$\times10^{-1}$ &1.5668$\times10^{-1}$        & --   &\multirow{6}{*}{0.07, 0.05}          & 1.0083$\times10^{-1}$          & --     & \multirow{6}{*}{1.52, 0.99}       \\
                 &32      & 1.96$\times10^{-1}$   & 4.9131$\times10^{-2}$        & 1.67      &    & 2.3972$\times10^{-2}$          & 2.07           \\
                  &64     & 9.82$\times10^{-2}$  & 1.0494$\times10^{-1}$        & -1.09      &    & 9.9420$\times10^{-3}$          & 1.27           \\
                  &128     & 4.91$\times10^{-2}$   & 6.0320$\times10^{-2}$        & 0.80     &    & 2.7054$\times10^{-3}$          & 1.88          \\
                   &256    &  2.45$\times10^{-2}$  & 8.8743$\times10^{-2}$        & -0.56    &    & 1.4489$\times10^{-3}$          & 0.90          \\ 
                   &320      & 1.96$\times10^{-2}$ & 9.3188$\times10^{-2}$        & -0.22    &   & 9.2746$\times10^{-4}$          & 2.00  \\ \hline
\end{tabular}
\end{table}

\begin{table}[]
\centering
\caption{\REVIEW{Error norm data for fluid outside the sphere case using Approach B considered in Sec.~\ref{sec_three_dimensional}}}
\label{tab_sphere_fo_apparoach_B}
\begin{tabular}{cllllllll}
\hline
\multicolumn{1}{l}{}  & \multirow{2}{*}{N} & \multirow{2}{*}{h} & \multicolumn{3}{l}{Continuous indicator function}       & \multicolumn{3}{l}{Discontinuous indicator function}    \\ \cline{4-9} 
\multicolumn{1}{l}{}  &                    &                    & Error    & Order & Fit ($m,R^2$)         & Error    & Order & Fit ($m,R^2$)         \\ \hline
\multirow{6}{*}{$\mathcal{E}^1$} &16 & 3.92$\times10^{-1}$ & 6.3869        & --    & \multirow{6}{*}{2.00, 1.00}       & 12.6388          & --   &\multirow{6}{*}{2.17, 1.00}          \\
                  &32     & 1.96$\times10^{-1}$  & 1.6017        & 2.00     &  & 2.3072            & 2.45    &     \\
                   &64    & 9.82$\times10^{-2}$  & 4.0250$\times10^{-1}$        & 1.99     &  & 3.9686$\times10^{-1}$          & 2.54         \\
                  &128     & 4.91$\times10^{-2}$  & 1.0047$\times10^{-1}$         & 2.00   &     & 1.0267$\times10^{-1}$          & 1.95         \\
                    &256   & 2.45$\times10^{-2}$  & 2.5157$\times10^{-2}$       & 2.00   &     & 2.8038$\times10^{-2}$          & 1.87         \\ 
                     &320      & 1.96$\times10^{-2}$ & 1.6090$\times10^{-2}$        & 2.00    &   & 1.7741$\times10^{-2}$          & 2.05  \\ \hline
\multirow{6}{*}{$\mathcal{E}^\infty$}     &16          & 3.92$\times10^{-1}$ &3.8521$\times10^{-2}$        & --   &\multirow{6}{*}{2.00, 1.00}          & 1.3960$\times10^{-1}$          & --     & \multirow{6}{*}{1.66, 0.97}       \\
                 &32      & 1.96$\times10^{-1}$   & 9.6372$\times10^{-3}$        & 2.00      &    & 3.6464$\times10^{-2}$          & 1.94           \\
                  &64     & 9.82$\times10^{-2}$  & 2.4095$\times10^{-3}$        & 2.00      &    & 5.3588$\times10^{-3}$          & 2.77           \\
                  &128     & 4.91$\times10^{-2}$   & 6.0239$\times10^{-4}$        & 2.00     &    & 3.1807$\times10^{-3}$          & 0.75          \\
                   &256    &  2.45$\times10^{-2}$  & 1.5059$\times10^{-4}$        & 2.00    &    & 1.0506$\times10^{-3}$          & 1.60          \\ 
                    &320      & 1.96$\times10^{-2}$ & 9.6382$\times10^{-5}$        & 2.00    &   & 9.9368$\times10^{-4}$          & 0.25  \\ \hline
\end{tabular}
\end{table}

\begin{table}[]
\centering
\caption{\REVIEW{Error norm data for fluid outside the sphere case using Approach C considered in Sec.~\ref{sec_three_dimensional}}}
\label{tab_sphere_fo_approach_C}
\begin{tabular}{cllllllll}
\hline
\multicolumn{1}{l}{}  & \multirow{2}{*}{N} & \multirow{2}{*}{h} & \multicolumn{3}{l}{Continuous indicator function}       & \multicolumn{3}{l}{Discontinuous indicator function}    \\ \cline{4-9} 
\multicolumn{1}{l}{}  &                    &                    & Error    & Order & Fit ($m,R^2$)        & Error    & Order & Fit ($m,R^2$)       \\ \hline
\multirow{6}{*}{$\mathcal{E}^1$} &16 & 3.92$\times10^{-1}$  & 11.4085        & --         & \multirow{6}{*}{0.87, 0.90}   & 11.3418          & --        & \multirow{6}{*}{2.13, 1.00}     \\
                &32      & 1.96$\times10^{-1}$  & 3.3852        & 1.75      &  & 1.9539            & 2.54         \\
                &64      & 9.82$\times10^{-2}$  & 1.9619        & 0.79     &  & 3.8731$\times10^{-1}$          & 2.33         \\
                &128      & 4.91$\times10^{-2}$  & 8.3413$\times10^{-1}$         & 1.23    &    & 9.6877$\times10^{-2}$          & 2.00         \\
                &256      & 2.45$\times10^{-2}$  & 6.6588$\times10^{-1}$       & 0.33      &  & 3.1445$\times10^{-2}$          & 1.62         \\
                &320      & 1.96$\times10^{-2}$ & 8.8538$\times10^{-1}$        & -1.28    &   & 1.5689$\times10^{-2}$          & 3.12         \\ \hline
\multirow{6}{*}{$\mathcal{E}^\infty$}       &16        & 3.92$\times10^{-1}$ & 3.7012$\times10^{-1}$        & --  & \multirow{6}{*} {0.41, 0.67}        & 1.1447$\times10^{-1}$          & --      & \multirow{6}{*} {1.60, 0.98}     \\
                 &32     & 1.96$\times10^{-1}$ & 1.4986$\times10^{-1}$       & 1.30      &    & 3.0673$\times10^{-2}$          & 1.90           \\
                 &64     & 9.82$\times10^{-2}$  & 1.4881$\times10^{-1}$        & 0.01    &      & 5.9817$\times10^{-3}$          & 2.36           \\
                &128      & 4.91$\times10^{-2}$  & 6.7817$\times10^{-2}$        & 1.13     &    & 2.8228$\times10^{-3}$          & 1.08         \\
                &256      & 2.45$\times10^{-2}$  & 9.2122$\times10^{-2}$        & -0.44     &  & 1.1423$\times10^{-3}$          & 1.31           \\
                 &320     & 1.96$\times10^{-2}$& 1.0649$\times10^{-1}$        & -0.65       &   & 9.3613$\times10^{-4}$          & 0.89          \\ \hline
\end{tabular}
\end{table}

\begin{table}[]
\centering
\caption{\REVIEW{Error norm data for the torus case using Approach C considered in Sec.~\ref{sec_three_dimensional}}}`
\label{tab_torus}
\begin{tabular}{cllllllll}
\hline
\multicolumn{1}{l}{}  & \multirow{2}{*}{N} & \multirow{2}{*}{h} & \multicolumn{3}{l}{Continuous indicator function}       & \multicolumn{3}{l}{Discontinuous indicator function}    \\ \cline{4-9} 
\multicolumn{1}{l}{}  &                    &                    & Error    & Order & Fit ($m,R^2$)        & Error    & Order & Fit ($m,R^2$)      \\ \hline
\multirow{7}{*}{$\mathcal{E}^1$} &16 & 3.92$\times10^{-1}$  & 2.0494       & --      & \multirow{6}{*} {1.25, 0.93}      & 2.7261         & --           & \multirow{6}{*} {1.45, 0.97}  \\
                  &32    & 1.96$\times10^{-1}$  & 5.1628$\times10^{-1}$        & 1.99       & & 1.2559            & 1.12         \\
                  &64    & 9.82$\times10^{-2}$  & 1.1226$\times10^{-1}$       & 2.20     &  & 2.1936$\times10^{-1}$          & 2.52         \\
                  &128    & 4.91$\times10^{-2}$  & 8.5599$\times10^{-2}$         & 0.39    &    & 8.7687$\times10^{-2}$          & 1.32         \\
                  &256    & 2.45$\times10^{-2}$  &4.7889$\times10^{-2}$       & 0.84      & & 5.2780$\times10^{-2}$          & 0.73         \\
                  &320    & 1.96$\times10^{-2}$ & 4.0430$\times10^{-2}$        & 0.76    &   & 3.9756$\times10^{-2}$          & 1.27         \\ \hline
\multirow{7}{*}{$\mathcal{E}^\infty$}       &16        & 3.92$\times10^{-1}$ & 7.1337$\times10^{-2}$        & --    & \multirow{6}{*} {0.63, 0.59}        & 1.8828$\times10^{-1}$          & --         & \multirow{6}{*} {0.95, 0.90}    \\
                  &32    & 1.96$\times10^{-1}$ & 1.6299$\times10^{-2}$       & 2.13      &   & 9.5677$\times10^{-2}$          & 0.98           \\
                  &64    & 9.82$\times10^{-2}$  & 4.7456$\times10^{-2}$        & 1.78      &    & 2.5818$\times10^{-2}$          & 1.89           \\
                  &128    & 4.91$\times10^{-2}$  & 8.9841$\times10^{-3}$        & -0.92    &     & 1.2753$\times10^{-2}$          & 1.01         \\
                 &256     & 2.45$\times10^{-2}$  & 7.6094$\times10^{-3}$        & 0.24    &   & 1.4638$\times10^{-2}$          & -0.20           \\
                 &320     & 1.96$\times10^{-2}$& 6.7303$\times10^{-3}$        & 0.55     &    & 1.1201$\times10^{-2}$          & 1.20          \\ \hline
\end{tabular}
\end{table}

\begin{table}[]
\centering
\caption{\REVIEW{Error norm data for the concentric annulus case with spatially constant Robin boundary conditions using Approach B considered in Sec.~\ref{sec_spatially_constant_robin}}}
\label{tab_torus}
\begin{tabular}{cllllllll}
\hline
\multicolumn{1}{l}{}  & \multirow{2}{*}{N} & \multirow{2}{*}{h} & \multicolumn{3}{l}{Continuous indicator function}       & \multicolumn{3}{l}{Discontinuous indicator function}    \\ \cline{4-9} 
\multicolumn{1}{l}{}  &                    &                    & Error    & Order & Fit ($m,R^2$)        & Error    & Order & Fit ($m,R^2$)        \\ \hline
\multirow{7}{*}{$\mathcal{E}^1$} &32 & 1.96$\times10^{-1}$ & 5.9218$\times10^{-1}$        & --      & \multirow{7}{*} {1.97, 1.00}      & 3.1462         & --        & \multirow{7}{*} {1.67, 0.95}     \\
                   &64   & 9.82$\times10^{-2}$  & 1.5389$\times10^{-1}$        & 1.94    &    & 1.5899            & 0.98         \\
                   &128   & 4.91$\times10^{-2}$  & 3.9792$\times10^{-2}$        & 1.95   &    & 4.3190$\times10^{-1}$          & 1.88         \\
                   &256   & 2.45$\times10^{-2}$  & 1.0035$\times10^{-2}$         & 1.99   &     & 4.4843$\times10^{-2}$          & 3.27         \\
                   &512   & 1.22$\times10^{-2}$  & 2.5860$\times10^{-3}$       & 1.96    &    & 4.9145$\times10^{-2}$          & -0.13         \\
                   &1024   & 6.13$\times10^{-3}$ & 6.7235$\times10^{-4}$        & 1.94   &    & 2.4537$\times10^{-2}$          & 1.00         \\ 
                   &2048    & 3.06$\times10^{-3}$ & 1.6595$\times10^{-4}$        & 2.02    &     & 2.0854$\times10^{-3}$          & 3.56           \\ \hline
\multirow{7}{*}{$\mathcal{E}^\infty$}      &32         & 1.96$\times10^{-1}$ & 7.3004$\times10^{-2}$        & --    & \multirow{7}{*} {1.97, 1.00}        & 3.3663$\times10^{-1}$          & --         & \multirow{7}{*} {1.45, 0.99}    \\
                  &64    & 9.82$\times10^{-2}$  & 1.9263$\times10^{-2}$        & 1.92     &     & 1.7209$\times10^{-2}$          & 0.97           \\
                  &128    & 4.91$\times10^{-2}$  & 4.9511$\times10^{-3}$        & 1.96    &      & 5.1747$\times10^{-2}$          & 1.73           \\
                 &256     & 2.45$\times10^{-2}$  & 1.2529$\times10^{-3}$        & 1.98    &     & 1.4644$\times10^{-2}$          & 1.82          \\
                 &512     & 1.22$\times10^{-2}$  & 3.1817$\times10^{-4}$        & 1.98    &    & 7.3577$\times10^{-3}$          & 0.99          \\
                 &1024     & 6.13$\times10^{-3}$ & 8.2175$\times10^{-5}$        & 1.95    &      & 3.2291$\times10^{-3}$          & 1.19          \\ 
                 &2048    & 3.06$\times10^{-3}$ & 2.0329$\times10^{-5}$        & 2.02    &     & 7.7661$\times10^{-4}$          & 2.06           \\ \hline
\end{tabular}
\end{table}

\begin{table}[]
\centering
\caption{\REVIEW{Error norm data for the concentric annulus case with spatially constant Robin boundary condition using Approach C considered in Sec.~\ref{sec_spatially_constant_robin}}}
\label{tab_torus}
\begin{tabular}{cllllllll}
\hline
\multicolumn{1}{l}{}  & \multirow{2}{*}{N} & \multirow{2}{*}{h} & \multicolumn{3}{l}{Continuous indicator function}       & \multicolumn{3}{l}{Discontinuous indicator function}    \\ \cline{4-9} 
\multicolumn{1}{l}{}  &                    &                    & Error    & Order & Fit ($m,R^2$)         & Error    & Order & Fit ($m,R^2$)         \\ \hline
\multirow{7}{*}{$\mathcal{E}^1$} &32 & 1.96$\times10^{-1}$ & 1.2927        & --      & \multirow{7}{*} {0.92, 0.56}      & 3.2195        & --        & \multirow{7}{*} {1.68, 0.95}     \\
                  &64    & 9.82$\times10^{-2}$  & 1.1437$\times10^{-1}$        & 3.50     &   & 1.6086            & 1.00         \\
                  &128    & 4.91$\times10^{-2}$  & 2.8402$\times10^{-2}$        & 2.01   &    & 4.3717$\times10^{-1}$          & 1.88         \\
                  &256    & 2.45$\times10^{-2}$  & 2.2619$\times10^{-1}$         & -2.99   &     & 4.5382$\times10^{-2}$          & 3.27         \\
                  &512    & 1.22$\times10^{-2}$  & 4.7478$\times10^{-2}$       & 2.25    &    & 4.9493$\times10^{-2}$          & -0.13         \\
                  &1024    & 6.13$\times10^{-3}$ & 3.8797$\times10^{-3}$        & 3.61   &    & 2.4447$\times10^{-2}$          & 1.02         \\ 
                 &2048     & 3.06$\times10^{-3}$ & 2.6897$\times10^{-2}$        & -2.79    &    & 2.0692$\times10^{-3}$          & 3.56    \\ \hline
\multirow{7}{*}{$\mathcal{E}^\infty$}       &32         & 1.96$\times10^{-1}$ & 1.5846$\times10^{-1}$        & --  & \multirow{7}{*} {0.71, 0.48}         & 3.3871$\times10^{-1}$          & --      & \multirow{7}{*} {1.45, 0.99}      \\
                 &64     & 9.82$\times10^{-2}$  & 2.1296$\times10^{-2}$        & 2.88     &     & 1.7262$\times10^{-2}$          & 0.97           \\
                 &128     & 4.91$\times10^{-2}$  & 7.1544$\times10^{-3}$        & 1.57     &     & 5.1901$\times10^{-2}$          & 1.73           \\
                &256      & 2.45$\times10^{-2}$  & 5.0560$\times10^{-2}$        & -2.82    &     & 1.4598$\times10^{-2}$          & 1.83          \\
                 &512     & 1.22$\times10^{-2}$  & 1.2498$\times10^{-2}$        & 2.02    &    & 7.3909$\times10^{-3}$          & 0.98          \\
                &1024      & 6.13$\times10^{-3}$ & 1.2720$\times10^{-3}$       & 3.30    &      & 3.2262$\times10^{-3}$          & 1.20          \\ 
                &2048      & 3.06$\times10^{-3}$ & 8.6497$\times10^{-3}$        & -2.77   &     & 7.7599$\times10^{-4}$          & 2.06    \\ \hline
\end{tabular}
\end{table}

\begin{table}[]
\centering
\caption{\REVIEW{Error norm data for the hexagram case with spatially varying Robin boundary conditions using Approach C considered in Sec.~\ref{sec_spatially_varying_robin}}}
\label{tab_torus}
\begin{tabular}{cllllllll}
\hline
\multicolumn{1}{l}{}  & \multirow{2}{*}{N} & \multirow{2}{*}{h} & \multicolumn{3}{l}{Continuous indicator function}       & \multicolumn{3}{l}{Discontinuous indicator function}    \\ \cline{4-9} 
\multicolumn{1}{l}{}  &                    &                    & Error    & Order & Fit ($m,R^2$)        & Error    & Order & Fit ($m,R^2$)        \\ \hline
\multirow{7}{*}{$\mathcal{E}^1$} &32 & 1.96$\times10^{-1}$ & 6.4450e-1        & --      & \multirow{7}{*} {0.75, 0.96}      & 8.4947e-1        & --        & \multirow{7}{*} {1.00, 0.96}     \\
                  &64    & 9.82$\times10^{-2}$  & 3.8669$\times10^{-1}$        & 0.74     &   & 5.5409e-1            & 0.62         \\
                  &128    & 4.91$\times10^{-2}$  & 2.1738$\times10^{-1}$        & 0.83   &    & 3.7542$\times10^{-1}$          & 0.56         \\
                  &256    & 2.45$\times10^{-2}$  & 8.3523$\times10^{-2}$         & 1.38   &     & 8.0590$\times10^{-2}$          & 2.21         \\
                  &512    & 1.22$\times10^{-2}$  & 6.1208$\times10^{-2}$       & 0.45    &    & 9.1875$\times10^{-2}$          & -0.19         \\
                  &1024    & 6.13$\times10^{-3}$ & 4.2044$\times10^{-2}$        & 0.54   &    & 3.0769$\times10^{-2}$          & 1.58         \\ 
                 &2048     & 3.06$\times10^{-3}$ & 3.4124$\times10^{-2}$        & 0.30    &    & 1.3971$\times10^{-3}$          & 1.14    \\ \hline
\multirow{7}{*}{$\mathcal{E}^\infty$}       &32         & 1.96$\times10^{-1}$ & 1.6739$\times10^{-1}$        & --  & \multirow{7}{*} {0.71, 0.95}         & 1.8868$\times10^{-1}$          & --      & \multirow{7}{*} {0.72, 0.95}      \\
                 &64     & 9.82$\times10^{-2}$  & 1.2126$\times10^{-1}$        & 0.47     &     & 1.4155$\times10^{-1}$          & 0.41           \\
                 &128     & 4.91$\times10^{-2}$  & 1.0764$\times10^{-1}$        & 0.17     &     & 1.2363$\times10^{-1}$          & 0.20           \\
                &256      & 2.45$\times10^{-2}$  & 3.4357$\times10^{-2}$        & 1.65    &     & 4.2007$\times10^{-2}$          & 1.56          \\
                 &512     & 1.22$\times10^{-2}$  & 3.7526$\times10^{-2}$        & -0.13    &    & 4.2895$\times10^{-2}$          & -0.03         \\
                &1024      & 6.13$\times10^{-3}$ & 1.5288$\times10^{-2}$       & 1.30    &      & 1.8924$\times10^{-2}$          & 1.18          \\ 
                &2048      & 3.06$\times10^{-3}$ & 9.5637$\times10^{-3}$        & 0.68   &     & 9.5005$\times10^{-3}$          & 0.99    \\ \hline
\end{tabular}
\end{table}

\begin{table}[]
\centering
\caption{\REVIEW{Error norm data for the hexagram with smooth exterior corners case using Approach C with the discontinuous indicator function considered in Sec.~\ref{sec_rounded_hexagram}}}
\label{tab_torus}
\begin{tabular}{cllllllll}
\hline
\multicolumn{1}{l}{}  & \multirow{2}{*}{N} & \multirow{2}{*}{h} & \multicolumn{3}{l}{Neumann problem}       & \multicolumn{3}{l}{Robin problem}    \\ \cline{4-9} 
\multicolumn{1}{l}{}  &                    &                    & Error    & Order & Fit ($m,R^2$)        & Error    & Order & Fit ($m,R^2$)        \\ \hline
\multirow{7}{*}{$\mathcal{E}^1$} &32 & 1.96$\times10^{-1}$ & 7.6209$\times10^{-1}$        & --      & \multirow{7}{*} {1.08, 0.91}      & 3.9049e-1        & --        & \multirow{7}{*} {1.26, 0.89}     \\
                  &64    & 9.82$\times10^{-2}$  & 1.4221$\times10^{-1}$        & 2.42     &   & 1.0627e-1            & 1.88          \\
                  &128    & 4.91$\times10^{-2}$  & 1.4759$\times10^{-1}$        & -0.05   &    & 1.0857$\times10^{-1}$          & -0.03         \\
                  &256    & 2.45$\times10^{-2}$  & 2.2453$\times10^{-2}$         & 2.72   &     & 1.4303$\times10^{-2}$          & 2.92         \\
                  &512    & 1.22$\times10^{-2}$  & 1.4344$\times10^{-2}$       & 0.65    &    & 2.9291$\times10^{-3}$          & 2.29         \\
                  &1024    & 6.13$\times10^{-3}$ & 1.5886$\times10^{-2}$        & -0.15   &    & 4.7328$\times10^{-3}$          & -0.69         \\ 
                 &2048     & 3.06$\times10^{-3}$ & 6.7013$\times10^{-3}$        & 1.25    &    & 2.9825$\times10^{-3}$          & 0.67    \\ \hline
\multirow{7}{*}{$\mathcal{E}^\infty$}       &32         & 1.96$\times10^{-1}$ & 2.3668$\times10^{-1}$        & --  & \multirow{7}{*} {0.95, 0.95}         & 1.2099$\times10^{-1}$          & --      & \multirow{7}{*} {1.00, 0.88}      \\
                 &64     & 9.82$\times10^{-2}$  & 4.8712$\times10^{-2}$        & 2.28     &     & 4.0233$\times10^{-2}$          & 1.59           \\
                 &128     & 4.91$\times10^{-2}$  & 4.3138$\times10^{-2}$        & 0.18     &     & 4.9885$\times10^{-2}$          & -0.31           \\
                &256      & 2.45$\times10^{-2}$  & 2.5443$\times10^{-2}$        & 0.76    &     & 1.0488$\times10^{-2}$          & 2.25          \\
                 &512     & 1.22$\times10^{-2}$  & 7.8196$\times10^{-3}$        & 1.70    &    & 3.4723$\times10^{-3}$          & 1.59         \\
                &1024      & 6.13$\times10^{-3}$ & 6.8699$\times10^{-3}$       & 0.19    &      & 1.9522$\times10^{-3}$          & 0.83          \\ 
                &2048      & 3.06$\times10^{-3}$ & 3.3329$\times10^{-3}$        & 1.04   &     & 3.4365$\times10^{-3}$          & -0.82    \\ \hline
\end{tabular}
\end{table}

\begin{table}[]
\centering
\caption{\REVIEW{Error norm data for the circle case with spatially varying Neumann boundary conditions using Approach C and D considered in Appendix~\ref{sec_other_ext_results}. The continuous indicator function is used here. A convergence rate of $\mathcal{O}(h^{0.98})$ (respectively, $\mathcal{O}(h^{0.95})$) with an $R^2$ value of 0.95 (respectively, 0.98) in the $L^{\infty}$ (respectively, $L^1$) norm is obtained using Approach C. In the case of Approach D with the top hat kernel function, a convergence rate of $\mathcal{O}(h^{-0.02})$ (respectively, $\mathcal{O}(h^{-0.01})$) with an $R^2$ value of 0.19 (respectively, 0.03) in the $L^{\infty}$ (respectively, $L^1$) norm is obtained. With the spline kernel function, Approach D achieves a convergence rate of  $\mathcal{O}(h^{-0.01})$ (respectively, $\mathcal{O}(h^{0.02})$) with an $R^2$ value of 0.69 (respectively, 0.55) in the $L^{\infty}$ (respectively, $L^1$) norm.}}
\label{tab_other_ext_results}
\begin{tabular}{clllllll}
\hline
\multicolumn{1}{l}{}  & \multirow{2}{*}{N} & \multicolumn{2}{l}{Approach C}       & \multicolumn{2}{l}{Approach D (top hat function)} & \multicolumn{2}{l}{Approach D (spline function)} \\ \cline{3-8} 
\multicolumn{1}{l}{}  &                     & Error    & Order        & Error         & Order             & Error           & Order         \\ \hline
\multirow{7}{*}{$\mathcal{E}^1$} & 32                          & 5.0137$\times10^{-2}$ & --    & 5.9980$\times10^{-1}$      & --         & 1.3436        & --                           \\
                      &64                             & 1.6119$\times10^{-2}$ & 1.64                    & 5.6372$\times10^{-1}$      & 0.9                           & 1.2869        & 0.06                         \\
                      & 128                          & 1.3963$\times10^{-2}$ & 0.21                    & 6.4900$\times10^{-1}$      & -0.20                          & 1.2772        & 0.01                         \\
                      & 256                         & 4.6304$\times10^{-3}$  & 1.59                    & 6.1927$\times10^{-1}$      & 0.07                          & 1.2737        & 4.0$\times10^{-3}$                         \\
                      & 512                         & 3.4634$\times10^{-3}$ &0.42                    & 6.0978$\times10^{-1}$     & 0.02                            & 1.2729       & 9.1$\times10^{-4}$                         \\
                      &  1024                      & 1.5665$\times10^{-3}$ & 1.14                    & 6.0404$\times10^{-1}$      & 0.01                            & 1.2720        & 1.1$\times10^{-3}$                         \\
                      &2048                        & 8.1497$\times10^{-4}$ & 0.94                    & 6.0404$\times10^{-1}$      & 0.00                           & 1.2714       & 6.5$\times10^{-4}$                         \\ \hline
\multirow{7}{*}{$\mathcal{E}^\infty$}               &32                        & 8.4450$\times10^{-3}$ & --     & 9.0222$\times10^{-2}$      & --             & 2.1357$\times10^{-1}$        & --                            \\
                      & 64                       & 3.7481$\times10^{-3}$  & 1.17                      & 1.0070$\times10^{-1}$      & -0.16                              & 2.1124$\times10^{-1}$        & 0.02                           \\
                      & 128                     & 5.0424$\times10^{-3}$ & -0.43                      & 1.1195$\times10^{-1}$      & -0.15                              & 2.1110$\times10^{-1}$        & 0.01                           \\
                      & 256                     & 1.0171$\times10^{-3}$ & 2.31                      & 1.0500$\times10^{-1}$      & 0.09                              & 2.1627$\times10^{-1}$        & 9.6$\times10^{-4}$                          \\
                      &  512                    & 6.5633$\times10^{-4}$ & 0.63                      & 1.0480$\times10^{-1}$      & 2.7$\times10^{-3}$        & 2.1855$\times10^{-1}$        & -0.03                     \\
                      &  1024                  & 3.1095$\times10^{-4}$ & 1.08                      & 1.0284$\times10^{-1}$      & 0.03                              & 2.1849$\times10^{-1}$        & -0.02                          \\
                      &  2048                  & 1.5776$\times10^{-4}$ & 0.98                    & 1.0284$\times10^{-1}$      & 0.00                              & 2.1842$\times10^{-1}$         & 4.6$\times10^{-4}$                           \\ \hline
\end{tabular}
\end{table}